\documentclass[12pt]{amsart}
\usepackage{amsmath}
\usepackage{amsthm}
\usepackage{amsfonts} 
\usepackage{hyperref}
\usepackage{epsfig}

\numberwithin{equation}{section}

\allowdisplaybreaks

\newtheorem{theorem}{Theorem}[section]

\newtheorem{proposition}[theorem]{Proposition}
\newtheorem{lemma}[theorem]{Lemma}
\newtheorem{remark}[theorem]{Remark}
\newtheorem{definition}[theorem]{Definition}

\renewcommand{\epsilon}{\varepsilon}

\newcommand{\abs}[1]{\left\vert #1\right\vert}
\newcommand{\1}[1]{{\mathbf 1}{\{#1\}}}
\newcommand{\R}{\mathbb{R}}
\newcommand{\N}{\mathbb{N}}
\newcommand{\Z}{\mathbb{Z}}
 
\newcommand{\PR}{\mathbb{P}_p}
\newcommand{\ES}{\mathbb{E}}

\newcommand{\ins}{{\rm INS}}

\newcommand{\extg}{\partial_{\rm ext}}
\newcommand{\extk}{\partial_{{\rm ext},K}}

\newcommand{\backtrack}{\mathcal{B}\mathcal{K}}
\newcommand{\Fsurf}{\mathcal{S}}
\newcommand{\dtrap}{\mathcal{T}}

\newcommand{\dist}{\vert\vert}
\newcommand{\xlower}{x_{{\rm low}}}
\newcommand{\bigxlower}{X_{{\rm low}}}
\newcommand{\const}{Q}

\newcommand{\hatxi}{\zeta}

\newcommand{\FLOOR}{{\rm Floor}}
\newcommand{\LID}{{\rm Lid}}

\title[The speed of biased walk in supercritical percolation]{Phase transition for the speed of the biased random walk on the supercritical percolation cluster
}
\date{}
\author[A.~Fribergh]{Alexander FRIBERGH}
\address{ Courant Institute of Mathematical Sciences,
                  251 Mercer Street,
                  New York University,
                  New York, 10012-1185, U.S.A.} 
\email{fribergh@cims.nyu.edu}

\author[A.~Hammond]{Alan HAMMOND
}
\address{ Department of Statistics, 
           University of Oxford,
                  1 South Parks Road,
                  Oxford, OX1 3TG, U.K. } 
\email{hammond@stats.ox.ac.uk}

\thanks{A.H. was supported in part by NSF
grants DMS-0806180 and OISE-0730136 and by EPSRC grant EP/I004378/1.} 

\begin{document}

\maketitle

\begin{abstract}
We prove the sharpness of the phase transition for the speed in biased random walk on the 
supercritical percolation cluster on $\Z^d$. 
That is, for each $d \geq 2$, and for any supercritical parameter $p > p_c$, we prove the existence of a critical strength for the bias, such that, below this value, the speed is positive, and, above the value, it is zero. We identify the value of the critical bias explicitly,
and, in the sub-ballistic regime, we find the polynomial order of the distance moved by the particle.
Each of these conclusions is obtained by investigating
the geometry of the traps that are most effective at delaying the walk. 

A key element in proving our results is to understand that, on large scales, the particle trajectory is essentially one-dimensional; we prove such a {\it dynamic renormalization} statement in a much stronger form than was previously known.
\end{abstract}
\section{Introduction}

A fundamental topic in the field of random walk in random environment is the long-term displacement of the particle, and how it is determined by fluctuations in the random landscape. Homogenization is a powerful technique that is often used to provide classical results, such as  laws of large numbers and central limit theorems; however, many natural examples exhibit anomalous behaviors.
A fundamental cause of such behavior is trapping. This phenomenon is present in many physical systems and is believed to be associated with certain universal limiting behaviors. (See \cite{Zaslavsky} for a physical discussion of anomalous transport arising in the presence of trapping effects.) In an attempt to find a random system  whose  scaling limit may be shared by many trapping models, Bouchaud proposed a simple idealized model. In {\em Bouchaud's trap model}, a random walk jumps on the vertices of a large graph,  and experiences  delays  whose law is determined by independent and identically distributed heavy-tail random variables attached to the vertices of the graph in question.
The lectures \cite{bacone} summarise the main results known on Bouchaud's trap model in $\Z^d$. In the case $d\geq 2$, certain common limiting behaviors appear, such as convergence of the suitably-rescaled ``internal clock" of the walker to the inverse of a stable subordinator, and the convergence of the rescaled process to
Brownian motion time-changed by this limiting clock process (a dynamics known as {\it fractional kinetics}).
As such, the limiting process manifests the phenomenon of {\it aging}, remaining static for random periods whose duration is comparable to the system age.  
 Several proofs presented in \cite{bacone} may be applied to more general models, which led the authors to suggest that these limiting behaviors should be a common feature of many trap models. The very specific 
one-dimensional case in  Bouchaud's trap model is studied in \cite{FIN}, which constructs its scaling limit, under which Brownian motion is time-changed by a spatially dependent singular process. 

An important challenge is to obtain derivations of scaling limit processes similar to fractional kinetics, and thus results on aging, in physically realistic settings. The article~\cite{Bouchaud} provides a physical overview in this regard. For the random energy model,~\cite{bbgone} and~\cite{bbgtwo} make such a derivation in the case of a caricature of Glauber dynamics (see also~\cite{BC} for a review); for the Sherrington-Kirkpatrick and $p$-spin models with ``random hopping time" dynamics,~\cite{bagun} identifies the mechanism of finding deepest traps on subexponential time-scales, while~\cite{babc} proves fractionally kinetic behavior on short exponential time-scales.

A fertile arena in which to seek physically realistic models for trapping and aging is the field of random walks in disordered media. 
Random walks in random environment (see~\cite{Zeitouni} for a review) provide some of the most interesting candidates for this inquiry. 
One regime in particular seems to be of great interest,  namely, when the walk is directionally transient and has zero speed.
Indeed,  \cite{kks} obtained stable scaling limits for some such models in one dimension,  and \cite{esz} refines this analysis to prove aging. These behaviors are characteristic of trapping, and these works raise the prospect of broadening this inquiry to the more physical and more demanding multi-dimensional models in this regime.

The last decade has witnessed a resurgence of mathematical interest in multi-dimensional
random walk in random environment, a subject that is deeply rooted in physics.
 As an illustrative model of diffusion in disordered media,
De Gennes~\cite{deGennes} proposed studying the random walk on the supercritical percolation cluster.
By now, the isotropic (or simple) random walk in this environment is very well understood, with progress including the derivation of a functional central limit theorem (see~\cite{BB} and~\cite{MP}) and a strong heat-kernel estimate~\cite{barlow}. 
However, the model does not exhibit the characteristic features of trapping. In contrast, when a constant external field is applied to a randomly evolving particle on the supercritical percolation cluster, physicists have long predicted anomalous behavior: see~\cite{BD},~\cite{Dhar} and~\cite{DS}. 
We will refer to this anisotropic model as the biased walk on the supercritical percolation cluster, leaving the formal definition for later. 
 
After a long time, two concurrent mathematical treatments,~\cite{BGP} and~\cite{Sznitman}, resolved some aspects of the physicists' predictions. 
Each of these works showed that, when the applied field is weak, the particle moves at positive speed, while, if the field is very strong, then the speed vanishes. This seemingly paradoxical effect is best understood in terms of trapping.
Indeed, as the walker is pushed by the applied field into new territory, it may sometimes encounter
 dead-ends that act as traps which detain the walker for long periods. This delaying mechanism becomes more powerful as a stronger field is applied,  eventually resulting in a sub-ballistic motion. Indeed, A-S. Sznitman mentions in the survey \cite{sznitmannew} that this mechanism appears to be similar to that responsible for aging in Bouchaud's trap model. 

The inference made in~\cite{BGP} and~\cite{Sznitman}  left unproved one of the key predictions, namely that the transition from the ballistic to the sub-ballistic regime is sharp. 
That is, for each $d \geq 2$ and $p > p_c$, there exists a critical value $\lambda_c$ such that, for values of the constant external field strictly less than $\lambda_c$, the walk has positive speed, while, for values strictly exceeding $\lambda_c$, the walk has zero speed.  This conjecture was also made by~\cite{BGP}. 
Evidence for the conjecture  is provided by its validity in the simpler context of biased random walk on a supercritical Galton-Watson tree~\cite{LPP}. In this case, a  detailed investigation of the trapping phenomenon has already been undertaken by \cite{FG}, \cite{BH} and \cite{H}, including
 the analysis of the particle's scaling limits in the sub-ballistic regime.

The task of establishing the phase transition is intimately related to understanding in greater detail the trapping phenomenon manifested by the biased walk in the supercritical percolation cluster: that is, to finding the stochastic geometry of the traps that detain the walker at late time, and to determining the sub-ballistic growth of walker displacement in the zero-speed regime. The present article is devoted to proving the conjecture and to elucidating trap geometry.
In undertaking this study, we develop a new technique for proving results on dynamic renormalization (as we will shortly describe).  This technique  forms a 
foundation for a deeper analysis of the trapping phenomenon in biased random walk on the supercritical percolation cluster, and also for other reversible disordered systems (see~\cite{FR} for such an application).

\subsection{Model definition.}

We now formally define the biased walk on the supercritical percolation cluster, using the definition in \cite{Sznitman}.

Firstly, we describe the environment. 
We denote by $E(\Z^d)$ the edges of the nearest-neighbor lattice $\Z^d$ for some $d\geq 2$. We fix $p \in (0,1)$ and perform a Bernoulli bond-percolation by picking a random configuration $\omega \in \Omega:=\{0,1\}^{E(\Z^d)}$ where each edge $e$ has probability $p$ of verifying $\omega(e)=1$, independently of the assignations made to all the other edges. We introduce the corresponding measure 
\[
P_{p}= (p \delta_1 + (1-p) \delta_0)^{\otimes E(\Z^d)}.
\]

An edge $e$ will be called open in the configuration $\omega$ if $\omega(e)=1$. The remaining edges will be called closed.  This naturally induces a subgraph of $\Z^d$ which will also be denoted by $\omega$, and it  yields a partition of $\Z^d$ into open clusters and isolated vertices. 

It is classical in percolation that for $p>p_c(d)$, where $p_c(d)\in (0,1)$ denotes the critical percolation probability of $\Z^d$ (see~\cite{Grimmett} Theorem~1.10 and Theorem~8.1), there exists a unique infinite open cluster $K_{\infty}(\omega)$, $P_p$ almost surely. Moreover, for each $x \in \Z^d$, 
the following event has positive $P_p$-probability: 
\[
\mathcal{I}_x= \Big\{ \text{there is a unique infinite cluster $K_{\infty}(\omega)$ and it contains $x$} \Big\}.
\]
We will use the shorthand $\mathcal{I}=\mathcal{I}_0$. We further define
\[
{\bf P}_p[~\cdot~] = P_{p}[~\cdot\mid \mathcal{I}].
\]

The bias $\ell=\lambda \vec \ell$ depends on two parameters: the strength $\lambda >0$, and the bias direction $\vec \ell \in S^{d-1}$ which lies in the unit sphere with respect to the Euclidean metric of $\R^d$. Given a configuration $\omega \in \Omega$, we consider the reversible Markov chain $X_n$ on $\Z^d$ with law $P_x^{\omega}$,  whose transition probabilities $p^{\omega}(x,y)$ for $x,y\in \Z^d$ are defined by
\begin{enumerate}
\item $X_0=x$, $P_x^{\omega}$-a.s.,
\item $p^{\omega}(x,x)=1$, if $x$ has no neighbor in $\omega$, and
\item $\displaystyle{p^{\omega}(x,y) =\frac{c^{\omega}(x,y)}{\sum_{z \sim x} c^{\omega}(x,z)}}$,
\end{enumerate}
where $x\sim y$ means that $x$ and $y$ are adjacent in $\Z^d$.  Here we set
\[
\text{for all $x,y \in \Z^d$,} \qquad \ c^{\omega}(x,y)=\begin{cases}
                                            e^{(y+x)\cdot\ell} & \text{ if } x\sim y \text{ and } \omega(\{x,y\})=1, \\
                                            0           & \text{ otherwise.}\end{cases} 
\]                                           

In the special case $x=0$, we use $P_0^{\omega}=P^{\omega}$.

We denote by $c^{\omega}(e)=c^{\omega}(x,y)$ the conductance of the edge $e=[x,y]$ in the configuration $\omega$, a notation which is natural in light of the relation between reversible Markov chains and electrical networks,
for a presentation of which the reader may consult
 ~\cite{DoyleSnell} or ~\cite{LP}. The above Markov chain is reversible with respect to the  invariant measure given by 
\[
\pi^{\omega}(x)=\sum_{y \sim x} c^{\omega}(x,y).
\]  

Finally, we define the annealed law of the biased random walk on the infinite percolation cluster by the semi-direct product $\PR = {\bf P}_p  \times  P_0^{\omega}[\,\cdot\,]$.

\subsection{Results}

 In order to state our results, let us recall some results proved in~\cite{Sznitman}. In all cases, i.e., for $d\geq2$ and $p > p_c$, the walk is directionally transient in the direction $\vec{\ell}$ (see Theorem 1.2 in~\cite{Sznitman}):
 \begin{equation}
 \label{def_dir_trans}
 \lim_{n\to \infty} X_n \cdot \vec{\ell} =\infty, \qquad \PR-\text{almost surely,}
 \end{equation}
 and verifies  the law of large numbers  (see Theorem 3.4 in~\cite{Sznitman}),
\[
\lim_{n\to \infty} \frac{X_n} n =v , \qquad \PR-\text{almost surely},
\]
where $v\in \R^d$ is a constant vector.
 
We now introduce the backtrack function of $x \in \Z^d$, which will be fundamental for gauging the extent of the slowdown effect on the walk: 
\begin{equation}
\label{def_backtrack}
\mathcal{B}\mathcal{K}(x)=\begin{cases} 0 & \text{ if } x\notin  K_{\infty}\\
                                              \displaystyle{\min_{(p_x(i))_{i\geq 0}\in \mathcal{P}_x}}\max_{i\geq 0} (x-p_x(i))\cdot \vec{\ell}  & \text{otherwise,}
\end{cases}
\end{equation}
where $\mathcal{P}_x$ is the set of all infinite open vertex-self-avoiding paths starting at~$x$. As Figure~$1$ indicates, connected regions where $\backtrack$ is positive may be considered to be traps: indeed, from points in such regions, the walk must move in the direction opposed to the bias in order to escape the region. The height of a trap may be considered to be the maximal value of $\backtrack (x)$ attained by the vertices $x$ in the trap.

\begin{figure}\label{defBK}
\centering\epsfig{file=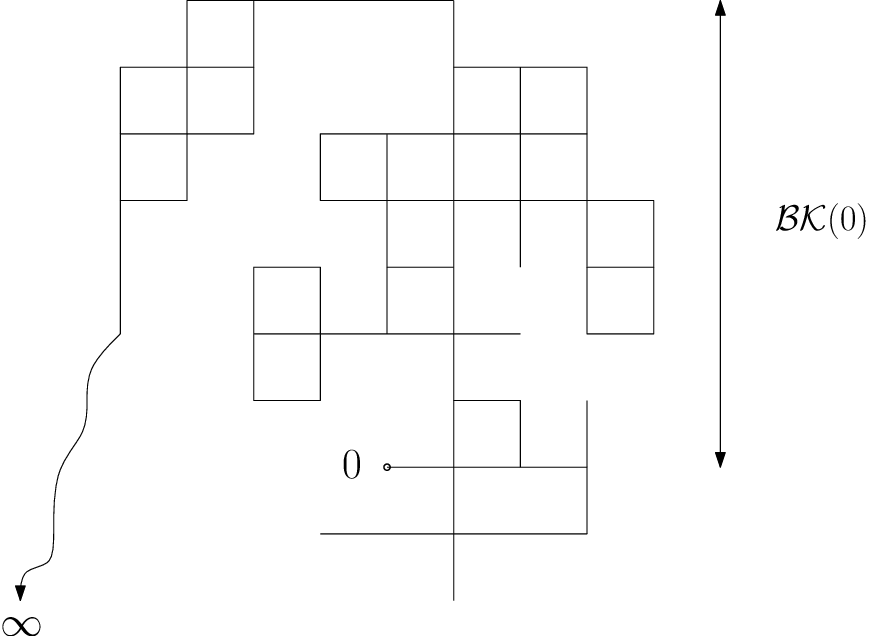, width=8cm}
\caption{Illustrating the value of $\mathcal{B}\mathcal{K}(0)$ in the case of an axial bias $\vec \ell$.}
\end{figure}

The next proposition is our principal percolation estimate. Section~\ref{sect_perco} 
is devoted to its proof.
\begin{proposition}
\label{backtrack_exponent}
For $d\geq 2$, $p > p_c$ and $\vec{\ell} \in S^{d-1}$, there exists $\hatxi(p,\vec{\ell},d)\in (0,\infty)$ such that 
\[
\lim_{n \to \infty}  n^{-1} \log {\bf P}_p[ \mathcal{B}\mathcal{K}(0)> n  ] = -\hatxi(p,\vec{\ell},d).
\]
\end{proposition}

We may then define the exponent $\gamma \in (0,\infty)$:
\begin{equation}
\label{def_gamma}
\gamma= \frac {\hatxi}{2\lambda}.
\end{equation}

Our first theorem confirms the conjecture of sharp transition from the ballistic to the sub-ballistic regime,
and identifies explicitly the critical value~$\lambda_c$ of the bias.
\begin{theorem}
\label{the_theorem}
For $d\geq 2$ and $p > p_c$, set $\lambda_c=\hatxi(p,\vec{\ell},d)/2$. We have that
\begin{enumerate}
\item if $\lambda <\lambda_c$, or, equivalently, $\gamma>1$, then $v\cdot \vec{\ell} >0$,
\item  if $\lambda >\lambda_c$, (or $\gamma<1$), then $v=\vec{0}$.
\end{enumerate}
\end{theorem}

\begin{remark}
We conjecture that for $\gamma=1$, we have that $v=0$. Slightly stronger percolation estimates are needed for this result. We feel that proving these estimates would be a distraction from our aims in the present work. It would also of interest to investigate the velocity direction $v/{\dist v \dist}$ as a function of bias direction $\vec{\ell} \in S^{d-1}$ when $\gamma > 1$. 
\end{remark}

We identify the regime where a functional central limit theorem holds:
\begin{theorem}
\label{the_theorem2}
Let $d \geq 2$ and $p > p_c$.  
Assume that $\gamma >2$.

 The $D(\R_+,\R^d)$-valued processes $B_{\cdot}^n = \frac 1 {\sqrt n}(X_{[\cdot n]}-[\cdot n] v)$ converge under $\PR$ to a Brownian motion with non-degenerate covariance matrix, where $D(\R_+,\R^d)$ denotes the space of right continuous $\R^d$-valued functions on $\R_+$ with left limits, endowed with the Skorohod topology, c.f.~Chapter 3 of \cite{EK}.
\end{theorem} 

G\'erard Ben Arous has pointed out to us that, in the part of the ballistic phase given by $\gamma \in (1,2)$, it is natural to expect fluctuations to be dominated by trap sojourns and thus to have anomalous rather than Gaussian behavior. 

Thirdly, in the sub-ballistic regime, we find the polynomial order of the magnitude of the walk's displacement.
\begin{theorem}
\label{the_theorem3}
Set $\Delta_n = \inf \big\{ m \in \N: X_m \cdot \vec{\ell} \geq n \big\}$. 

Let $d\geq 2$ and $p > p_c$.  
If $\gamma \leq 1$ then
\[
\lim \frac{\log \Delta_n}{\log n} = 1/\gamma, \qquad \PR-\text{almost surely},
\]
and
\[
\lim \frac{\log X_n \cdot \vec{\ell}}{\log n} =\gamma, \qquad \PR-\text{almost surely}.
\]
\end{theorem}

In dimension two, the critical bias has a rather explicit expression.
\begin{remark}\label{lkjh}
In the case $d=2$ (with $p > p_c = 1/2$), we have that 
\[
\hatxi(p,\vec\ell,2)  = 2 \inf_{\vec\ell' \in S^1} \frac{\xi^{(\vec\ell')}}{\vec\ell' \cdot \vec\ell}\,.
\]
Here, for $\vec\ell' \in S^1$,  
$\xi^{(\vec\ell')}$ 
denotes the $\ell'$-direction inverse correlation length in the subcritical percolation $P_{1-p}$, given by 
$\xi^{(\vec\ell')} = - \lim_{n \to \infty} \frac{\log P_{1-p}\big( 0 \leftrightarrow \lfloor n \ell' \rfloor \big)}{ n}$, where $\lfloor \cdot \rfloor$ denotes the component-wise integer part.
\end{remark}
\begin{remark}
In the case considered in~\cite{BGP}, where the conductances are $c^{\omega}(x,y)=\beta^{\max (x\cdot e_1,y\cdot e_1)} \1{\omega(\{x,y\})=1}$ with $\beta>1$, we may obtain the same results with $\gamma=\hatxi/\log \beta$ by the same methods. This was conjectured in~\cite{BGP}.
\end{remark}

\subsection{Trap geometry} 
In the zero-speed regime $\lambda >\lambda_c$ identified in Theorem~\ref{the_theorem}, the particle at late time typically resides near the base of a large trap. Our investigation of the sharp transition prompts us to examine the typical geometry of this trap. 

In the two-dimensional case, a trap is surrounded by a path of open dual edges, while, in higher dimensions, this surrounding is a surface comprised of open dual plaquettes. We refer to this as the trap surface. In either case, the surface exists in the subcritical phase, and is costly to form.

Our study reveals a significant difference in trap geometry between the cases $d =2$ and $d \geq 3$: see Figure~$2$. In any dimension, the typical trap is a long thin object oriented in some given direction. When $d\geq 3$, the trap surface is typically uniformly narrow, in the sense that it has a bounded intersection with most hyperplanes orthogonal to the trap's direction.  
On the other hand, in the two-dimensional case, the trap surface may be viewed as two disjoint subcritical dual paths that meet just below the base of the trap. These two dual paths each travel in the trap direction a distance given by the trap height. Rather than remaining typically at a bounded distance from one another, the two paths separate. This means that, in contrast to the higher dimensional case, the trap surface is not uniformly narrow.

This distinction means that, in the case that $d \geq 3$, there is a \lq\lq trap entrance\rq\rq, which is a single vertex, located very near the top of the trap, through which the walk must pass in order to fall into the trap. In contrast, in dimension $d=2$, the top of the trap is comprised of a line segment orthogonal to the trap direction and of length of the order of the square-root of the trap height. 

Finally, we wish to emphasize that for generic non-axial $\vec{\ell}$, the trap direction does not coincide with $\vec{\ell}$. In the two-dimensional case, it is given by the minimizer $\vec{\ell '}$ in the infimum appearing in Remark~\ref{lkjh}. It is a simple inference from the theory of subcritical percolation connections (see~\cite{civ}) that this minimizer is unique. We do not need this fact so we will not prove it. We also mention that the factor of two appearing in Remark~\ref{lkjh} arises from the two disjoint subcritical dual paths that form the trap.

\begin{figure}[h]
\centering
\epsfig{file=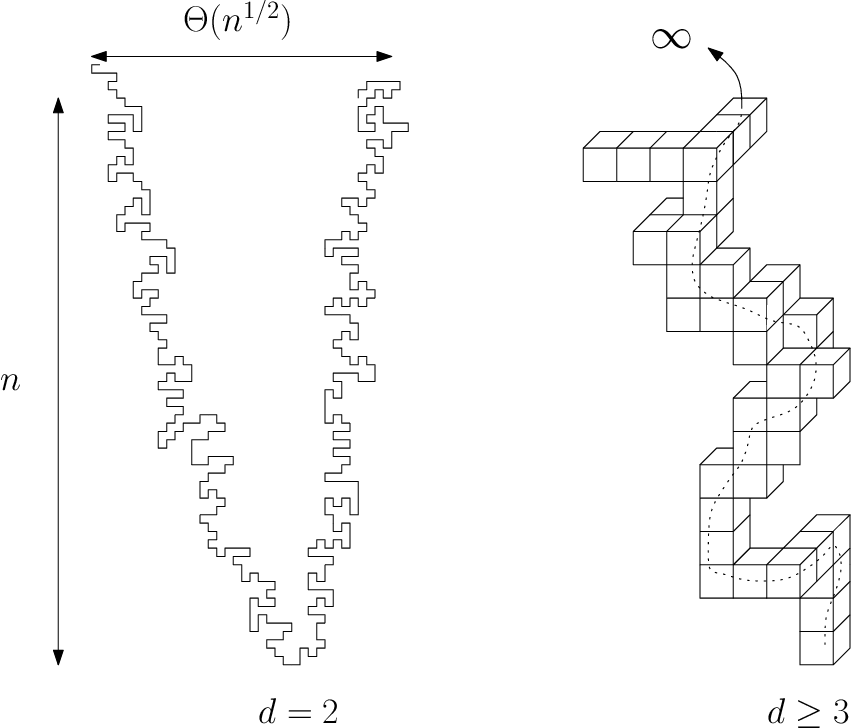, width=13cm}
\caption{The dual trap surfaces associated with typical traps.}
\end{figure}

\subsection{Heuristic argument for the value of the critical bias.} \label{heuristic}
We now present a heuristic argument that makes plausible Theorems ~\ref{the_theorem},~\ref{the_theorem2} and~\ref{the_theorem3} and highlights the trapping mechanism.

Suppose that the particle has fallen into a trap and is now at  its base, (a vertex that we will label $b$). How long will it take to return to the top of the trap? To answer this question, let us assume, for the sake of argument, that the top of the trap is composed of a unique vertex $t$. (As we mentioned previously, this assumption is essentially accurate for $d\geq 3$. The assumption is only a matter of convenience, and the conclusions that we will reach are also valid for $d=2$.)

Denoting by $p(x,y)$ the probability of reaching $y$ starting from $x$ without return to $x$, a simple reversibility argument yields 
\[
\pi(b) p(b,t)= \pi(t) p(t,b).
\]

Typical traps being long and thin, the term $p(t,b)$, which is the probability of reaching the base of the trap starting from the top, is typically bounded away from 0 uniformly. We now notice that, for a trap of height $h$, $\pi(t)/\pi(b)$ is equal to $\exp(-2\lambda h)$ up to a multiplicative constant, and thus
\[
p(b,t)\approx \exp(-2\lambda h).
\]

Hence, the typical time $T$ to reach the top of a trap of height $h$ is of the order of $\exp(2\lambda h)$. Recalling Proposition~\ref{backtrack_exponent}, the time to exit a randomly chosen trap is given by
\[
\PR[T\geq n]\approx \PR[\exp(2\lambda \backtrack (b))> n] \approx n^{-\gamma}.
\]

Actually, we will need a closely related statement, giving an upper bound on the tail of the time spent in a box.
\begin{proposition}
\label{BL2_PF}

For every $\epsilon>0$ and $\alpha>0$, for all $t$ large enough and for $L:[0,\infty) \to (0,\infty)$ satisfying $\lim \frac{\log L(t)} t =0$,  we have that
\[
\PR[T^{\text{ex}}_{B(L,L^{\alpha})}\geq t]\leq  CL^{d\alpha(2\gamma+1)}t^{-(1-\epsilon)\gamma}+e^{-cL},
\]
for some positive constants $c$ and $C$.
\end{proposition}

Now, if we were to model $\Delta_n$ as the total time spent in a sequence of $n$ independently selected traps, we would reach the conclusion of Theorem~\ref{the_theorem}. This heuristic appears to be accurate since it is natural to expect the trajectory of the particle to be linear on a macroscopic scale. However, there are substantial difficulties in justifying it:
\begin{enumerate}
\item there might be large lateral displacements orthogonal to the preferred direction. This may have two problematic effects: the walk might encounter many more than $n$ traps before $\Delta_n$; moreover, during such displacements the walk might spent a lot of time outside of traps and be significantly slowed down;
\item there might be strong correlations between traps;
\item finally, a walk that, having fallen to the base of a trap, has returned to its top, might be prevented from escaping from the trap, by making a large number of lengthy excursions into the trap. 
\end{enumerate}

Resolving these three potential problems is a substantial technical challenge, one that lies at the heart of this paper. All three difficulties will be overcome by using strong estimates on regeneration times obtained through a {\em dynamic renormalization} result which we now discuss.

\subsection{Dynamic renormalization} 
The key tool for our analysis in this work is to show that the particle will exit a large box in the preferred direction with overwhelming probability.
In order to state our conclusion to this effect,  
we need a little notation. 

\begin{definition}\label{defone}
We set $f_1:=\vec{\ell}$ and then choose further vectors so that $(f_i)_{1\leq i \leq d}$ forms an orthonormal basis of $\R^d$.

Given a set $V$ of vertices of $\Z^d$ we denote by $\abs{V}$ its cardinality, and by $E(V)=\{ [x,y] \in E(\Z^d): ~x,y\in V\}$ its set of edges. We set
\[
\partial V= \{x \notin V : y\in  V,~x\sim y\}.
\]
For any $L,L' \geq 1$, we set
\[
B(L,L')=\Bigl\{z\in \Z^d : \vert z \cdot \vec{\ell} \, \vert < L,\ \abs{z\cdot f_i} < L' \text{ for $i \in [2,d]$} \Bigr\},
\]
and
\[
\partial^+ B(L,L')=\{z\in \partial B(L,L') : z \cdot \vec{\ell} \geq  L \}.
\]
\end{definition}
\begin{theorem}
\label{BL}
For any $d\geq 2$ and $p>p_c$ and $\ell \in \R^d$. For any $\alpha>d+3$, there exist  constants $C > c > 0$ such that, for each $L \in \N$,
 \begin{equation}\label{eqbl}
  \PR[T_{\partial B(L,L^{\alpha})}\neq T_{\partial^+ B(L,L^{\alpha})}]  \leq C e^{-cL}.
  \end{equation}
\end{theorem}
Theorem~\ref{BL} resolves the three substantial difficulties to which we have alluded because, alongside an analysis using regeneration times, it proves that the walk moves quickly in relatively open space, and that its long-term trajectory has a uni-dimensional character.

The proof of Theorem~\ref{BL} is a rather intuitive one that can be extended to other reversible models. As already mentioned, similar methods can be used in the context of biased random walks in independently and identically distributed non-zero random conductances, see~\cite{FR}. From these methods follow the counterpart of Theorem~\ref{BL} in this model, a result which is one of the main pieces needed to characterize the regime where the asymptotic speed is positive.

In the course of the proof of Theorem~\ref{BL}, we introduce a natural decomposition of the percolation configuration into small components in which the walk is liable to move atypically, and a large complementary region that is, in effect, free space. This partition is reminiscent of the Harris decomposition of a supercritical Galton-Watson tree, under which the tree is partitioned into an infinite leafless subtree, and a collection of subcritical trees that hang off this subtree. In the context of biased random walk, the infinite tree may be considered to be a backbone, on which the walk advances at linear speed, while the subcritical trees are the traps that may delay the walk. The papers~\cite{FG}, ~\cite{BH} and ~\cite{H} undertake a detailed examination of trapping for biased random walks on such trees. While involved, these works depend in an essential way on the unambiguous backbone-trap decomposition available for such environments. One of the big challenges in the case of supercritical percolation is to find a useful analogue of this decomposition. The partition presented in the proof of Theorem~\ref{BL} provides a very natural candidate for this decomposition. It is our aim to pursue an inquiry into scaling limits of particle trajectory 
for biased walk in supercritical percolation; for this,
a fundamental role will be played both by the strong control on particle motion stated in Theorem~\ref{BL}, and by the tractable backbone-trap decomposition that is central to the proof of this result.  

A-S.~Sznitman has previously introduced two criteria, conditions $(T)$ and $(T^\gamma)$, which are closely related to Theorem~\ref{BL}. These conditions essentially correspond to~(\ref{eqbl}), replacing $B(L,L^{\alpha})$ by $B(L,CL)$ for some $C$ (and the right-hand side by $\exp(-L^{\gamma})$ in the case of $(T^\gamma)$); see~\cite{SZ2} for a formal definition. Sznitman was motivated to introduce these conditions in light of a celebrated conjecture that any uniformly elliptic random walk in random environment that is directionally transient has positive speed, (whether or not the walk is reversible). He succeeded in showing that any uniformly elliptic RWRE that satisfies condition $(T)$ has positive speed, and it is conjectured that the converse is true. This reinforces the sentiment that Theorem~\ref{BL} is key to the analysis of the asymptotic speed.

\subsection{Applications of the techniques and related open problems}
\subsubsection{Scaling limits} The sharp phase transition and sub-ballistic velocity exponent having been identified, a natural problem raised by the progress of this article is to understand the limiting behavior of the scaled displacement $n^{-\gamma}\vert X_n \vert$, an inquiry whose analogue for supercritical Galton-Watson trees has been completed. The problem on trees has an interesting structure, with the constant bias walk having no scaling limit due to a logarithmic periodicity effect \cite{FG}, but with certain randomly biased walks having stable limit laws \cite{H}. It appears that these two behaviors may coexist in the biased walk on supercritical percolation, the latter arising as a rational resonance effect; see Section 1.10 of~\cite{H} for further discussion of this point. In such a coexistence, fractionally kinetic behaviour might be expected to arise for ``irrational'' choices of the bias direction, while a persistent effect of periodicity would disrupt this behaviour for other choices; this dichotomy indicates how anisotropic motion in discrete disordered systems may furnish a rich array of trapping behaviours.

In regard to pursuing the analysis for trees in the more physically natural setting of percolation, we anticipate that several techniques introduced by this article will play a critical role. Notably, the dynamic renormalization Theorem~\ref{BL} will play such a role, 
since it provides a very strong control of the walk away from trap structures. Furthermore,  understanding of trap geometry has deepened, especially for $d\geq 3$, where a regeneration structure inside the trap is now apparent: see Remark~\ref{decompo_regen}.

\subsubsection{Speed as a function of bias} 
Perhaps the most important conceptual distinction between backbone and traps is the notion that increasing the bias strength should lead to an increase in particle velocity on the backbone but also to an increase in the delay caused by visits to traps. 
\begin{figure}[h]\label{speedconj}
\centering
\epsfig{file=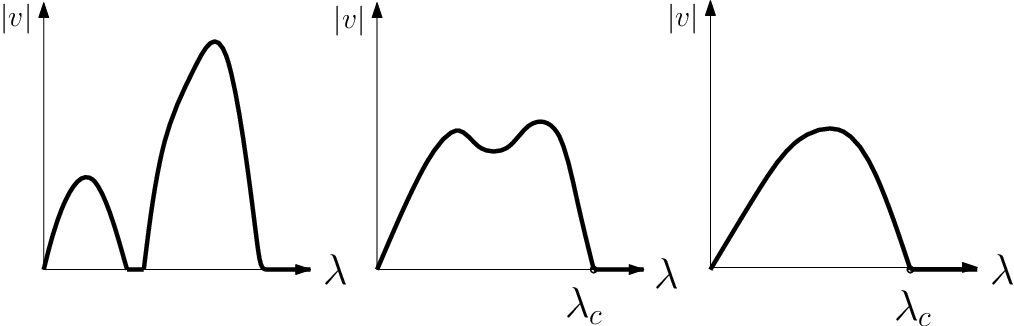, width=12cm}
\caption{The speed as a function of the bias. Sznitman and Berger, Gantert and Peres established positive speed at low $\lambda$, but their works left open the possibility depicted in the first sketch. Our work rules this out, though the behavior of the speed in the ballistic regime depicted in the second sketch remains possible. The third sketch shows the unimodal form predicted physically.}
\end{figure}
It is an attractive aim  then to establish the unimodal nature of the speed function depicted in the third graph in Figure~3. Theorem~\ref{the_theorem} establishes the existence and value of a critical point, 
while the technique of proof sheds much light on trap geometry,
but this does not mean that the conjecture is necessarily close at hand. Indeed, even the problem for trees is unresolved. The backbone for the environment considered in \cite{FG} is a supercritical Galton-Watson tree without leaves. For this environment, we recall an open problem from \cite{LPart}, to show that particle speed is an increasing function of the bias parameter. This has only been proved for large enough biases, see~\cite{BAFS}. It is very interesting to notice that an explicit formula for the speed has recently been obtained in~\cite{Aidekon} in the case for biased random walks on supercritical Galton-Watson trees (even with leaves). Nonetheless, at this point, it does not allow to solve the open problem from \cite{LPart}.

We mention also that the derivative at $\lambda = 1$ is predicted to exist, and indeed is known physically as the mobility. In fact, the Einstein relation links biased walk on the supercritical percolation cluster to its unbiased counterpart, because it
predicts that the mobility is proportional to diffusivity, which is the variance of the Brownian motion appearing in the central limit theorem for simple random walk on the supercritical percolation cluster. We refer the reader to~\cite{GMP} and \cite{LR} for some mathematical results on the Einstein relation in different models.

\subsubsection{The critical bias as a function of the percolation parameter} 
Fixing $\ell \in S^{d-1}$, we may view the critical bias $\lambda_c:[p_c,1] \to (0,\infty)$ as a function of the supercritical percolation parameter. Traps becoming much more usual as $p_c$ is approached from above, it is natural to conjecture that $\lambda_c(p) \to \infty$ as $p \uparrow 1$ and that  $\lambda_c(p) \to 0$ as $p \downarrow p_c$; indeed, the first of these is easily established for any $d \geq 2$ and $\ell \in S^{d-1}$ by standard Peierls arguments.
Furthermore, using Remark 1.2, it is easy to see that, for $d=2$ and for any given $\ell \in S^1$, $\lambda_c(p) = (p-p_c)^{\nu + o(1)}$, where $\nu$ is the correlation length exponent (provided that this exponent exists). It is predicted that $\nu = 4/3$ for bond percolation on the square lattice, and this has been established for site percolation on the triangular lattice \cite{SW}. Our arguments apply without significant change to the natural analogue of our model to the latter case (where movement to closed sites is forbidden), so that we indeed obtain $\lambda_c(p) = (p - 1/2)^{4/3 + o(1)}$ as $p \downarrow 1/2$ for this variant of the model. 

\begin{figure}[h]\label{criticalbias}
\centering
\epsfig{file=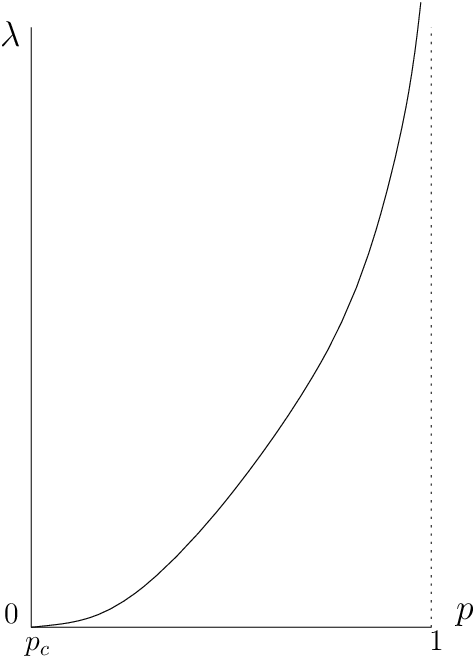, width=5cm}
\caption{The critical curve $\big(p,\lambda_c(p)\big)$ is illustrated in $\Z^d$.}
\end{figure}

In dimension $d=2$, that $\lambda_c:(p_c,1) \to (0,\infty)$ is monotone increasing is observed from Remark 1.2, and monotonicity in the standard coupling of percolation models at differing parameters. This is illustrated in Figure ~4. Although it is natural to believe that this monotonicity holds in every dimension, the formation of a trap - its open interior and its closed surface - presents two requirements that have opposite monotonicities, so that a genuine argument may be needed. We are grateful to Amir Dembo for posing this question. 

\subsection{Trapping and the Anderson model}
Finally, we mention a further important physical context where a phase transition from localized behaviour to transport is anticipated.

For a Schr\"odinger operator with a random ergodic potential and with a spectrum bounded from below,
it is predicted \cite[Section 4.2]{LGP} that the eigenstates associated to eigenvalues close to the bottom of the spectrum
are localized in regions where the potential fluctuates, while, at least in dimension at least three, high energies should be associated to extended states. The phase transition from an {\it insulator region}
of localized states to a {\it metallic region} of extended states is called the {\it metal-insulator transition}.
See \cite{germinetkleinsurvey} for a survey of the mathematical literature relating to this transition in the case of potentials whose randomness has the property of independence at a distance, and \cite{germinetklein} for a characterization of the transition in terms of a local transport exponent. 

In the case of the parabolic Anderson model, a direct probabilistic representation is available, since the model  may then be interpreted as a population of random walkers reproducing or dying at rates determined by a potential that varies randomly in space. The model exhibits massive imbalances in the location of surviving particles at late time, even if the potential has only a small element of randomness. In an effect reminiscent of trapping, most particles cluster in small ``intermittent islands'' at late time. See \cite{HKM}
for a study of different regimes of this intermittency, \cite{morters} for a survey of recent progress on the model when the random potential is heavy-tailed, when an even stronger form of localization occurs, and \cite{MOS} for progress on aging in the parabolic Anderson model. 
Related effects are witnessed in 
Brownian motion in an environment of Poissonian obstables that either kill or repel the motion on contact: see \cite{sznitmanbook}
for a review of progress in this area.

\subsection{Organization of the paper}

Let us now say a word on the order in which we present our results. 

The central percolation estimate is Proposition~\ref{backtrack_exponent}. Its proof having a different flavor from most of the arguments of the paper, we choose to present it in the final Section~\ref{sect_perco}. 

We begin by proving the
part of  Theorem~\ref{the_theorem3} concerned with the lower bound on $\Delta_n$.   The proof, which makes rigorous the outline
presented in Section~\ref{heuristic},   appears in Section~\ref{sect_LB_hit}. The approach of this section sheds enough light on the mechanics of the model 
to explain why the paper's first three theorems should be valid.

We introduce the important concept of regeneration times in Section~\ref{sect_first_regen}. In this section, we give basic properties, and then state and prove the main deduction about regeneration times, Theorem~\ref{expo_ballistic}, which relies on  two results: the exit probability distribution of large boxes
(Theorem~\ref{BL}) and the exit time of large boxes (Proposition~\ref{BL2_PF}).

In Section~\ref{proofs_theorem}, we use Theorem~\ref{expo_ballistic} to prove Theorems~\ref{the_theorem}, \ref{the_theorem2} and~\ref{the_theorem3}.

The exit-time Proposition~\ref{BL2_PF} is proved in Section~\ref{sect_time} using spectral gap estimates and the percolation estimate Proposition~\ref{backtrack_exponent}.

In Section~\ref{sect_BL}, we present the proof of the fundamental dynamic renormalization Theorem~\ref{BL}. 
The proof is substantial and the section begins by sketching its main elements.
A central idea is to introduce a modified walk in order to improve the spectral gap estimate obtained in Section~\ref{sect_time}.
This implies that the new walk leaves a large box quickly, and this, in turn, yields very strong control of exit probabilities.

We have chosen an order to the paper that we believe is most conducive to understanding the results and their proofs; nonetheless, the sections may be read independently of one another, for, when results are quoted from elsewhere in the paper, their statements have almost always been discussed in the introduction.

First of all, we present a section that introduces notation that we will often need.

\noindent{\bf Acknowledgments.}
We would like to thank G\'erard Ben Arous for his essential role in this project, which included 
initiating the research and continual support in conducting it.
We thank the Mittag-Leffler Institute, whose spring 2009 semester on discrete probability 
both authors attended, and which provided an excellent environment to pursue the research.
We further thank Stephen Buckley for comments on an earlier version of the manuscript, Oren Louidor for useful discussions regarding Ornstein-Zernike theory, and a referee for thoroughly reviewing the paper. 

\section{Notation}

Let us denote by $\nu=(e_i)_{i=1\ldots d}$ an orthonormal basis of $\Z^d$ such that $e_1 \cdot\vec \ell \geq e_2\cdot \vec \ell \geq \cdots \geq e_d\cdot \vec \ell \geq 0$; in particular, we have that $e_1\cdot\vec \ell\geq  1 /\sqrt d$. In the case that $d=2$, we use the special notation $\vec{\ell}^{\perp}$ for a unit vector orthogonal to $\vec{\ell}$.

Given a set $V$ of vertices of $\Z^d$, recall from Definition~\ref{defone} that we denote $E(V)=\{ [x,y] \in E(\Z^d): x,y\in V\}$ and
$\partial V= \{x \notin V :  y\in  V,~x\sim y\}$.
This allows us to introduce, for any $L,L' \geq 1$,
\[
\partial^- B(L,L')=\{z\in \partial B(L,L') :  z \cdot \vec{\ell} \leq  -L \},
\]
and
\[
\partial^{\text{s}} B(L,L')= \partial B(L,L')\setminus ( \partial^+ B(L,L') \cup \partial^- B(L,L')) .
\]

We write $d_{\omega}(x,y)$ for the graphical distance in $\omega$ between $x$ and $y$. Define for $x\in \omega$ and $r>0$
\[
B_{\omega}(x,r)=\{y\in \omega:~d_{\omega}(x,y)\leq r\}.
\]

For $A$ a set of vertices in a graph $\omega$, let $\omega\setminus A$ denote the configuration formed from $\omega$ by closing all edges adjacent with at least one endpoint in $A$.

We set
\[
\mathcal{H}^+(k)=\{x\in \Z^d: x\cdot \vec{\ell} \geq k\} \text{ and } \mathcal{H}^-(k)=\{x\in \Z^d: x\cdot \vec{\ell} \leq k\},
\]
as well as
\[
\mathcal{H}^+_x=\mathcal{H}^+(x\cdot \vec{\ell}) \text{ and }\mathcal{H}^-_x=\mathcal{H}^-(x\cdot \vec{\ell}) .
\]

Let us denote for $x\in \Z^d$
\begin{align*}
\mathcal{I}_{x}^-&=\{\omega\in \Omega:~x \text{ belongs to an infinite cluster of } \mathcal{H}^-_x\\ 
 & \text{ induced by the restriction of $\omega$ to edges between vertices in } \mathcal{H}^-_x\},
\end{align*}
and $\mathcal{I}_{x}^+$ using the same definition replacing $\mathcal{H}^-_x$ with $\mathcal{H}^+_x$. Furthermore, we use $\mathcal{I}^-$ (resp.~$\mathcal{I}^+$) for $\mathcal{I}^-_0$ (resp.~$\mathcal{I}^+_0$).

We denote by $\{x \leftrightarrow y \}$ the event that $x$ and $y$ are connected in $\omega$. Accordingly, we denote 
by $K^{\omega}(x)$ the cluster (or connected component) of~$x$ in~$\omega$.

We introduce the following notation. For any set of vertices $A$ of a certain graph on which a random walk $X_n$ is defined, we denote
\[
T_A=\inf\{n\geq 0:~X_n\in A\},~T_A^+=\inf\{n\geq 1:~X_n\in A \} 
\]
\[
\text{ and } \, \, T_A^{\text{ex}}=\inf\{n\geq 0:~X_n\notin A\}.
\]
For $x$ a vertex of the graph, we write $T_x$ in place of $T_{\{ x \}}$ and similarly write $T_x^+$ and $T_x^{\text{ex}}$.
Note that the hitting time $\Delta_n$ defined in Theorem \ref{the_theorem3}
is given by $\Delta_n=T_{\mathcal{H}^+(n)}$.

Furthermore, the pseudo-elliptic constant $\kappa_0=\kappa_0(\ell,d)>0$ will denote
\begin{equation}
\label{kap0}
\kappa_0=\frac{e^{-e_1\cdot \ell}}{\sum_{e\in \nu} e^{-e\cdot \ell}},
\end{equation} 
which is the minimal non-zero transition probability.

A path is understood to be nearest-neighbor unless otherwise stated. A path will be called simple if it visits no vertex more than once.

In this paper, constants are denoted by $c$ or $C$ without emphasizing their dependence in $d$, $\lambda$, $p$ or $\vec{\ell}$. Moreover, the value of those constants may change from line to line.

\section{Lower bound on the hitting time}\label{sect_LB_hit}

We discussed heuristically in Section~\ref{heuristic} how to obtain an estimate on the time spent in traps. A lower bound on $\Delta_n$, the time to reach distance $n$ in the direction $\vec{\ell}$, can be obtained by considering only certain traps. It turns out that it is not difficult to make an explicit choice of such traps so that the resulting lower bound is sharp. This is what we do in this section. The main result is now given. 
\begin{proposition}
\label{LB}
We have that
\[
\liminf \frac{\log \Delta_n}{\log n} \geq 1/\gamma,\qquad \PR\text{-almost surely}.
\]
\end{proposition}

The traps that we will consider are all one-headed traps, each having a unique entrance through which the walk must pass to fall inside. These traps are convenient because certain electrical resistance formulae are available to describe how the walk moves inside them.

We say that there is a one-headed trap $\mathcal{T}^1(x)$ with head $x\in K_{\infty}$, if 
\begin{enumerate}
\item $[x,x+e_1]$ is open,
\item $\vert K^{\omega([x,x+e_1])}(x+e_1) \vert <\infty $, and
\item $(x+e_1)\cdot \vec{\ell} \leq y\cdot \vec{\ell}$ for $y\in K^{\omega([x,x+e_1])}(x+e_1)$,
\end{enumerate}
where $\omega([x,x+e_1])$ is the environment coinciding with $\omega$ except that the edge $[x,x+e_1]$ is closed. 
We set $\mathcal{T}^1(x)=K^{\omega([x,x+e_1])}(x+e_1)$.
If these conditions are not verified,
 we set $ \mathcal{T}^1(x) = \emptyset$. See Figure~5.

\begin{figure}[h]
\label{oneheadpic} 
\centering\epsfig{file=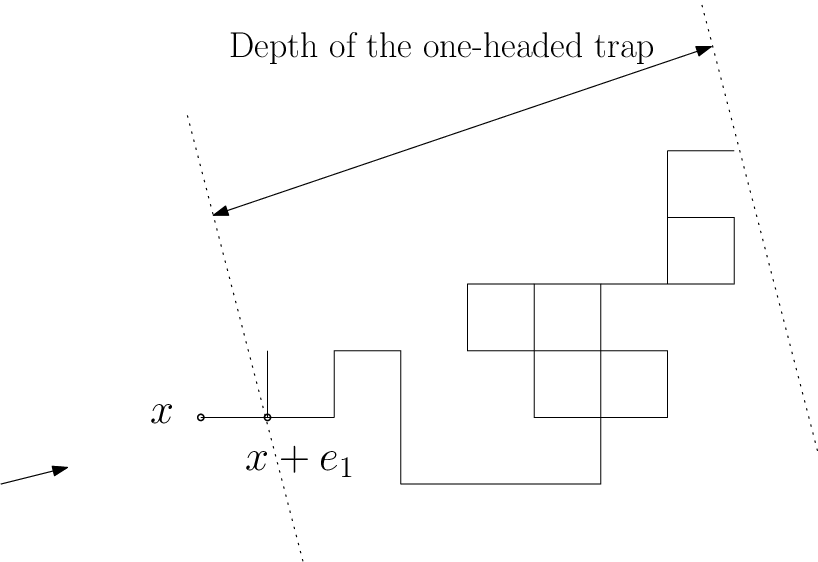, width=8cm}
\caption{An example of a one-headed trap.}
\end{figure}

The furthest distance $\mathfrak{D}(x)$
that the walk may move in the direction of the drift inside a trap from a given point $x$ is given by
\[
\mathfrak{D}(x)=\begin{cases} 0 & \text{ if } \omega \in \mathcal{I}_x^+ \\ 
                                   \max_{\mathcal{P}\in \mathcal{P}_x^+}\max_{y\in \mathcal{P}}  (y-x)\cdot \vec \ell & \text { otherwise,}
                                   \end{cases}
 \]
 where $\mathcal{P}_x^+$ is the set of simple paths that start at $x$ and remain in $\mathcal{H}_x^+$. (Note that $\mathfrak{D}(x) = 0$ if $\omega \in \mathcal{I}_x^+$ because $x$ does not belong to a trap in this case.) We call $\mathfrak{D}(x+e_1)$ the depth of the trap $\mathcal{T}^1(x)$. Note that the trap depth may differ from the trap height that we discussed in the introduction (after~(\ref{def_backtrack})).
 
We recall that $\hatxi$ was defined in Remark~\ref{def_backtrack}. We will use the following estimate on the depth of a trap.
\begin{lemma}
\label{estim_height}
For any $\epsilon > 0$, and for all $t$ sufficiently high,
\[
P_p[\mathcal{T}^1(x) \neq \emptyset \text{ and } \mathfrak{D}(x+e_1)\geq t ]\geq \exp(-(1+\epsilon)\hatxi t).
\]
\end{lemma}

Although the proof of this lemma is quite standard, we defer it to the end of Section~\ref{D3} (for $d\geq 3$) and of Section ~\ref{D2} (for $d=2$), since the result is a percolation estimate.

We now quickly sketch the proof of Proposition \ref{LB}. 
In order to reach $\mathcal{H}^+(n)$, the walk must travel through  $n/(\log n)^{3}$ slabs of width $(\log n)^{3}$. By Lemma~\ref{estim_height}, we see that, typically, there are several slabs in which the walk will reach the entrance of  a one-headed trap whose depth is at least  $(1-\epsilon)\hatxi^{-1} \log n$. With reasonable probability, the walk will fall to the base of any given one of these traps, and will then spend at least $n^{(1-\epsilon)\gamma}$ units of time in that trap. Electrical network calculations will be used to reach these conclusions.

\begin{figure}[h]
\centering
\epsfig{file=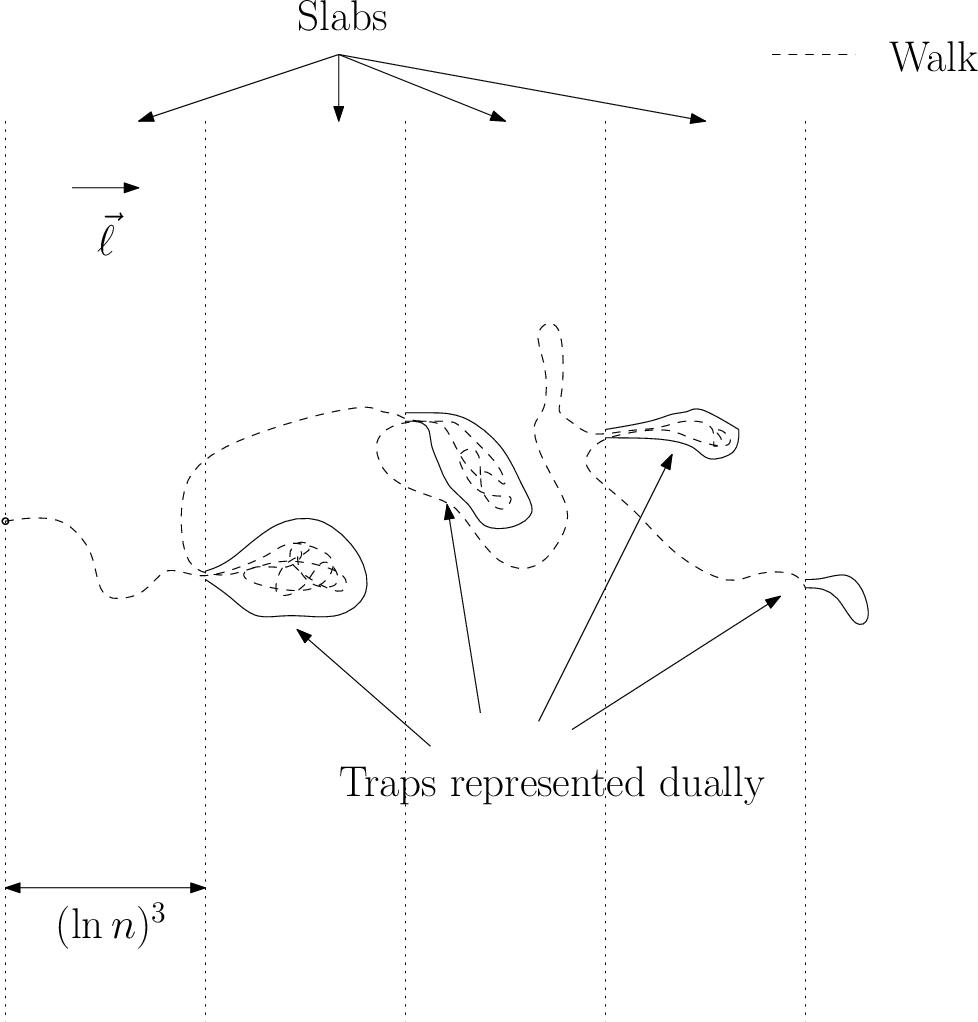, width=10cm}
\caption{The traps encountered in successive slabs in the proof of Proposition~\ref{LB}.}
\end{figure}

\begin{remark}
\label{rem_anteperco}
To prove Proposition~\ref{LB}, it is enough to derive the same statement for the law
 $\tilde{\mathbb{P}}_p :=P_p[~\cdot \mid \mathcal{I}^-] \times P^{\omega}[~\cdot~]$. 
To explain rigorously why this is so, it is useful to invoke the notion of regeneration times for the walk. The next section is devoted to introducing this notion, and, in Remark \ref{rem_perco}, we will provide the needed explanation. 
Informally, the reason that the altered form suffices is that, on arriving for the first time at a high coordinate in the direction $\vec\ell$,
it is highly likely that the walk lies in an infinite open cluster contained in the past half-space, so that the almost sure deduction made about the law $P_p[~\cdot \mid \mathcal{I}^-] \times P^{\omega}[~\cdot~]$ yields the same conclusion for $\PR$.
\end{remark}

\noindent{\bf Proof of Proposition \ref{LB}.}
As we have noted, we need only prove the result for the law 
 $\tilde{\mathbb{P}}_p =P_p[~\cdot \mid \mathcal{I}^-] \times P^{\omega}[~\cdot~]$. 
Note that the configuration in the hyperplane $\mathcal{H}^+(1)$ is unconditioned percolation.

We fix $\epsilon>0$. 

\vspace{0.5cm}

{\it Step 1: The number of one-headed traps met.}

\vspace{0.5cm}

For each $i \geq 1$, we define the slab
\[
S_i=\Big\{x\in \Z^d:~x\cdot \vec{\ell} \in \big[i(\log n)^{3}, (i+1)(\log n)^{3}\big)\Big\}.
\]

As noted in (\ref{def_dir_trans}), the walk $X_n$ is transient in the direction $\vec{\ell}$. 
Hence, we may define 
 \[
 T_i=T_{S_i}<\infty\text{ and }Y_i=X_{T_i}.
 \]
 
We also denote by
  \[
  Z_i=\1{ \mathcal{T}^1(Y_i)\neq \emptyset \text{ and } \mathfrak{D}(Y_i+e_1)\geq (1-\epsilon) \hatxi^{-1} \log n \text{ and }\abs{\mathcal{T}^1(Y_i)} \leq (\log n)^{2}}
\]
the event that, as the walk enters the $i$-th slab, it encounters 
a one-headed trap that is of the ``critical"  depth 
and that is not too large. The latter requirement is made in order that the trap not  enter the next slab.

Note that, at time $T_i$, 
the region $\mathcal{H}^+(Y_i)$ is  unvisited by the walk, and has the distribution of  
unconditioned percolation. Hence, the random variables $Z_i$ are independent. They are identically distributed, with common distribution
\[
\tilde{\PR}[Z_i=1]=P_p[  \mathcal{T}^1(0)\neq \emptyset \text{ and } \mathfrak{D}(e_1)\geq (1-\epsilon)\hatxi^{-1} \log n \text{ and }\abs{\mathcal{T}^1(0)} \leq (\log n)^{2}].
\]

Recall that $\omega_{[0,e_1]}$ denotes the environment formed from $\omega$ by insisting that the edge 
$[0,e_1]$ be closed. 
Note that $\abs{\mathcal{T}^1(0) }> (\log n)^{2}$ under $\omega$ implies that, under $\omega_{[0,e_1]}$, the vertex $e_1$ lies in an open cluster that is finite and has size at least $(\log n)^{2}$. Hence,
\begin{align*}
P_p[\abs{\mathcal{T}^1(0)} > (\log n)^{2}]& \leq  P_{p}[ (\log n)^{2} \leq \abs{K^{\omega_{[0,e_1]}}(e_1)} <\infty ] \\
                                     &= \frac 1{1-p} P_{p}[[0,e_1]\text{ is closed}, (\log n)^{2} \leq \abs{K^{\omega_{[0,e_1]}}(e_1)} <\infty ] \\
                                       &\leq  \frac 1{1-p} P_{p}[(\log n)^{2} \leq \abs{K(e_1)} <\infty ]    \leq \exp(-c(\log n)^{2}),
\end{align*}
the final inequality
by Theorem 6.75 in~\cite{Grimmett}. By Lemma~\ref{estim_height}, we see that
\begin{align*}
\tilde{\PR}[Z_i=1] & \geq P_p[ \mathcal{T}^1(0)\neq \emptyset \text{ and } \mathfrak{D}(e_1)\geq (1-\epsilon)\hatxi^{-1} \log n ] - P_p[(\log n)^{2}>\abs{\mathcal{T}^1(0)}]\\
                     & \geq cn^{\epsilon-1}.
\end{align*}

This means that $ \text{Bern}(cn^{\epsilon-1})$ is stochastically dominated by $Z_i$. Hence, denoting 
\begin{equation}
\label{def_an}
A(n)=\{\text{card}\{i \leq \lfloor n/(\log n)^{3} \rfloor:~Z_i=1\}\geq n^{\epsilon/2}\},
\end{equation}
and recalling that the sequence $Z_i$ is independent, we find that
\begin{equation}
\label{LB_LBB}
\tilde{\PR}[A(n)^c] \leq P[\text{Bin}(\lfloor n/(\log n)^{3} \rfloor,n^{\epsilon-1})\leq n^{\epsilon/2}]=o(1).
\end{equation}

That is,  before  reaching $\mathcal{H}^+(n)$, the walk will typically
 enter at least $n^{\epsilon/2}$ slabs in which it will 
 encounter a one-headed trap whose depth is at least $(1-\epsilon) \hatxi^{-1} \log n$.

\vspace{0.5cm}

{\it Step 2: The probability of entering and exiting the trap.}

\vspace{0.5cm}

Consider a one-headed trap $\mathcal{T}^1(x)$ such that $\mathfrak{D}(\mathcal{T}^1(x)) \geq (1-\epsilon)\hatxi^{-1} \log n$ and with $\abs{\mathcal{T}^1(x)}\leq (\log n)^2$.

Let $\delta_x$ denote the point maximizing $y\cdot \vec{\ell}$ among $y \in K^{\omega([x,x+e_1])}(x+e_1)$. In the case of there being several such points, we pick one according to some predetermined order on $\Z^d$. The point $\delta_x$  represents the base of the trap. We aim to show that the walk is likely to reach the base of the trap, and that, from there, it is hard to escape the trap.

Let us denote $T_x=\inf\{n\geq 0,~X_n=x\}$ and $T_x^+=\inf\{n \geq 1,~X_n=x\}$ for any $x\in \Z^d$. In terms used within the theory of finite electrical networks (see~\cite{DS} or~\cite{LP}), we know that
\begin{equation}
\label{LB_eqa}
P^{\omega([x,x+e_1])}_{x+e_1}[T_{\delta_x}<T_{x+e_1}^+]=\frac{C^{\omega([x,x+e_1])}(x+e_1\leftrightarrow \delta_x)}{\pi^{\omega([x,x+e_1])}(x+e_1)},
\end{equation}
and
\begin{equation}
\label{LB_eqb}
P^{\omega([x,x+e_1])}_{\delta_x}[T_{x+e_1}<T_{\delta_x}^+]=\frac{C^{\omega([x,x+e_1])}(x+e_1\leftrightarrow \delta_x)}{\pi^{\omega([x,x+e_1])}(\delta_x)},
\end{equation}
where $C^{\omega([x,x+e_1])}(x+e_1\leftrightarrow \delta_x)$ denotes the effective conductance between $x+e_1$ and $\delta_x$ as defined in \cite[Section 2.2]{LP}. Rayleigh's monotonicity principle \cite[Section 2.4]{LP} states that the effective conductance increases when individual conductances in the graph are increased. In particular, collapsing vertices together increases effective conductances.

Estimating $C^{\omega([x,x+e_1])}(x+e_1\leftrightarrow \delta_x)$ is our next step. In the graph obtained by  collapsing the vertices of $K^{\omega([x,x+e_1])}(x+e_1)\setminus \{x+e_1\}$ to $\delta_x$, the vertex $\delta_x$ is connected to $x+e_1$ by some number $i\in[1,2d]$ of edges in parallel each of whose conductance is at most $C \exp(2\lambda x \cdot \vec{\ell})$. Hence, by Rayleigh's monotonicity principle, 
\[
C^{\omega([x,x+e_1])}(x+e_1,\delta_x) \leq  C \exp(2\lambda (x+e_1) \cdot \vec{\ell}),
\]
so that~(\ref{LB_eqb}) becomes
\begin{equation}
\label{LB_eq1}
P^{\omega([x,x+e_1])}_{\delta_x}[T_{x+e_1}<T_{\delta}^+] \leq C \exp(2\lambda (x+e_1-\delta_x) \cdot \vec{\ell}) \leq C n^{-(1-\epsilon)\gamma}.
\end{equation}

This statement quantifies the assertion that the walk that reaches the base of the trap finds the trap hard to exit.

We denote by $(\mathcal{P}_{x+e_1}(i))_{0 \leq i \leq l_x'}$ a simple open path from $x+e_1$ to $\delta_x$ in $\mathcal{T}^1(x)$. Closing all edges except for those that connect the consecutive elements of $\mathcal{P}_{x+e_1}$ cannot increase any effective conductance. Hence,
\[
C^{\omega([x,x+e_1])}(x+e_1,\delta_x) \geq C^{\omega([x,x+e_1])}(\mathcal{P}_{x+e_1}),
\]
where $C^{\omega([x,x+e_1])}(\mathcal{P}_{x+e_1})$ is the effective conductances of the path $\mathcal{P}_{x+e_1}$.

However, $\mathcal{P}_{x+e_1}$ is composed of $l_x'\leq  \abs{\mathcal{T}^1(Y_i)} \leq C (\log  n)^{2}$ edges in series. Each of them having conductance at least $C \exp(2\lambda x \cdot \vec{\ell})$, we obtain by network reduction that
\[
C^{\omega([x,x+e_1])}(x+e_1,\delta_x) \geq  \frac 1 {(\log n)^{2}} \exp(2\lambda x \cdot \vec{\ell}),
\]
so that~(\ref{LB_eqa}) becomes
\begin{equation}
\label{LB_eq2}
P^{\omega([x,x+e_1])}_{x+e_1}[T_{\delta_x}<T_{x+e_1}^+] \geq \frac c {(\log n)^{2}},
\end{equation}
a bound that shows that, in departing from the trap entrance, the walk has a significant probability of 
 reaching  the trap base.

\vspace{0.5cm}

{\it Step 3: Time spent in a one-headed trap.}

\vspace{0.5cm}

Assuming that there exists a one-headed trap at $x$ such that $\mathfrak{D}(\mathcal{T}^1(x)) \geq (1-\epsilon)\hatxi^{-1} \log n$ and $\abs{\mathcal{T}^1(x)}\leq (\log n)^2$, we would like to give a lower bound on the time spent in that trap; more precisely, we want to estimate $ P^{\omega}_x[X_1=x+e_1,~X_2\in\mathcal{T}^1(x)\text{ and }T_{x+e_1}^+ \geq n^{(1-2\epsilon)\gamma}]$.

Writing $\theta_i$ for the time-shift by $i$, we see that, on the previously described event, if we have 
\begin{enumerate}
\item $X_1=x+e_1$ and $X_2\in\mathcal{T}^1(x)$,
\item $T_{\delta_x}\circ \theta_2 <T_x^+\circ \theta_2$,
\item $\text{card} \{i \in [T_{\delta_x}\circ \theta_2,T_x^+\circ \theta_2],~ X_i=\delta_x \} \geq n^{(1-2\epsilon)\gamma}$,
\end{enumerate}
then $X_1=x+e_1$, $X_2\in\mathcal{T}^1(x)$ and $T_x^+ \geq n^{(1-2\epsilon)\gamma}$.

 By~(\ref{kap0}) and~(\ref{LB_eq2}), we have that
\[
P^{\omega}_x[X_1=x+e_1,~X_2\in\mathcal{T}^1(x)] \geq c \text{ and } P^{\omega}_x[T_{\delta_x}<T_x^+] \geq \frac c {(\log n)^{2}},
\]
 and, by~(\ref{LB_eq1}),
 \[
 P^{\omega}_{\delta_x}[\text{card} \{i \leq T_x^+,~ X_i=\delta_x \} \geq n^{(1-2\epsilon)\gamma}] \geq c.
 \]

Hence, by Markov's property,
 \begin{equation}
\label{LB_LB}
 P^{\omega}_x[X_1=x+e_1,~X_2\in\mathcal{T}^1(x)\text{ and }T_{x+e_1}^+ \geq n^{(1-2\epsilon)\gamma}] \geq  \frac c {(\log n)^{2}},
 \end{equation}
 for $x$ such that $\mathfrak{D}(\mathcal{T}^1(x)) \geq (1-\epsilon) \hatxi^{-1} \log n$ and $\abs{\mathcal{T}^1(x)}\leq (\log n)^2$.

\vspace{0.5cm}

{\it Step 4: Conclusion.}

\vspace{0.5cm}

 For all $j\leq n/(\log n)^{3}-1$ such that $Z_j=1$
\begin{align}
\label{stepistep}
& P^{\omega}_{Y_{j_i}}[ T^{\text{ex}}_{S_j}\circ \theta_{T_j} \geq n^{(1-2\epsilon)\gamma}] \\\nonumber
\geq  & P^{\omega}_{Y_{j_i}}[X_{T_j+1}=Y_j+e_1, ~ X_{T_j+2}\in \mathcal{T}^1(Y_j),~T_{Y_j+e_1}^+\circ \theta_{T_j+2} \geq n^{(1-2\epsilon)\gamma}],
\end{align}
since if the walk enters a one-headed trap encountered directly as it enters the $j$-th slab, the time spent in that slab is larger than the time spent to exit the trap.

Consider the event
\[
B(n)=\{ \text{for all $j\leq n/(\log n)^{3}-1$, } T^{\text{ex}}_{S_j}\circ \theta_{T_{S_j}} \leq n^{(1-2\epsilon)\gamma} \},
\]
and note that
\begin{equation}
\label{bn_delta}
\text{on $B(n)^c$,} \qquad \Delta_n\geq n^{(1-2\epsilon)\gamma}.
\end{equation}

On the event  $A(n)$ that was defined in~(\ref{def_an}), we may find indices  $j_1\leq \ldots \leq j_{n^{\epsilon/2}} \leq n/(\log n)^{3}-1$ for which $Z_{j_i}=1$. We see then that
\begin{align*}
\tilde{\mathbb{P}}_p[A(n),B(n)] & \leq \tilde{\mathbb{P}}_p[  \text{for all $i\leq n^{\epsilon/2}$, } T^{\text{ex}}_{S_{j_i}}\circ \theta_{T_{S_{j_i}}} \leq n^{(1-2\epsilon)\gamma}] \\
 & \leq  E_p \bigl[\prod_{i\leq n^{\epsilon/2}} P^{\omega}_{Y_{j_i}}[ T^{\text{ex}}_{S_{j_i}}\circ \theta_{T_{S_{j_i}}} \leq n^{(1-2\epsilon)\gamma}]\bigr] \\
 & \leq \Bigl(1-c {(\log n)^{-2}}\Bigr)^{n^{\epsilon/2}} = o(1),
\end{align*}
where we used Markov's property and~(\ref{LB_LB}). Now, by~(\ref{LB_LBB}), we see that $\tilde{\mathbb{P}}_p [B(n)]=o(1)$. By~(\ref{bn_delta}), this means that
\[
\tilde{\mathbb{P}}_p \Bigl[\liminf \frac{ \log \Delta_n}{\log n} \geq (1-2\epsilon)/\gamma\Bigr]\to 1,
\]
and letting $\epsilon$ go to 0, we obtain the lemma for the measure $\tilde{\mathbb{P}}_p$ and thus for $\PR$ as we previously argued.
\qed

\section{Estimates on the tail of the first regeneration time}
\label{sect_first_regen}

\subsection{Definition of the regenerative structure}

Directionally transient random walks in random environment usually have a regenerative structure, which provides a key tool for their analysis, see~\cite{SZ}. 

Informally, we may define this structure as follows. Consider the first time at which the particle reaches a new maximum in the direction $\vec{\ell}$ which is also a minimum in this direction  for the future trajectory of the random walk; we call this time $\tau_1$. Even though $\tau_1$ is not a stopping time, it has the interesting property of separating the past and the future of the random walk into two independent parts in an annealed sense.

Iterating this construction, we actually obtain a sequence $(\tau_k)_{k\in \N}$ of regeneration times, which separates the walk into independent blocks. This property can be extremely useful, because it reduces the problem of understanding such a random walk in random environment to a sum of independent and identically distributed random variables.  The walk's behavior may be understood particularly well if a good bound on the tail of the common distribution of the  $\tau_i$ is available.

Due to minor and standard issues, to define the $(\tau_k)_{k\in\N}$ formally, we actually need to choose a subsequence of the above construction. The definition is a little cumbersome, and the reader uninterested in these details might wish to skip to the next subsection (on the key properties of regeneration times).

We use the same renewal structure that A-S.~Sznitman used in~\cite{Sznitman}. This construction takes care of correlations in the forward space induced by the conditioning $\mathcal{I}$ by initially waiting for a regeneration time with an infinite cluster behind the particle. The construction also addresses the other minor problem, which is that, due to the $1$-dependence of transition probabilities, the adjacent regeneration blocks are not entirely independent.

For $x\in \Z^d$, $\omega \in \mathcal{I}_x$ and $j\geq 0$, we consider 
$\mathcal{P}_{j,x}(\omega)$, namely, the set of open simple $K_{\infty}$-valued paths of the form  $(\pi(i))_{i\geq 0}$ with 
 $\pi(0)=x$  and $\pi(i)\in \mathcal{H}^+_x$  $i\geq j$. 

The parameter $p$ exceeding $p_c$, percolation takes place on $\mathcal{H}^-_x$ \cite{BGN}, and, $\PR$-a.s.,
 on $\mathcal{I}_x$, $\mathcal{P}_{j,x}(\omega)$ is not empty for large enough $j$. We thus define:
\[
J_x(\omega) =\begin{cases} \inf\{j\geq 0,~\mathcal{P}_{j,x}(\omega)\neq \emptyset \}, &\text{ if }\omega \in \mathcal{I}_x,\\
                          \infty,& \text{ if } \omega \notin \mathcal{I}_x. \end{cases}
\]
                        
We then introduce the sequence of configuration-dependent stopping times and the corresponding successive maxima of the walk in the direction $\vec{\ell}$:
\begin{align*}
& W_0=0,~m_0=J_{X_0}(\omega) \leq \infty, \text{ and by induction,} \\
& W_{k+1}=2+T_{\mathcal{H}^+(m_k)}\leq \infty, m_{k+1}=\sup\{X_n\cdot \vec{\ell},n\leq W_{k+1}\}+1\leq \infty, \\
 & \qquad \qquad \qquad  \qquad \qquad \qquad   \qquad \qquad \qquad   \qquad \qquad \qquad  \text{ for all }k\geq 0.
\end{align*}

We define the collection of edges
\begin{align*}
&  B=\{b\in \mathbb{E}_d:~b=[-e_1,e-e_1] \text{ with $e$ any unit vector of } \\
& \qquad \text{$\Z^d$ such that $e\cdot \vec{\ell} =e_1\cdot \vec{\ell}$}\},
\end{align*}
as well as the stopping times
\begin{align*}
S_1&=\inf\{W_k: ~k\geq 1,~X_{W_k}=X_{W_k-1}+e_1=X_{W_k-2}+2e_1, \text{ and } \\
       &\qquad \omega(b)=1 \text{ for all } b\in B+X_{W_k} \}.
       \end{align*}

We further define
\[
D=\inf\{n\geq 0: X_n\cdot \vec{\ell}<X_0\cdot \vec{\ell}\},
\]
as well as the stopping times $S_k$, $k\geq 0$, $R_k$, $k\geq 1$, and the levels $M_k$, $k\geq 0$:
\begin{align*}
S_0&=0,~M_0=X_0\cdot \vec{\ell},\text{ and for $k\geq 0$}, \\
S_{k+1}&=S_1\circ \theta_{T_{\mathcal{H}^+(M_k)}}\leq \infty,~R_{k+1}=D\circ\theta_{S_{k+1}}+S_{k+1}\leq \infty, \\
M_{k+1}&=\sup_{n\leq R_{k+1}} X_n\cdot \vec{\ell} +1\leq \infty .
\end{align*}

Finally, we define the basic regeneration time
\[
\tau_1=S_K, \text{ with } K=\inf\{k \geq 1:~ S_k<\infty \text{ and } R_k=\infty\}.
\]

Then let us define the sequence $\tau_0=0<\tau_1<\tau_2< \cdots <\tau_k< \ldots$, via the following procedure:
\begin{equation}
\label{regen_struct}
\tau_{k+1}=\tau_1+\tau_k(X_{\tau_1+ \cdot }-X_{\tau_1},~\omega(\cdot+X_{\tau_1})),~k\geq 0.
\end{equation}

That is, $k+1$-st regeneration time is the $k$-th such time after the first one.

We also introduce $\mathcal{O}_B$,  the event that all the edges of $ B$ are open. Under $\PR$, the events $\mathcal{I}_{-e_1}^-$ and $\mathcal{O}_B \cap \{D=\infty \}$ are independent and have positive probability. The following conditional measure, that plays an important role for the renewal property, is thus well-defined:
\begin{equation}
\label{def_prob_regen}
\PR[~\cdot \mid 0-\text{regen}]:=\PR[~\cdot~\mid \mathcal{I}_{-e_1}^-,\mathcal{O}_B,D=\infty].
\end{equation}

\subsection{Key properties of regeneration times}

The first main result is that the regeneration structure exists and is finite. (See Lemma 2.2 in~\cite{Sznitman}.)
\begin{lemma}
\label{tau_finite}
For any $k\geq 1$, we have ${\bf P}_p$-a.s., for all $x\in K_{\infty}(\omega)$,
\[
\tau_k<\infty, \qquad \qquad P_x^{\omega}-\text{a.s.}
\]
\end{lemma}

The fundamental renewal property is now stated: see  Theorem 2.4 in~\cite{Sznitman}.
\begin{theorem}
\label{tau_indep}
Under $\PR$, the processes $(X_{\tau_1\wedge\cdot}), (X_{(\tau_1+\cdot)\wedge \tau_2}-X_{\tau_1}),\cdots, (X_{(\tau_k+\cdot)\wedge \tau_{k+1}}-X_{\tau_k}),\ldots$ are independent and, except for the first one, are distributed as $(X_{\tau_1\wedge \cdot} )$ under $\PR[~\cdot \mid 0-\text{regen}]$.
\end{theorem}

That the conditional distribution of the particle's future has the specific form $\PR[~\cdot \mid 0-\text{regen}]$ (defined in~(\ref{def_prob_regen})) is included only for completeness and will not be used, so that the reader need not be too concerned about its definition. 

We now turn to the main property that we will prove regarding regeneration times.

\begin{theorem}
\label{expo_ballistic}
Let $d \geq 2$ and $p > p_c$.
For all $\epsilon > 0$, there exists $C  > 0$ such that, for all $t > 0$,
\[
 \PR[\tau_2-\tau_1 \geq t]  \leq C t^{-(1-\epsilon)\gamma}.
 \]
\end{theorem}

As we have seen, the main difficulty in proving Theorem~\ref{the_theorem} is to establish the upper bound. This will be  a straightforward consequence of 
Theorem \ref{expo_ballistic}. Theorem \ref{expo_ballistic} is also fundamental for Theorems \ref{the_theorem2} and \ref{the_theorem3}. 

Theorem \ref{expo_ballistic} is, in essence, a consequence of Proposition \ref{BL2_PF} and Theorem \ref{BL}. To derive it, we will also use the following assertion regarding regeneration times, which itself is a consequence of Theorem \ref{BL}.

\begin{proposition}
\label{tail_X_tau}
There exist $C > c > 0$ such that, for all $t > 0$,
\[
\PR[(X_{\tau_2}-X_{\tau_1})\cdot \vec{\ell}\geq t] \leq Ce^{-ct}.
\]
\end{proposition}

We may now justify rigorously the reduction described in Remark~\ref{rem_anteperco}.
\begin{remark}\label{rem_perco}
Note that, under each law $\PR$ and $P_p[~\cdot \mid \mathcal{I}^-] \times P^{\omega}[~\cdot~]$, the first regeneration time $\tau_1$ is almost surely finite. As such, we may couple the two laws by identifying the half-space $\mathcal{H}^+$ and future walk trajectories after these first regeneration times. Clearly, under the coupling, the values of the infimum limit appearing in Proposition~\ref{LB} are almost surely equal in each marginal.
\end{remark}

\subsection{Proofs of the results}

Firstly, we give an estimate on the tail of $(X_{\tau_2}-X_{\tau_1})\cdot \vec{\ell}$.

Let $M$ denote a random variable which has the distribution of $\sup_{n\leq T^+_{\mathcal{H}^-(0)}} (X_n-X_0)\cdot \vec{\ell}$ under the law $\PR[~\cdot \mid \mathcal{I}_{-e_1}^-,~\mathcal{O}_B,~D<\infty]$.

Recall Proposition 2.5 in \cite{Sznitman}, a rather general fact that relates the tail of $(X_{\tau_2}-X_{\tau_1})\cdot  \vec{\ell}$ to the tail of an excursion of $X_n$. 
\begin{lemma}
\label{tail_sznitman_lemma}
The random variable  $(X_{\tau_2}-X_{\tau_1})\cdot  \vec{\ell}$ is stochastically dominated by $\Sigma_J$, where $\Sigma_0=0$, and, for $k\geq 1$,
 \[
 \Sigma_k=(2+\overline{M}_1+4\overline{H}_1)+\cdots (2+\overline{M}_k+4\overline{H}_k).
 \]
Here,  $\overline{M}_i$, $\big\{ \overline{H}_i: i\geq 1\big\}$ and $J$ are independent, with $\overline{M}_i$ having the distribution of $M$, and with  $\overline{H}_i$ having the distribution of 
 \[
  H \text{ a geometric variable on $\N^*$ of parameter $1-p^{\abs{B}+1}/(2d)^2$},
  \]
 while $J$ is a geometric variable on $\N^*$ of parameter $\PR[T^+_{\mathcal{H}^-(0)}<\infty \mid \mathcal{O}_B]$.
\end{lemma}

\noindent {\bf  Proof of Proposition~\ref{tail_X_tau}.}
Fix an $\alpha>d+3$. We have that
\begin{align*}
&\PR[2^k\leq M< 2^{k+1}]\\
\leq & \frac 1 {\PR[T_0^+<\infty]}\Bigl[\PR[T_{\partial B(2^k,2^{\alpha k})}\neq T_{\partial^+B(2^k,2^{\alpha k})}] \\
                            & + \PR[X_{T_{\partial B(2^k,2^{\alpha k})}} \in \partial^+B(2^k,2^{\alpha k}), T_{\mathcal{H}^-(0)}^+ \circ \theta_{ T_{\partial^+B(2^k,2^{\alpha k})}}<T_{\mathcal{H}^+(2^{k+1})}\circ \theta_{ T_{\partial^+B(2^k,2^{\alpha k})}}]\Bigr].
\end{align*}

By a union bound on the $2^{\alpha k}$ possible positions of $X_{ T_{\partial^+B(2^k,2^{\alpha k})}}$, and by using translation invariance, we find that
\[
\PR[2^k\leq M< 2^{k+1}]\leq   \frac 1 {\PR[T_0^+<\infty]}\Bigl[  2^{\alpha k} \PR[T_{\partial B(2^{k+1},2^{\alpha (k+1)})}\neq T_{\partial^+B(2^{k+1},2^{\alpha (k+1)})}]+e^{-(2^k)}\Bigr],
 \]
as a consequence of Theorem~\ref{BL}.

Hence, using Theorem~\ref{BL} again, we obtain
\[
\PR[2^k\leq M< 2^{k+1}]\leq   c 2^{\alpha k}e^{-c(2^k)},
 \]
          and
\[
\PR[M\geq t] \leq C e^{-ct}.
\]

Writing $A \preceq B$ to indicate that $B$ stochastically dominates $A$, we see that Lemma~\ref{tail_sznitman_lemma} implies
\[
(X_{\tau_2}-X_{\tau_1})\cdot \vec{\ell} \preceq \sum_{i=1}^{\text{Geom}(\PR[T^+_{\mathcal{H}^-(0)}\mid \mathcal{O}_B])} S_i,
\]
where $S_i=2+\overline{M}_i+4\overline{H}_i$ is an independent and identically distributed sequence. Moreover, $\PR[S_i \geq t ] \leq C e^{-ct}$, by the similar estimate on $M$. The statement of the proposition follows.
\qed \medskip

\noindent {\bf  Proof of Theorem~\ref{expo_ballistic}.}
Denote by $\chi$ the smallest integer such that $\{X_i-X_{\tau_1},~i\in [\tau_1,\tau_2]\}\subseteq B(\chi,\chi^{\alpha})$.
The idea of the proof is that, if the time spent in a {\em regeneration box} $B(\chi,\chi^{\alpha})$ is large, then 
\begin{enumerate}
\item either the regeneration box is large, which is unlikely by Theorem~\ref{BL},
\item or the walk spends a lot of time in a small area, an event whose probability is controlled by~Proposition~\ref{BL2_PF}.
\end{enumerate}

 We have then that
\begin{align}
\label{expo_step_1}
\PR[\tau_2-\tau_{1}\geq t]&\leq \PR[\chi \leq k,~T^{\text{ex}}_{X_{\tau_{1}}+B(k,k^{\alpha})}\circ \tau_{1}\geq t] +\PR[\chi \geq k] \\ \nonumber 
                                &\leq\PR[T^{\text{ex}}_{B(k,k^{\alpha})}\geq t\mid 0-\text{regen}] +\PR[\chi \geq k] \\ \nonumber
                                &\leq \frac 1 {\PR[ 0-\text{regen}]} \PR[T^{\text{ex}}_{B(k,k^{\alpha})}\geq t] +\PR[\chi \geq k]  .
                                \end{align}

{\it Step 1: The size of the regeneration box.}

Firstly, we notice that
\begin{align}
\label{expo_step_12}
 \PR[\chi \geq k]  \leq & \PR[(X_{\tau_2}-X_{\tau_1})\cdot \vec{\ell}\geq k]  \\ \nonumber
                                      &  +\PR[(X_{\tau_2}-X_{\tau_1})\cdot \vec{\ell} \leq k,~\max_{\tau_1\leq i,j \leq \tau_2} \max_{l\in [2,d]} \abs{ (X_j-X_i)\cdot f_l} \geq k^{\alpha}].
                                      \end{align}
                                      
The first term on the right-hand side may be bounded above  by Proposition~\ref{tail_X_tau}.
The second may be so bounded by the following reasoning: on the event that $(X_{\tau_2}-X_{\tau_1})\cdot \vec{\ell}\leq  k$ and $\max_{j\neq 1}\max_{\tau_{1}\leq j_1,j_2 \leq \tau_{2}} \abs{ (X_{j_1}-X_{j_2})\cdot f_j } \geq  k^{\alpha}$, the walk $X_n$ does not exit the box $B(k,k^{\alpha})$ on the side $\partial^+ B(k,k^{\alpha})$. This implies that
\begin{align*}
&\PR[X_{\tau_2}-X_{\tau_1}\leq k,~\max_{\tau_1\leq i,j \leq \tau_2} \max_{l\in [2,d]} \abs{ (X_j-X_i)\cdot f_l} \geq k^{\alpha} ]\\
 \leq & \frac 1 {\PR[0-\text{regen}]}\PR[T_{\partial B(k,k^{\alpha})}\neq T_{\partial^+ B(k,k^{\alpha})}]\leq Ce^{-ck},
\end{align*}
by Theorem~\ref{BL}.

Applying these bounds to~(\ref{expo_step_12}), we obtain
\begin{equation}
\label{expo_step_15}
 \PR[\chi \geq k]\leq Ce^{-ck}.
 \end{equation}

{\it Step 2: The time spent in a box.}

Proposition~\ref{BL2_PF} implies that, for any $\epsilon>0$, and for $t > 0$ sufficiently high and $k\leq t$
\begin{equation}
\label{expo_step_3}
\PR[T^{\text{ex}}_{B(k,k^{\alpha})}\geq t]\leq  C k^{d\alpha(2\gamma +1)}t^{-(1-\epsilon)\gamma}+e^{-ck}.
\end{equation}

{\it Step 3: Conclusion.}

Finally, recalling~(\ref{expo_step_1}),~(\ref{expo_step_15}) and~(\ref{expo_step_3}), we find that, for any $\epsilon>0$ and  $k\leq t$,
\[
\PR[\tau_2-\tau_1 \geq t]  \leq  C e^{-ck}+  Ck^{d\alpha(2\gamma +1)}t^{-(1-\epsilon)\gamma}, 
\]
so that we obtain, using $k=(\log t)^{2}$,
\[
\PR[\tau_2-\tau_1 \geq t]  \leq C t^{-(1-\epsilon)\gamma}.
\]
\qed \medskip

\section{The proofs of the three main theorems}
\label{proofs_theorem}
We now provide the proofs of Theorems \ref{the_theorem}, \ref{the_theorem2} and \ref{the_theorem3}. 
Each result follows directly from Theorem \ref{expo_ballistic}.

\noindent{\bf Proof of Theorem~\ref{the_theorem3}.} 
By the definition of the $n$-th regeneration time, we necessarily have that $X_{\tau_n}\cdot e_1\geq n$. It follows that 
\[
\Delta_n \leq \tau_n = \tau_1 +\sum_{i=0}^{n-1} (\tau_{i+1} -\tau_i).
\]

Recalling that $\gamma \in  (0,1)$, and by invoking Lemma~\ref{tau_finite}, Theorem~\ref{tau_indep} and Theorem~\ref{expo_ballistic}, a standard argument gives that, for any $\epsilon>0$ and for $n$ large enough,
\[
\tau_n \leq Cn^{1/\gamma +\epsilon}.
\]

The quantity $\epsilon > 0$ being arbitrary, we obtain
 \[
\limsup \frac{\log \Delta_n}{\log n} \leq 1/\gamma, \qquad \PR\text{-a.s.}
\]

Recalling Proposition~\ref{LB}, we see that
\begin{equation}\label{eqdeltan}
\lim \frac{\log \Delta_n}{\log n} = 1/\gamma,  \qquad \PR\text{-a.s.}
\end{equation}

We will now perform a classical inversion argument. We set $\overline{X}_n=\max_{i\leq n} X_i\cdot \vec{\ell}$, and note that
\begin{equation}\label{eqoverxn}
0 \leq \overline{X}_n - X_n\cdot \vec \ell \leq \max_{i \in [0,n-1]} \big( X_{\tau_{i+1}}-X_{\tau_i} \big) \cdot \vec\ell \, \, \text{ for } n\geq \tau_1,
\end{equation}
as well as 
\[
\overline{X}_n \geq m \, \,  \textrm{if and only if} \, \, \Delta_m \leq n
\]
for each $n,m \in \N$.   Using this property alongside (\ref{eqdeltan}), we find that, for any $\epsilon > 0$,
\begin{equation}\label{blablablabla}
n^{\gamma-\epsilon} \leq \overline{X}_n \leq n^{\gamma+ \epsilon}, \qquad \PR\text{-a.s.}
\end{equation}
 for $n$ large enough.

Moreover, by Proposition~\ref{tail_X_tau} and the Borel-Cantelli lemma, we find that, for $n$ large enough,
\[
\max_{i \in [0,n-1]} \big( X_{\tau_{i+1}}-X_{\tau_i} \big)  \leq (\log n)^2 \qquad \PR\text{-a.s.}
\]

Using this bound along with~(\ref{blablablabla}) in~(\ref{eqoverxn}), and recalling that $\tau_1 < \infty$, we see that, for $n$ large enough,
\[
n^{\gamma-\epsilon} \leq X_n \cdot \vec \ell \leq n^{\gamma+ \epsilon}, \qquad \PR\text{-a.s.},
\]
for any $\epsilon>0$. Hence, by taking $\epsilon$ to $0$,
\[
\lim \frac{\log X_n\cdot \vec{\ell}}{\log n} = \gamma,  \qquad \PR\text{-a.s.}
\]
\qed \medskip

\noindent{\bf Proof of Theorem~\ref{the_theorem}.} 
By the existence of the regenerative structure, we know that, by Theorem 3.4 in~\cite{Sznitman},
\begin{equation}
\label{LLNNN}
\lim \frac{X_n}n = v = \frac{\ES[X_{\tau_2}-X_{\tau_1}]}{\ES[\tau_2-\tau_1]}.
\end{equation}

By our definition of regeneration times, it is clear that $\ES[X_{\tau_2}-X_{\tau_1}]\cdot \vec{\ell} \geq 1$, and thus, if $\ES[\tau_2-\tau_1]<\infty$, then $v \cdot \vec{\ell} >0$.

 For $\gamma>1$, it follows by taking $\epsilon > 0$ small enough in Theorem~\ref{expo_ballistic} that
\[
\PR[\tau_2-\tau_1 \geq n] \leq Cn^{-(1+\epsilon)} \text{ and } \ES[\tau_2-\tau_1]<\infty.
\]

Hence, if $\gamma>1$ then $v \cdot \vec{\ell} >0$. Conversely, for $\gamma<1$, we see that, by Proposition~\ref{LB},
 $\liminf \Delta_n / n=\infty$. Moreover, by Proposition~\ref{tail_X_tau}, we know that $\ES[X_{\tau_2}-X_{\tau_1}]\cdot \vec{\ell}<\infty$, and so, by a standard inversion argument,
the limit 
$\lim n^{-1}\Delta_n$ exists and is given by  $\frac{\ES[\tau_2-\tau_1]}{\ES[X_{\tau_2}-X_{\tau_1}]\cdot \vec{\ell}}$.
This means that $\ES[\tau_2-\tau_1]=\infty$, so that, in light of~(\ref{LLNNN}),  $v$ equals $\vec 0$. 
\qed \medskip

\noindent{\bf Proof of Theorem~\ref{the_theorem2}.} 
For $\gamma>2$, it follows by taking $\epsilon > 0$ small enough in Theorem~\ref{expo_ballistic} that
\[
\PR[\tau_2-\tau_1 \geq n] \leq Cn^{-(2+\epsilon)} \text{ and } \ES[(\tau_2-\tau_1)^2]<\infty.
\]

This integrability condition, coupled with arguments presented in the proof of Theorem 3.4 of~\cite{Sznitman}, yields the result.
\qed \medskip

\section{Exit time of a box}
\label{sect_time}

Here we prove Proposition~\ref{BL2_PF}. The method used to prove this proposition is a minor variation of an argument presented by A-S.~Sznitman in~\cite{Sznitman}.

We remind the reader that $P^{\omega}$, the transition operator, maps $L^2(\pi^{\omega})$ into $L^2(\pi^{\omega})$. We are interested in the principal Dirichlet eigenvalue  $\Lambda_{\omega}(B(L,L^{\alpha}))$ of $I-P^{\omega}$ in $B(L,L^{\alpha})\cap K_{\infty}$, which is given by
\begin{equation}
\label{prop_dirichlet}
\Lambda_{\omega}(B(L,L^{\alpha}))=\begin{cases} & \inf\{\mathcal{E}(f,f), f_{\mid (B(L,L^{\alpha})\cap K_{\infty})^c}=0, \abs{\abs{f}}_{L^2(\pi)}=1\}, \\ & \qquad  \text{ when $B(L,L^{\alpha})\cap K_{\infty}\neq \emptyset$,} \\ 
                                  &  \infty, \text{ by convention when } B(L,L^{\alpha})\cap K_{\infty} =\emptyset , \end{cases}
                                  \end{equation}               
  where the Dirichlet form is defined for $f,g\in L^2(\pi)$ by
  \[
  \mathcal{E}(f,g)=(f,(I-P^{\omega})g)_{\pi} =\frac 12 \sum_{\abs{x-y}=1} (f(y)-f(x))(g(y)-g(x)) c^{\omega}([x,y]).
  \]                               
               
               Define
\begin{equation}
\label{BL2_def_H_max}
\mathcal{B}\mathcal{K}(L,L^{\alpha})=\max_{x\in B(L,L^{\alpha})} \mathcal{B}\mathcal{K}(x).
\end{equation}
                   
  Lemma~\ref{BL2_LB_Dirichlet} will provide a lower bound on $\Lambda_{\omega}(B(L,L^{\alpha}))$ in terms of maximal trap height $\mathcal{B}\mathcal{K}(L,L^{\alpha})$. The Perron-Frobenius theorem will then  provide an upper bound on  the exit time of $B(L,L^{\alpha})$ in terms of $\Lambda_{\omega}(B(L,L^{\alpha}))$ and so complete the proof of Proposition~\ref{BL2_PF}. 


\begin{lemma}    
\label{BL2_LB_Dirichlet}
For $\omega$ such that $0\in K_{\infty}$,
 \[
 \Lambda_{\omega}(B(L,L^{\alpha})) \geq c  L^{-2d\alpha}e^{-2\lambda\mathcal{B}\mathcal{K}(L,L^{\alpha})}.
 \]
\end{lemma}
\noindent{\bf Proof.}
By the definition of $\mathcal{B}\mathcal{K}(L,L^{\alpha})$, we may construct, for every $x\in B(L,L^{\alpha})\cap K_{\infty}$,
a simple open path $p_x = (p_x(i))_{0\leq i\leq l_x}$ with the property that $p_x(0) = x$, $p_x(i) \in B(L,L^{\alpha})\cap K_{\infty} $ for $0 < i < l_x$ and
  $p_x(l_x) \in \partial B(L,L^{\alpha})$, and such that
\[
\max_{i\leq l_x} (x-p_x(i))\cdot \vec{\ell}\leq \mathcal{B}\mathcal{K}(L,L^{\alpha}).
\]

We use a classical argument of Saloff-Coste~\cite{Saloff-Coste}. For $f$ such that  $\abs{\abs{f}}_{L^2(\pi)}=1$ and $f_{\mid (B(L,L^{\alpha})\cap K_{\infty})^c}=0$,
we have that
\begin{align}
\label{arg_SC}
1&=\sum_x f^2(x)\pi_{\omega}(x) = \sum_x \Bigl[\sum_i f(p_x(i+1)) -f(p_x(i))\Bigr]^2\pi_{\omega}(x) \\ \nonumber
  & \leq \sum_x l_x\Bigl[\sum_i \bigl(f(p_x(i+1))-f(p_x(i))\bigr)^2\Bigr]\pi_{\omega}(x) .
  \end{align}

Now, since  $\pi_{\omega}(x) /c_{\omega}([p_x(i+1), p_x(i)]) \leq e^{2\lambda \mathcal{B}\mathcal{K}(L,L^{\alpha})}/\kappa_0$, 
\[
1\leq \frac{e^{2\lambda \mathcal{B}\mathcal{K}(L,L^{\alpha})}}{\kappa_0} \sum_{b=\{y,z\}} (f(z)-f(y))^2 c_{\omega}(b) \times \max_b \sum_{x\in K_{\infty}\cap  B(L,L^{\alpha}),b\in p_x} l_x 
 \]
 where $b\in p_x$ means that $b=\bigl[p_x(i),p_x(i+1)\bigr]$ for some $i \in [0,l_x - 1]$. By virtue of
 \begin{enumerate}
 \item $l_x\leq C L^{1+\alpha(d-1)}$ for any $x\in B(L,L^{\alpha})$, and the fact that
 \item $b=[x,y]\in E(\Z)^d$ may only be crossed by paths  $p_z$ for which $z\in B(L,L^{\alpha})$,
 \end{enumerate}
 we have that
 \[
 \max_b \sum_{x\in K_{\infty}\cap  B(L,L^{\alpha}),b\in p_x} l_x  \leq \rho_d L^{2(1+\alpha(d-1))},
 \]
 and also that
 \[
 1\leq \kappa_0^{-1} \rho_d e^{2\lambda \mathcal{B}\mathcal{K}(L,L^{\alpha}) }   L^{2(1+\alpha(d-1))} \sum_{b=\{y,z\}} (f(z)-f(y))^2 c_{\omega}(b).
 \]
 
Hence, by (\ref{prop_dirichlet}),
 \[
 \Lambda_{\omega}(B(L,L^{\alpha})) \geq c L^{-2d\alpha}e^{-2\lambda \mathcal{B}\mathcal{K}(L,L^{\alpha})}.
 \]
\qed \medskip

\begin{remark}
This lemma confirms our intuition that the key element of the geometry of traps is the height of the trap. Indeed, the Dirichlet eigenvalue (which is in essence the inverse of the typical exit time), is exponentially decaying in the height of traps, whereas all other effects are only polynomial and can be neglected for our purpose. 
\end{remark}

\begin{lemma}
\label{BL2_Lambda}
For every $\epsilon>0$ and for $L: (0,\infty) \to (0,\infty)$ 
satisfying $\lim t^{-1} \log L(t)  =0$,  we have
that, for $t$ sufficiently high,
\[
{\bf P}[ \Lambda_{\omega}(B(L,L^{\alpha})) \leq t^{-1} ] \leq C t^{-(1-\epsilon)\gamma } L^{d\alpha(2\gamma+1)}.
\]
\end{lemma}

\noindent{\bf Proof.}
We may use a union bound argument along with Proposition~\ref{backtrack_exponent} and Lemma~\ref{BL2_LB_Dirichlet} to prove that, for all $\epsilon>0$, there exists $C<\infty$, such that, for any $t>1$,
\[
{\bf P}_p[t< \mathcal{B}\mathcal{K}(L,L^{\alpha})] \leq C L^{1+(d-1)\alpha} e^{-(1-\epsilon)t \hatxi}.
\]

Thus, we obtain
\begin{align*}
 {\bf P}[ \Lambda_{\omega}(B(L,L^{\alpha})) \leq  t^{-1}] 
   \leq &  {\bf P}[\mathcal{B}\mathcal{K}(L,L^{\alpha})> (1/2\lambda) \log (c t L^{-2d\alpha}) ] \\
                                                                                \leq &  C t^{-(1-\epsilon)\gamma } L^{d\alpha (2\gamma+1)} , 
 \end{align*}
since $\gamma=\hatxi/(2\lambda)$.
\qed \medskip

\noindent{\bf Proof of Proposition \ref{BL2_PF}.} 
Note that, for $\omega\in \mathcal{I}$,  the operator on $L^2(\pi_{\omega})$ given by $\1{B(L,L^{\alpha}) \cap K_{\infty}} P^{\omega} \1{B(L,L^{\alpha}) \cap K_{\infty}}$ is associated with the walk killed on leaving $B(L,L^{\alpha})$. Recall that Perron-Frobenius' theorem states that this operator has norm given by the operator's maximum positive eigenvalue, which is $1-\Lambda_{\omega}(B(L,L^{\alpha}) )$. This yields that, for any $x\in K_{\infty} \cap B(L,L^{\alpha})$,
\begin{align*}
& \pi_{\omega}(x)P_x^{\omega}[T^{\text{ex}}_{B(L,L^{\alpha})}>n] \\ 
= & (\1{x}, (\1{B(L,L^{\alpha}) \cap K_{\infty}} P^{\omega} \1{B(L,L^{\alpha}) \cap K_{\infty}} )^n \1{B(L,L^{\alpha})})_{L^2(\pi_{\omega})} \\  
\leq& (\pi_{\omega}(x) \pi_{\omega}({B(L,L^{\alpha})}))^{1/2}(1-\Lambda_{\omega}(B(L,L^{\alpha}) ))^{n}.
\end{align*}

Thus, using  $\pi_{\omega}(B(L,L^{\alpha}))\leq C L^{1+(d-1)\alpha} e^{2\lambda L}$, we have that
\[
P_0^{\omega}[T^{\text{ex}}_{B(L,L^{\alpha})}>t] \leq  C L^{(1+(d-1)\alpha)/2}\exp(\lambda L-t\Lambda_{\omega}(B(L,L^{\alpha}) )).
\]

Hence
\[
\PR[T^{\text{ex}}_{B(L,L^{\alpha})} > t]\leq \PR[\Lambda_{\omega}(B(L,L^{\alpha}) ) <2 \lambda L/t]+ C\exp(-\lambda L/2).
\]

Note that  Lemma~\ref{BL2_Lambda} implies that
\[
 \PR[\Lambda_{\omega}(B(L,L^{\alpha}) ) >2 \lambda L/t] \leq  
C L^{d\alpha(2\gamma+1)}t^{-(1-\epsilon)\gamma}.
\]

The statement of the proposition follows from the two preceding assertions. \qed

\section{Exit probability of large boxes}
\label{sect_BL}

This section is devoted to the proof of the dynamic renormalization Theorem~\ref{BL}. Although our argument is quite technical, 
it is inspired by a rather direct approach which, 
in an effort to aid the reader's comprehension, we will begin by outlining. 
A key tool for the approach 
is the Carne-Varopoulos bound~\cite{Carne}, which states the following for $X$ a reversible Markov chain on a finite graph $G$ with invariant measure $\pi$.
Let $x$ and $y$ be any two points in $G$. Then, for each $n \in \N$,
\begin{equation}
\label{CVB}
P_x[X_n=y] \leq 2 \Bigl(\frac{\pi(y)}{\pi(x)}\Bigr)^{1/2} \exp\Bigl(-\frac{d(x,y)^2}{2n}\Bigr),
\end{equation}
where $d(\cdot,\cdot)$ is the graphical distance in $G$. 

This bound is very useful when $n$ is small. For our specific problem, it implies that, for each $n \in \N$,
\begin{eqnarray}
& & P_0^{\omega}[X_{T^{\text{ex}}_{ B(L,L^{\alpha})}} \notin \partial^+ B(L,L^{\alpha}), ~T^{\text{ex}}_{B(L,L^{\alpha})}<n ] \label{eqheur} \\
&\leq & C(p,\lambda) \sum_{k\leq n} \bigg( L^{\alpha (d-1)}e^{-\lambda L}+ 2(d-1)L^{1+\alpha(d-2)}\exp\bigl(-\frac{L^{2\alpha}}k+\lambda L\bigr) \bigg) \nonumber \\
& \leq & C(p,\lambda) n \Bigl( L^{\alpha (d-1)}e^{-\lambda L}+L^{1+\alpha(d-2)}\exp\bigl(-\frac{L^{2\alpha}}n+\lambda L\bigr)\Bigr). \nonumber
\end{eqnarray}

(In the second line, the index $k$ corresponds to the escape time of the walk, and the two summands to the exit of $B(L,L^{\alpha})$ on the side opposite to the direction of the drift and on one of the lateral sides.)
In this way, we see that, if the exit time of $B(L,L^{\alpha})$ is typically small, (by which we mean, much smaller that $L^{2\alpha}$), then the previous equation implies that walk is likely to exit on the preferred side $ \partial^+ B(L,L^{\alpha})$. This reasoning was used by A-S.~Sznitman in~\cite{Sznitman}. He controlled the tail of the exit of $B(L,L^{\alpha})$ using essentially the same arguments as in the proof of Proposition~\ref{BL2_PF}. He obtained Theorem~\ref{BL} with a polynomial decay, and required the hypothesis that the bias of the walk be small.

However, it is clear that, due to trapping, high values of the bias will cause the exit time of the box to be high, specifically, much higher than $L^{2\alpha}$. From this, we see that it is hopeless to try to obtain Theorem~\ref{BL} by a direct approach using the Carne-Varopoulos bound.

This difficulty creates an impasse, but there is a glimmer of hope. Although traps increase the exit time of a box, they do not much change the exit probabilities. The one-headed traps used in Section~\ref{sect_LB_hit} are an excellent example. Indeed, sealing off the heads of all one-headed traps contained in the box would leave unaffected the box-exit probabilities; however, since this operation blocks the entrance to the one-headed traps it might serve to lower significantly the exit time. Our goal will actually be to introduce a modified walk which, in a certain sense, ignores the time spent in all important traps. The new walk would move much more quickly than the old one, and we might
 hope to apply to it the argument just provided in outline. 

More precisely, we will define a \lq\lq bad\rq\rq subset of $\Z^d$ that contains all large traps, and which is composed of finite clusters. The complement of this bad set will correspond to relatively open space. It is natural to view this partition of $\Z^d$ as being a backbone-trap decomposition, where the traps are the clusters that comprise the bad zone.

Our modified walk will be the original walk recorded only in the complement of the bad zone. This walk will remain reversible, but will make non-nearest-neighbor jumps corresponding to sojourns of the original walk in the bad clusters. 
It will transpire that bad clusters are usually small, so that these jumps will be short. This being the case, the chemical distance for the modified walk will be close to that for the original one, meaning that we will be able to adapt the argument involving the Carne-Varopoulos  bound to the modified walk.

After the bad region and the modified walker have been defined, we will apply the proof strategy explained at the beginning of this section to obtain Theorem~\ref{BL}. 




\subsection{A first candidate for the backbone-trap decomposition}

\subsubsection{Definitions}
\label{sect_def}

We will introduce some notations often used for renormalization arguments, see~\cite{Grimmett}. Fix $K\in\N$ and set $B(K)=[-5K/8,5K/8]^d$. The largest open cluster of $B(K)$ will be denoted by $M(K)$. 

We say that $M(K)$ crosses $B(K)$ in the $i$-th direction if $M(K)$ contains an open path having end-vertices $x,y$ satisfying $x\cdot e_i=-\lfloor 5K/8\rfloor $ and $y\cdot e_i=\lfloor 5K/8 \rfloor$.

The cluster $M(K)$ is called a crossing cluster for $B(K)$ if $M(K)$ crosses $B(K)$ in the $i$-th direction for each $i \in [1,\ldots, d]$. 

Let $\eta > 0$. We say that $B(K)$ is a $(K,\eta)$-open box if
\begin{enumerate}
\item $M(K)$ is a crossing cluster for $B(K)$, and
\item $M(K)$ is the unique open cluster $C$ of $B(K)$ satisfying $\abs{C \cap B(K)} \geq \eta K$.
\end{enumerate}

We write $\overline{x} = (x_1,\ldots,x_d)$, so that in fact $\overline{x}$ is equal to $x$. We make this distinction because we will employ the $\overline{x}$ notation to indicate the index of a $K$-box (rather than an actual lattice point).

 We set $B_{\overline{x}}(K)=(K x_1,\ldots,K x_d)+B(K)$, and call $B_{\overline{x}}(K)$
the $K$-box at $\overline{x}$.
We say that $B_{\overline{x}}$
is $(K,\eta)$-open if  its biggest cluster $M_{\overline{x}}(K)$ verifies the two properties mentioned above.

By Theorem 7.61 in~\cite{Grimmett}, we have the following result.
\begin{lemma}
\label{BL_open}
For $d\geq 2$, $p>p_c$, $\overline{x}\in \Z^d$ and $\eta>0$, we have that
\[
{\bf P}_p[\text{$B_{\overline{x}}(K)$ is $(K,\eta)$-open}] \to 1 \qquad \text{as } K\to \infty.
\]
\end{lemma}

In the sequel, we choose $\eta = 3^{-d}/8$ and do not display $\eta$ in notation.

Two vertices $u,v \in \Z^d$ are said to be $*$-neighbors if $\abs{\abs{u-v}}_{\infty}=1$. This topology induces a natural notion of $*$-connected component on vertices.

 The $K$-box at $\overline{y}$ is said to be super-open if the $K$-box at $\overline{x}$ is $(K,\eta)$-open for every $\overline{x}$ satisfying  $\abs{\abs{\overline{x}-\overline{y}}}_{\infty}\leq 1$. If a box is not super-open, we call it weakly-closed.
 
 Denoting  $Y_{\overline{x}}^K=\1{\text{the box at $\overline{x}$ is super-open}}$, we see that $Y_{\overline{x}}^K$ is independent of $Y_{\overline{y}}^K$ for $\overline{y}$ such that $\abs{\abs{\overline{x}-\overline{y}}}_1 \geq 4$.

We say that the $K$-box at $\overline{x}\in \Z^d$ is good (respectively~super-good), if there exists an infinite directed path of 
open (respectively super-open) $K$-boxes starting at $\overline{x}$. That is, we have $\{\overline{x}=\overline{x}_0,\overline{x}_1,\overline{x}_2,\overline{x}_3,\ldots  \}$ such that, for each $i\geq 0$,
\begin{enumerate}
\item we have $\overline{x}_{2i+1} -\overline{x}_{2i} = e_1$ and $\overline{x}_{2i+2}-\overline{x}_{2i+1}\in  \{e_1,\ldots,e_d\}$,
\item there is an open (respectively super-open) $K$-box at $\overline{x_i}$.
\end{enumerate}
If a box is not good, it is said to be bad. If it is not super-good, we call it weakly bad. 

The key property of a good point $x$ will be that there exists a open path $(x_i)_{i\geq 0}$ with $x_0 = x$ such that $x_0\cdot \vec{\ell}<x_1\cdot\vec{\ell}\leq x_2 \cdot\vec{\ell}<x_3\cdot \vec{\ell}\leq \ldots$, and $(x_i-x_0)\cdot \vec{\ell} \geq c_d i$, $i \geq 0$, for some $c_d > 0$.

Let us denote by $T^K(\overline{x})$ the bad *-connected component of $\overline{x}$ for the dependent site percolation defined by $\1{B_{\overline{x}}(L)\text{ is good}}$. 

We now have a provisional candidate for the backbone-trap decomposition. Indeed, we  introduce the {\it provisional backbone} $\mathfrak{B}_{\text{prov}}$ which consists of all points contained in some good box, and the {\it provisional trapping zone} $\mathfrak{T}_{\text{prov}}$, which is the complement of  $\mathfrak{B}_{\text{prov}}$. The candidate has the promise that  from  each point of  $\mathfrak{B}_{\text{prov}}$ issues an infinite path that hardly backtracks.

\subsubsection{Properties}

We now provide two results showing that the proposal for the backbone-trap decomposition has some attractive properties: the clusters comprising the bad zone are typically small; and the good zone is, in essence, comprised of infinite open paths that move very directly in the drift direction.

The first of these results is stated in terms of the width of a subset $A \subseteq \Z^d$, which we now define to be 
\[
W(A)=\max_{1 \leq i \leq d}\Bigl(\max_{y\in A} y\cdot e_i - \min_{y\in A} y\cdot e_i\Bigr).
\]

\begin{lemma}
\label{BL_size_closed_box}
There exists $K_0 \in \N$ such that, for any given $K \geq K_0$, and for any $\overline{x}\in \Z^d$, the cluster $T^K(\overline{x})$ is finite $P_p[~\cdot~]$-almost surely (and ${\bf P}_p[~\cdot~]$-almost surely). Moreover, for $K \geq K_0$,
\[
P_p[  W(T^K(\overline{x})) \geq n] \leq C(K)\exp(-M(K)n),
\]
where $M(K) \to \infty$  as $K \to \infty$.
\end{lemma}

\noindent{\bf Proof.}
Consider the dependent site percolation defined by $\1{B_{\overline{x}}(L)\text{ is super-open}}$ and denote its law by $P_{p,\text{super}}$. Let us denote by $S^K(\overline{x})$ the weakly bad 
connected component of $\overline{x}$ for $P_{p,\text{super}}$.

 It is plain that, 
if $\overline{y} \in \Z^d$ is such that the $K$-box at $\overline{y}$ is super-good, and 
 $\overline{z} \in \Z^d$ satisfies $\abs{\abs{\overline{z}-\overline{y}}}_{\infty}\leq 1$, then the box at $\overline{z}$ is good. Hence, if $\overline{z}$ is bad, then every box $\overline{y}$ with $\abs{\abs{\overline{y}-\overline{z}}}_{\infty}\leq 1$ are weakly bad. This implies that the  bad *-connected component of $\overline{x}$ ($T^K(\overline{x})$) is included in $S^K(\overline{x})$, so
 \[
W(T^K(\overline{x}))\leq W(S^K(\overline{x})).
\]
Hence,  we only need to show the lemma for $W(S^K(\overline{x}))$. 

 We call two vertices 2-connected if $\abs{\abs{u-v}}_1=2$ so that we may define the even weakly bad cluster $S^K_e(\overline{x})$ of $\overline{x}$ as the 2-connected component of weakly bad boxes containing the box at $\overline{x}$. It is then clear that any element of $S^K(\overline{x})$ is a neighbor of $S^K_e(\overline{x})$ so that $W(S^K(\overline{x}))\leq W(S^K_e(\overline{x}))+2$.

Consider now the site percolation model on the even lattice $\Z^d_{\text{even}}=\{\overline{v}\in \Z^d,~\abs{\abs{\overline{v}}}_1\text{ is even} \}$ where $\overline{y}$ is even-open if and only if, in the original model, the boxes at $\overline{y}$, $\overline{y+e_1}$ and $\overline{y+e_1+e_i}$  ($i\leq d$) are super-open. An edge $[\overline{y},\overline{z}]$ is even-open if, and only if, $\overline{y}$ and $\overline{z}$ are even-open.  Note the following.
\begin{enumerate}
\item This model is a $6$-dependent oriented percolation model in the sense  that the status of edges of $\Z^d_{\text{even}}$ at distance more than six is independent. We will denote its law by   $P_{p,\text{orient}}$.
\item For $K$ large enough, by Lemma~\ref{BL_open}, the probability that a box is super-open is arbitrarily close to $1$; hence, the probability that an edge is even-open in $P_{p,\text{orient}}$ can be made arbitrarily close to $1$.
\end{enumerate}

Fix $p'$ close to 1.  By using results in~\cite{LSS}, the law $P_{p,\text{orient}}$ dominates an independent and identically distributed bond percolation with parameter $p'$ for $K$ large enough. Let us introduce the outer edge-boundary $\partial_E^e S^K_e(\overline{x})$ of $S^K_e(\overline{x})$ (in the graph $\Z_{\text{even}}^d$ with the following notion of adjacency: $\overline{x}$ and $\overline{y}$ are adjacent if $\overline{x}-\overline{y}\in \{ \pm (e_i \pm e_j),~i\neq j \leq d\}$).

We describe how to do the proof for $d=2$.  We will assign an arrow to any edge $[\overline{y},\overline{z}]\in \partial_E^e S^K_e(\overline{x})$, with $\overline{y}\in S_e^K(\overline{x})$ and $\overline{z}\notin S_e^K(\overline{x})$, in the following way:
\begin{enumerate}
\item $\swarrow$, if $\overline{y}-\overline{x}=e_1-e_2$,
\item $\nwarrow$, if $\overline{y}-\overline{x}=e_1+e_2$,
\item $\nearrow$, if $\overline{y}-\overline{x}=-e_1+e_2$,
\item $\searrow$, if $\overline{y}-\overline{x}=-e_1-e_2$.
\end{enumerate}

This boundary is represented dually in Figure 7. By an argument similar to that of~\cite{Durrett} (p.1026), we see that $n_{\nearrow} +n_{\searrow} =n_{\swarrow}+n_{\nwarrow}$, where $n_{\nearrow}$, for example,  is the number of edges labelled $\nearrow$ in $\partial_E^e S_e^K(\overline{x})$.

\begin{figure}[h]
\label{bij}
\centering
\epsfig{file=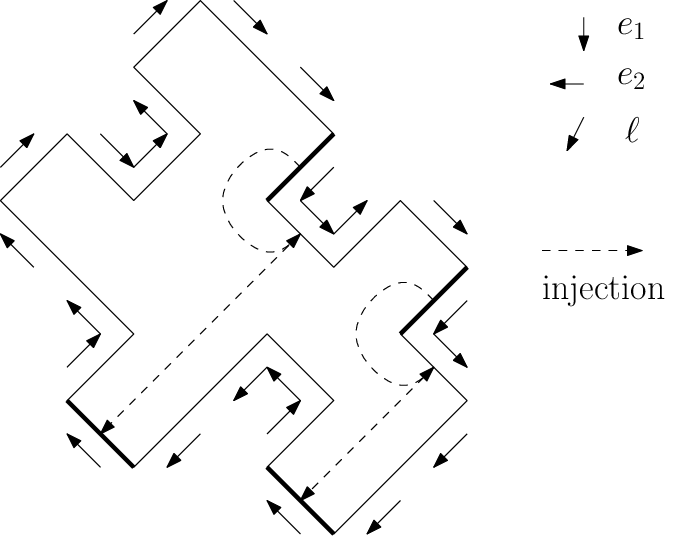, width=10cm}
\caption{The outer edge-boundary $\partial_E^e S_e^K(\overline{x})$ of $S_e^K(\overline{x})$ represented dually on the even lattice when $d=2$.}
\end{figure}

\begin{enumerate}
\item  Any $\nwarrow$ edge of $\partial_E^e  S^K_e(\overline{x})$ has one endpoint, say $\overline{y}$, which is bad, and one, $\overline{y+e_1+e_2}$, which is good. This implies that $\overline{y}$ is even-closed. Hence, if $n_{\nwarrow} \geq n_{\swarrow}/2$, then at least one sixth  of the edges of $\partial_E^e S_e^K(\overline{x})$ are even-closed, since  $n_{\nearrow} +n_{\searrow} =n_{\swarrow}+n_{\nwarrow}$.
\item Otherwise, let us assume that $n_{\swarrow}\geq 2n_{\nwarrow}$. We may notice any $\swarrow$ edge followed (in the sense of the arrows) by an $\searrow$ edge can be mapped in an injective manner to an $\nwarrow$ edge. This injection is indicated in Figure 7 (by considering the bold edges). This injection, with $n_{\swarrow}\geq 2n_{\nwarrow}$, means that at least half of the $\swarrow$ edges are not followed by an $\searrow$ edge. So, using that $n_{\nearrow} +n_{\searrow} =n_{\swarrow}+n_{\nwarrow}$, we see that at least one sixth of the edges of $\partial_E^d  S_e^K(\overline{x})$ are $\swarrow$ edges that are not followed by an $\searrow$ edge. Consider such an $\swarrow$ edge, we can see that the endpoint $\overline{y}$ of the $\swarrow$ edge which is inside $S_e^K(x)$ verifies that $\overline{y+2e_1}$ is not in $S_e^K(\overline{x})$. Hence for any such $\swarrow$, there is one endpoint $\overline{y}$ which is bad and such that $\overline{y+2e_1}$ is good, and hence $\overline{y}$ is even-closed. Once again at least one sixth of the edges of $\partial_E^e S_e^K(\overline{x})$ are even-closed.
\end{enumerate}

This means that at least one sixth of the edges of $\partial_E^e  S_e^K(\overline{x})$ are even-closed. The outer boundary is a minimal cutset, as described in~\cite{Babson}. The number of such boundaries of size $n$ is bounded (by Corollary 9 in~\cite{Babson}) by $\exp(Cn)$. Hence, if $p'$ is close enough to $1$, a counting argument allows us to obtain the desired exponential tail for $W( S_e^K(\overline{x}))$ under $P_{p'}$, and hence under $P_{p,\text{orient}}$ (since the latter is dominated by the former).

For general dimension, we note that there exists $i_0\in[2,d]$ such that a proportion at least $1/d$ of the edges of $\partial_E^e  S_e^K(\overline{x})$  are edges of the form $[\overline{y},\overline{y\pm e_{i_0}}]$ and $[\overline{y},\overline{y\pm e_{1}}]$ for some $\overline{y}\in \Z^d$. We may then apply the previous reasoning in every plane $\overline{y}+\Z e_1+\Z e_{i_0}$ containing edges of $\partial_E^e S_e^K(\overline{x})$ to show that at least a proportion $1/6$ of those edges are even-closed. Thus, at least a proportion $1/6d$ of the edges of $\partial_E^e  S_e^K(x)$ verify the same property. Repeating the same counting argument as in the previous paragraph allows us to infer the lemma.
\qed \medskip

\begin{lemma}
\label{BL_self_avoid}
Fix a vertex $x$ which lies in the largest cluster of a good box. Then there exists an infinite open simple path $x=p_x(0),\ldots,$ such that
\[
\text{for all $i$, }p_x(i)\in M_{\overline{y_i}}(K) \text{ for some $\overline{y_i}\in \Z^d$ such that } B_{\overline{y_i}}(K) \text{ is good.}
\]
and
\[
\max_{i\in \N} (x-p_x(i))\cdot \vec{\ell}\leq CK,
\]
as well as
\[
\min\{i \geq 0: (P_x(i)-x) \cdot \vec \ell \geq n\} \leq 2(n+1) \abs{B_0(K)}.
\]
\end{lemma}
\noindent{\bf Proof.}
If $x$ is in the largest cluster of a good box, there exists a directed path $\overline{z_0},\overline{z_1},\ldots$ of good boxes, where $x\in B_{\overline{z_0}}(K)$.

Now, let us consider two of these boxes that are adjacent and whose indices are $\overline{z_i}$ and $\overline{z_{i+1}}$. We can show that $M_{\overline{z_i}}(K)\cap M_{\overline{z_{i+1}}}(K)\neq \emptyset$.
Indeed, say $\overline{z_{i+1}}=\overline{z_i+e}$ for some $e\in \nu$; then $B_{\overline{z_i}}(K)\cap B_{\overline{z_{i+1}}}(K)=\{ u\in \Z^d,~ (u-K\overline{z_i}) \cdot e \in [3K/8,5K/8], ~ (u-K\overline{z_i})\cdot f \in [-5K/8,5K/8] \text{ for $f\in \nu \setminus \{e\}$} \}$. The set $M_{\overline{z_i}}(K)$ being a crossing cluster of $B_{\overline{z_i}}(K)$ in the direction $e$, it induces an open cluster of size at least $K/4$ in $B_{\overline{z_{i+1}}}(K)$, which means, by the second property of good boxes, that $M_{\overline{z_i}}(K)\cap M_{\overline{z_{i+1}}}(K)\neq \emptyset$.

Using this property, it is possible to construct an infinite simple path $x=p_x(0),\ldots$ such that
\[
\text{ if, for some $i,j$, } \, p_x(i)\in M_{\overline{z_j}}(K), \text{ then } p_x(i+1) \in \cup_{k\geq j} M_{\overline{z_k}}(K).
\]

Hence, the first property is verified.

Moreover, since the path $p_x$ is directed, this implies that 
\[
\max_{i\in \N} (x-p_x(i))\cdot \vec{\ell}\leq \max_{e\in \nu} \max_{a\in B_{\overline{0}}(K),~b\in B_{\overline{e}}(K)} (a-b)\cdot \vec{\ell} \leq 2K,
\]
which is the second desired property. 

For the third one, notice that a directed path of $2(n+1)$ $K$-boxes $B_{\overline{z_i}}(K)$ starting at $x$ will intersect $\mathcal{H}^+(x\cdot \vec{\ell}+n)$. Hence the path $(p_x(i))_{i\in\N}$ constructed above can only visit $2(n+1)\abs{B_0(K)}$ before entering $\mathcal{H}^+(x\cdot \vec{\ell}+n)$, and since the path is simple our third property is verified.
\qed \medskip

\subsubsection{The actual backbone-trap decomposition}
Recall that our aim is to prove an analogue - indeed, an improved version - of the exit-time
Proposition~\ref{BL2_PF} for the modified walk associated to the provisional backbone-trap decomposition defined at the end of Subsection~\ref{sect_def}. To do so, we must find an analogue of the Saloff-Coste argument that we gave in~(\ref{arg_SC}). For that argument to work well, we must construct, for each $x \in \mathfrak{B}_{\text{prov}}$, a semi-infinite open path emanating from $x$ that moves in the walk's preferred direction $\vec{\ell}$ with little backtracking. If $x$ lies in the largest cluster of a good box, it is Lemma~\ref{BL_self_avoid} that furnishes such a path.

It is now that we may see how our provisional choice of backbone-trap decomposition is problematic. 
Consider Figure~8.  The dotted polygon indicates the boundary of
one of the bad zones  
$T^K(\overline{x})$ (for some $\overline{x} \in \Z^d$ for which $B_{\overline{x}}(L)$ is bad).
The jagged solid lines depict a closed surface that causes the indicated vertex $v$ to lie in a finite open cluster of the graph formed by removing  $T^K(\overline{x})$  from the percolation. It may certainly be the case that $x \in K_{\infty}$, in which case, $x \in \mathfrak{B}_{\text{prov}}$. In this case, then, the need to find a viable form of the Saloff-Coste argument entails finding a forward-moving infinite open path that starts at $v$.
However, any such path must move through the bad zone  $T^K(\overline{x})$, and, in this region, we have not obtained control over the structure of long open paths.

This difficulty noticed, its solution is simple enough: we will choose a backbone-trap decomposition whose trapping zone is given by adding to $\mathfrak{T}_{\text{prov}}$  all of the finite open clusters created by the removal of $\mathfrak{T}_{\text{prov}}$ from the percolation configuration. 
There will then be no need to construct forward-going paths from points such as $v$.


\begin{figure}[h]
\label{figbing}
\centering
\epsfig{file=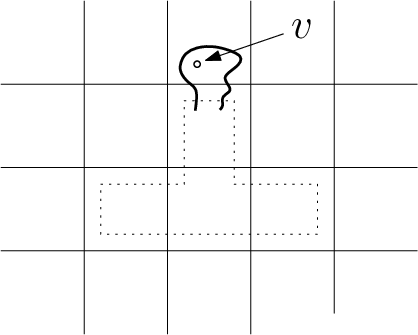, width=13cm}
\caption{Every infinite open path from the vertex $v$ passes through the adjacent bad cluster.}
\end{figure}

In the sequel, $K$ will be chosen to exceed $K_0$ in Lemma~\ref{BL_size_closed_box}, so that this lemma may be applied.

\subsubsection{Definition}

We now formally introduce the bad zone $\text{BAD}_{K}(x)$ containing $x$. In the case that all boxes containing $x$ are open, we have $\text{BAD}_{K}(x)=\emptyset$. If one of the $K$-boxes that contains $x$ is bad,  (say $B_{\overline{z}}(K)$ for some $\overline{z}\in \Z^d$),   we enumerate $\overline{x_1},\overline{x_2},\ldots, \overline{x}_{\abs{T^K(\overline{z})}}$ the points of $T^K(\overline{z})$ (see after Lemma~\ref{BL_open} above for the definition). We notice that the *-connected component $T^K(\overline{z})$ does not depend on the choice of $\overline{z}$ as long as the box at $\overline{z}$ is bad, so $T^K(\overline{z})$ is well-defined. This allows us to define
\begin{enumerate}
\item $\text{BAD}_{K}'(x)=\cup_{i \leq \abs{T^K(\overline{z})}} \{y\in B_{\overline{x_i}}(K),~ y \text{ is not in any good box}\},$
\item $\text{BAD}_{K}''(x)=\cup_{y\in H}  K^{\omega\setminus \text{BAD}_{K}'(x)}(y) $ where $H$ is the subset of points $y\in \partial \text{BAD}_{K}'(x)$ such that $y$ is connected to $\text{BAD}_{K}'(x)$ in $\omega$ and $\abs{K^{\omega\setminus \text{BAD}_{K}'(x)}(y)} <\infty $,
\item $\text{BAD}_{K}(x)=\text{BAD}_{K}'(x) \cup \text{BAD}_{K}''(x)$.
\end{enumerate}

Finally, let us define the $K$-bad zone $\text{BAD}_K=\cup_{x\in \Z} \text{BAD}_K(x)$ which is our trapping zone. Associated to this, we have the backbone:
\begin{equation}
\label{def_good}
\text{GOOD}_K=\Z^d \setminus \text{BAD}_K = \Z^d \setminus \cup_{x\in \Z} \text{BAD}_K(x).
\end{equation}

We choose to use the notation $\text{BAD}_K$ and $\text{GOOD}_K$  traditional for a percolation argument.

\subsubsection{Properties}

Firstly,  we notice the following.
\begin{lemma}
\label{remarkkk}
 If $y\notin \text{BAD}_K(x)$ and $y\sim z \in \text{BAD}_{K}(x)$ with $[y,z]$ open, then there exists an infinite open simple path in $\omega \setminus \text{BAD}_{K}(x)$ starting at $y$.
\end{lemma}

\noindent{\bf Proof.}
{\it Fact: }If $y\sim z \in \text{BAD}_{K}(x)$ with $[y,z]$ open, then $\abs{K^{\omega\setminus \text{BAD}_{K}'(x)}(y)}=\infty$.
We first see that $z\notin  \text{BAD}_K''(x)$. Indeed, if it were the case that $z\in  \text{BAD}_K''(x)$ (by definition of $\text{BAD}_K''(x)$) any $z'\sim z$ with $[z,z']$ open would have to be in $\text{BAD}_{K}'(x) \cup \text{BAD}_{K}''(x)= \text{BAD}_K(x)$, in particular $y$ would be in $ \text{BAD}_K(x)$, which is a contradiction.
 
  So necessarily $z\in \text{BAD}_K'(x)$ and since $y\notin  \text{BAD}_K''(x) \subseteq \text{BAD}_K(x)$, the definition of $ \text{BAD}_K'(x)$ implies that $\abs{K^{\omega\setminus \text{BAD}_{K}'(x)}(y)}=\infty$. 

Now that this {\it fact} has been noticed, it is obvious that 
\[
K^{\omega\setminus \text{BAD}_{K}'(x)}(y) \cap K^{\omega\setminus \text{BAD}_{K}'(x)}(u) = \emptyset,
\]
for any $u\in \partial \text{BAD}_{K}'(x)$ such that $\abs{K^{\omega\setminus \text{BAD}_{K}'(x)}(u)} <\infty $, meaning that
\[
K^{\omega\setminus \text{BAD}_{K}'(x)}(y) \cap \text{BAD}_{K}''(x) =\emptyset,
\]
so $K^{\omega\setminus \text{BAD}_{K}(x)}(y)= K^{\omega\setminus \text{BAD}_{K}'(x)}(y)$ and, thus, is infinite. This proves the lemma.
\qed \medskip

We wish to prove an analogue of Lemma~\ref{BL_size_closed_box}  for the augmented bad clusters. Indeed, we have the following result.

\begin{lemma}
\label{BL_trap_size}
For any $K$ large enough and $x\in  \Z^d$, we have that $\text{BAD}_{K}(x)$ is finite $P_p[~\cdot~]$-a.s. (and ${\bf P}_p[~\cdot~]$-a.s.),
\[
\text{BAD}_{K}(x) \subseteq \text{BAD}_{K}'(x)+[-K/8,K/8]^d,
\]
 and finally
 \[
{\bf P}_p [W(\text{BAD}_{K}(x)) \geq n] \leq C(K)e^{-M(K)n},
\]
where $M(K) \to \infty$ as $K \to \infty$.
\end{lemma}

\noindent{\bf Proof.} {\it Step 1: Proving the lemma under an assumption.}

Let us prove the result under the assumption that, for each $y\in \partial \text{BAD}_{K}'(x)$ such that $\abs{K^{\omega\setminus \text{BAD}_{K}'(x)}(y)} <\infty $, the bound $\abs{K^{\omega\setminus \text{BAD}_{K}'(x)}(y)} \leq K/8$ holds.

The inclusion
\[
\text{BAD}_{K}(x) \subseteq \text{BAD}_{K}'(x)+[-K/8,K/8]^d,
\]
is a direct consequence of that assumption. Choose a bad $K$-box at $\overline{y}$ containing $x$; if there is none, we replace $T^K(\overline{y})$ by $\emptyset$ in the sequel. We have that
\[
W(\text{BAD}_{K}'(x)) \leq \frac{5K}4 W(T^K(\overline{y})),
\]
whence
\[
W(\text{BAD}_{K}(x)) \leq  \frac{5K}4 W(T^K(\overline{y})) +K/4.
\]
Using the estimate provided by Lemma~\ref{BL_size_closed_box}, we complete the proof of Lemma~\ref{BL_trap_size}.
 
 \vspace{0.5cm}

{\it Step 2: Proof of the assumption.}

\vspace{0.5cm}

 Firstly, we see that, if $\abs{K^{\omega\setminus \text{BAD}_{K}'(x)}(y)} <\infty $ and $y\in \partial \text{BAD}_K'(x)$, then
\begin{enumerate}
\item  by definition of $\text{BAD}_{K}'(x)$,
 the cluster $K^{\omega\setminus \text{BAD}_{K}'(x)}(y) $ (which has no points of $\text{BAD}_{K}'(x)$) can only intersects good boxes, and 
  \item $K^{\omega\setminus \text{BAD}_{K}'(x)}(y)$ cannot intersect the biggest cluster of a good box, for, otherwise, Lemma~\ref{BL_self_avoid} would imply that this cluster is infinite.
\end{enumerate}

We take a $K$-box at $\overline{y}$ containing $y$. Using this, we will show that 
\[
K^{\omega\setminus \text{BAD}_{K}'(x)}(y) \subseteq \cup_{\abs{\abs{\overline{y_1}-\overline{y}}}_{\infty}\leq 1} B_{\overline{y_1}}(K).
\]

 The set $K^{\omega\setminus \text{BAD}_{K}'(x)}(y) $ does not intersect the biggest cluster of a good box, so that we know that $\abs{K^{\omega\setminus \text{BAD}_{K}'(x)}(y) \cap B_{\overline{y_1}}(K)}\leq 3^{-d}K/8$ for all $\overline{y_1}$ such that $\abs{\abs{\overline{y}-\overline{y_1}}}_{\infty}\leq 1$ and the $K$-box at $\overline{y_1}$ is good. Using the fact that $K^{\omega\setminus \text{BAD}_{K}'(x)}(y) $ intersects only good boxes, we find that
\[
\abs{K^{\omega\setminus \text{BAD}_{K}'(x)}(y) \cap \bigcup_{ \abs{\abs{\overline{y}-\overline{y_1}}}_{\infty}\leq 1} B_{\overline{y_1}}(K) } \leq K/8.
\]

However, $x$ is at distance at least $3K/4$ from any box $B_{y_1}(K)$ with $\abs{\abs{\overline{y}-\overline{y_1}}}_{\infty}\geq 2$. Hence, $K^{\omega\setminus \text{BAD}_{K}'(x)}(y) \subseteq \cup_{\abs{\abs{\overline{y}-\overline{y_1}}}_{\infty}\leq 1} B_{\overline{y_1}}(K)$. It follows that 
 \[
 \abs{K^{\omega\setminus \text{BAD}_{K}'(x)}(y)} \leq K/8,
 \]
  so that, by Step $1$, we are done. \qed \medskip

It is necessary to argue that, in the choice of backbone-trap decomposition, enough problematic areas have been excised from the backbone  that Lemma~\ref{BL_self_avoid}  is applicable to every point in the backbone. This we will show in Proposition~\ref{BL_prop_walk_1}, for which the next two results are preliminaries.

\begin{lemma}
\label{boxes}
For any $K$ large enough and $x,y \in \Z^d$, we have that
\begin{enumerate}
\item if $\text{BAD}_K(x) \cap \text{BAD}_K(y)\neq \emptyset$, then $\text{BAD}_K(x)=\text{BAD}_K(y)$,
and that
\item if $\text{BAD}_K(x) \cap \text{BAD}_K(y)=\emptyset$, then $d_{\Z^d}(\text{BAD}_K(x),\text{BAD}_K(y)) \geq K/2$.
\end{enumerate}
\end{lemma}
\noindent{\bf Proof.}
Assume that $\text{BAD}_K(x)$ and $\text{BAD}_K(y)$ are non-empty. By Lemma~\ref{BL_trap_size}, we have that
\[
\text{BAD}_{K}(x) \subseteq (KT^K(\overline{x})+ [-5K/8,5K/8]^d)+[-K/8,K/8]^d,
\]
where $\overline{x}$ is a bad $K$-box containing $x$, and
\[
 \text{BAD}_{K}(y) \subseteq  (KT^K(\overline{y})+ [-5K/8,5K/8]^d)+[-K/8,K/8]^d,
\]
where $\overline{y}$ is a bad $K$-box containing $y$.

{\it Case 1}: if for every $\overline{z_1}\in T^K(\overline{x})$ and $\overline{z_2}\in T^K(\overline{y})$, we have $\abs{\abs{\overline{z_1}-\overline{z_2}}}_{\infty}>1$, then 
\[
d_{\Z^d}^{\infty}(KT^K(\overline{x})+ [-5K/8,5K/8]^d,KT^K(\overline{y})+ [-5K/8,5K/8]^d) \geq 3K/4,
\]
where here $d_{\Z^d}^{\infty}(\cdot,\cdot)$ denotes the distance in the infinity norm, so that 
\[
d_{\Z^d}^{\infty}(\text{BAD}_K(x),\text{BAD}_K(y)) \geq 3K/4-2(K/8)=K/2.
\]

{\it Case 2}: if there exists $\overline{z}_1 \in T^K(\overline{x})$ and $\overline{z}_2 \in T^K(\overline{y})$ such that $\abs{\abs{\overline{z_1}-\overline{z_2}}}_{\infty}\leq 1$, then the $K$-boxes at $\overline{z}_1$ and $\overline{z}_2$ are bad and $*$-connected. By definition of $ T^K(\overline{x})$ and $ T^K(\overline{y})$ as $*$-connected components of bad $K$-boxes, we have that $T^K(\overline{x})=T^K(\overline{y})$, so that $\text{BAD}_{K}'(x)=\text{BAD}_{K}'(y)$ and thus $\text{BAD}_{K}(x)=\text{BAD}_{K}(y)$.

Hence, we see that either 
\[
\text{BAD}_{K}(x)=\text{BAD}_{K}(y),  \text{ or }d_{\Z^d}(\text{BAD}_K(x),\text{BAD}_K(y)) \geq K/2,
\]
 which proves the lemma. \qed \medskip

Recalling that $ \text{GOOD}_K$ was defined in~(\ref{def_good}), we have the following.
\begin{lemma}
\label{good_connection}
Take $K$ large enough. If $x\in \text{GOOD}_K\cap K_{\infty}$ there exist an open path from  $x$  to the largest cluster of a good box of length less than $K/2$.
\end{lemma}

\noindent{\bf Proof.}
From each $x\in K_{\infty}$ begins an open simple path $(\mathcal{P}(i))_{0\leq i\leq K/4}$ of length $K/4$. 

{\it Case 1}: if $\mathcal{P}$ is included in good boxes, then $\mathcal{P}$ intersects the largest cluster of a good box. If it did not, we could apply a reasoning similar to {\it step 2} of Lemma~\ref{BL_trap_size} to learn that $\abs{\mathcal{P}} \leq K/8$, which contradicts the fact that $\mathcal{P}$ is self-avoiding of length $K/4$. So in this case the lemma is true.

{\it Case 2}: Otherwise, let us denote $i_0=\min\{i\leq K/4:~ \mathcal{P}(i) \in \text{BAD}_K\}$ and fix $z\in \Z^d$ such that $x_2:=\mathcal{P}(i_0)\in\text{BAD}_K(z)$. By means of $x\in \text{GOOD}_K$, we find that $i_0\geq 1$, and so we may consider $x_1:=\mathcal{P}(i_0-1)\in  \text{GOOD}_K $. We have that $x_1\in \text{GOOD}_K$ is adjacent to  $x_2 \in\text{BAD}_K(z)$ and $[x_1,x_2]$ is open, and so, by Lemma~\ref{remarkkk}, we know that there is an infinite open simple path $\mathcal{P}'$ in $\Z^d\setminus \text{BAD}_K(z)$ starting at $x_1$. The path $\mathcal{P}'$ starts at distance one from $\text{BAD}_K(z)$, so, by the second part of Lemma~\ref{boxes}, the first $K/4$ steps of $\mathcal{P}'$ do not intersect $\text{BAD}_K(y)$ for $y\neq z$. This means that the first $K/4$ steps of $\mathcal{P}'$ do not intersect $\text{BAD}_K$. 
Using the reasoning of the {\it Case 1}, the path $\mathcal{P}'$ intersects the largest cluster of a good box. Finally, concatenating the path from $x$ to $x_1$ provided by $\mathcal{P}$ and the path from $x_1$ to the largest cluster of a good box provided by $\mathcal{P}'$, we finish the proof of the lemma. \qed \medskip

The fundamental property we will need about $ \text{GOOD}_K$ is the following.

\begin{proposition}
\label{BL_prop_walk_1}
Let $K$ be large enough. From any vertex $x\in K_{\infty}\cap \text{GOOD}_K$, there exists a nearest-neighbor (in $\Z^d$) open path $x=x_0, x_1 \ldots$  in $\text{GOOD}_K$ such that
\[
\max_{i\geq 0} (x_0-x_i)\cdot \vec \ell \leq C K,
\]
and, for any $n$,
\[
\min\{i \geq 0: (x_i-x_0) \cdot \vec \ell \geq n\} \leq (n+1) 2\abs{B_0(K)}+CK.
\]
\end{proposition}

\noindent{\bf Proof.}
By Lemma~\ref{good_connection}, $x$ is at distance less than $K/2$ from some $y$ lying in the largest cluster of a good box. Now, by Lemma~\ref{BL_self_avoid}, we can find from $y$ an infinite open path $y=p_y(0),\ldots$ such that
 \[
\max_{i\in \N} (y-p_y(i))\cdot \vec{\ell}\leq CK,
\]
which, by construction, is included in the union over the good boxes $B$ in $\omega$ of the biggest cluster in $B$; note that  this union lies in $\text{GOOD}_K$. Concatenating the path from $x$ to $y$ and this infinite path in $\text{GOOD}_K$, we have constructed an infinite open path in $\text{GOOD}_K$ starting from $x$ that verifies the properties mentioned in the proposition. \qed \medskip

Henceforth, $K$ is fixed large enough for Proposition~\ref{BL_prop_walk_1} to be verified. We will not  display the dependence of constants on $K$ anymore.

\subsection{Constructing the modified walk}

We now construct the modified walk that ignores the trapping zones.

 Let  $A\subseteq \Z^d$ be composed of nearest-neighbor components which are finite. We call the walk induced by $X_n$ on $A$, the walk $Y_n$ defined to be $Y_n=X_{\rho_n}$ where
\[
\rho_0=T_{A} \text{ and } \rho_{i+1}=T^+_{ A}\circ \theta_{\rho_i}.
\]

We now give an intrinsic construction of the walk using conductances. We denote by $\omega^*_{A}$ the graph obtained from $\omega$ by the following transformation. The vertices of $\omega^*_{A} $ are the vertices of $\Z^d\setminus A$, and the edges of $\omega^*_{A} $ are
\begin{enumerate}
\item  $\{[x,y],~x, y \in \Z^d\setminus A, \text{ with } x\notin \partial A \text{ or } y\notin \partial A\}$ and have conductance $c^{\omega^*_{A}}([x,y]) := c^{\omega}([x,y])$,
\item $ \{[x,y],~ x,y \in \partial A\}$ (including loops) which have conductance
\begin{align*}
c^{\omega^*_{A}}([x,y]) :&= \pi^{\omega}(x)P^{\omega}_x[X_1\in  A\cup \partial A,~T_y^+=T_{\partial A}^+] \\
     &= \pi^{\omega}(y)P^{\omega}_y[X_1\in  A\cup \partial A,~T_x^+=T_{\partial A}^+].
\end{align*}
\end{enumerate}
The second equality is a consequence of reversibility, and ensures symmetry for the conductances. This graph allows us to define a reversible random walk since there are only a finite number of edges with non-zero conductance at each vertex.

Now we define $\omega_K=\omega^*_{\text{BAD}_K(\omega)}$, which is possible since $\text{BAD}_K(\omega)$ has only finite connected components for $K$ large enough by Lemma~\ref{BL_trap_size} and the first part of Lemma~\ref{boxes}. We notice that the vertices of $\omega_K$ are in $\text{GOOD}_K$. The modified walk to which we have been referring is the walk in $\omega_K$.

\begin{proposition}
\label{BL_prop_walk}
The reversible walk defined by the conductances $\omega_K$ verifies the two following properties.
\begin{enumerate}
\item  It is reversible with respect to $\pi^{\omega}(\cdot)$.
\item If started at $x\in \omega_K$, it has the same law as  the walk induced by $X_n$ on $\text{GOOD}_K(\omega)$ started at $x$.
\end{enumerate}
\end{proposition}
\noindent{\bf Proof.}
To prove the first part, 
note that the invariant measure for $\pi^{\omega_K}$ satisfies by definition: for any $x\in \omega_K$,\begin{align*}
\pi^{\omega_K}(x)=&\sum_{y \text{ an } \omega_K-\text{neighbor of }x} c^{\omega_K}(x,y) \\
=&\sum_{\substack{y\in \partial \text{BAD}_K \text{ and is} \\ \text{ an } \omega_K-\text{neighbor of }x}} c^{\omega_K}(x,y) \, +
\sum_{\substack{y \notin \partial\text{BAD}_K \text{ and is}\\ \text{ an } \omega_K-\text{neighbor of }x}} c^{\omega_K}(x,y) \\
=& 
\sum_{\substack{y \in \partial\text{BAD}_K \text{ and is}\\ \text{ an } \omega_K-\text{neighbor of }x}} c^{\omega_K}(x,y) \,  +
\sum_{\substack{y\notin \partial\text{BAD}_K \text{ and is}\\ \text{ an } \omega_K-\text{neighbor of }x}} c^{\omega}(x,y).
\end{align*}

Furthermore,
\begin{align*}
& \sum_{\substack{ y \in \partial\text{BAD}_K \text{ and is} \\ \text{ an } \omega_K-\text{neighbor of }x}} c^{\omega_K}(x,y) \\
=&\sum_{\substack{y \in \partial\text{BAD}_K\text{ and is} \\ \text{ an } \omega_K-\text{neighbor of }x}}  \pi^{\omega}(x)P^{\omega}_x[X_1\in \text{BAD}_K,~T_y^+=T_{\partial \text{BAD}_K}^+]  \\
=& \pi^{\omega}(x)P^{\omega}_x[X_1\in  \text{BAD}_K]  =\sum_{\substack{y \in \partial\text{BAD}_K \text{ and is}\\ \text{ an } \omega-\text{neighbor of }x}} c^{\omega}(x,y),
\end{align*}
and since $\{y\notin \partial\text{BAD}_K \text{ and is} \text{ an } \omega-\text{neighbor of }x\}=\{y\notin \partial\text{BAD}_K \text{ and is} \\ \text{ an } \omega_K-\text{neighbor of }x\}$, we have that
\[
\pi^{\omega_K}(x)=\sum_{\substack{y\notin \partial\text{BAD}_K \text{ and is} \\ \text{ an } \omega-\text{neighbor of }x}} c^{\omega}(x,y) +
\sum_{\substack{y \in \partial\text{BAD}_K \text{ and is} \\ \text{ an } \omega-\text{neighbor of }x}} c^{\omega}(x,y) = \pi^{\omega}(x).
\]

The second point follows directly from the definition, since $\pi^{\omega_K}=\pi^{\omega}$. \qed \medskip

Let us notice the following technical detail.
\begin{lemma}
\label{cond_inegalite}
For $x,y \in \text{GOOD}_K$ that are nearest-neighbors, we have that $c^{\omega_K}([x,y])\geq c^{\omega}([x,y])$.
\end{lemma}
\noindent{\bf Proof.}
All cases are clear except when $x,y \in \partial \text{BAD}_K$. In that case, since $x\sim y$, we have that
\begin{align*}
c^{\omega_K}([x,y])&= \pi^{\omega}(x)P^{\omega}_x[X_1\in  \text{BAD}_K\cup \partial \text{BAD}_K,~T_y^+=T_{\partial \text{BAD}_K}^+]  \\
 & \geq \pi^{\omega}(x)P^{\omega}_x[X_1=y]=c^{\omega}([x,y]).
\end{align*}
\qed \medskip

\subsection{ Spectral gap estimate in $\omega_K$}


We are now ready to adapt the proof of Proposition~\ref{BL2_PF} for our modified walk. For some technical reasons, we introduce the notation
\[
\tilde{B}(L,L^{\alpha})=\{x\in \Z^d,~-L\leq x\cdot \vec{\ell} \leq 2L \text{ and } \abs{x\cdot f_i} \leq L^{\alpha} \text{ for $i\geq 2$}\}.
\]

In this context, the principal Dirichlet eigenvalue of $I-P^{\omega_K}$ in $\tilde{B}(L,L^{\alpha})\cap K_{\infty} \cap \text{GOOD}_K$ is
\begin{equation}
\label{prop_dirichlet_1}
\Lambda_{\omega_K}(\tilde{B}(L,L^{\alpha}))=\begin{cases} & \inf\{\mathcal{E}(f,f), f_{\mid (\tilde{B}(L,L^{\alpha})\cap K_{\infty} \cap \text{GOOD}_K)^c}=0, \abs{\abs{f}}_{L^2(\pi(\omega_K))}=1\}, \\ & \qquad  \text{ when $\tilde{B}(L,L^{\alpha})\cap K_{\infty} \cap \text{GOOD}_K \neq \emptyset$,} \\ 
                                  &  \infty, \text{ by convention when } \tilde{B}(L,L^{\alpha}) \cap K_{\infty} \cap \text{GOOD}_K=\emptyset , \end{cases}
                                  \end{equation}                                  
  where the Dirichlet form is defined for $f,g\in L^2(\pi^{\omega_K})$ by
  \begin{align*}
  \mathcal{E}_K(f,g)& =(f,(I-P^{\omega_K})g)_{\pi^{\omega_K}} \\
                              &=\frac 12 \sum_{x,y \text{ neighbors in }\omega_K} (f(y)-f(x))(g(y)-g(x)) c^{\omega_K}([x,y]).
  \end{align*}

We have the following lemma.
\begin{lemma}    
\label{BL_Dirichlet}
There exists $c>0$, such that for $\omega$ verifying $\tilde{B}(L,L^{\alpha})\cap K_{\infty} \cap \text{GOOD}_K \neq \emptyset$,
 \[
 \Lambda_{\omega_K}(\tilde{B}(L,L^{\alpha})) \geq c L^{-(d+1)}.
 \]
\end{lemma}
\noindent{\bf Proof.}
From any vertex $x\in K_{\infty} \cap \text{GOOD}_K$, we may choose an open nearest-neighbor path (using only edges of $\Z^d$) $p_x = (p_x(i))_{0\leq i\leq l_x}$ with the property that $p_x(0) = x$, $p_x(i) \in B(L,L^{\alpha})\cap K_{\infty} \cap \text{GOOD}_K$ for $0 < i < l_x$ and
  $p_x(l_x) \in \partial B(L,L^{\alpha})\cap \text{GOOD}_K$ verifying the condition of Proposition~\ref{BL_prop_walk_1}. Using Proposition~\ref{BL_prop_walk} and Lemma~\ref{cond_inegalite}, we see that
\begin{equation}
\label{BL_backtrack}
 \max_{i\leq l_x} \frac{ \pi_{\omega_K}(x) }{c_{\omega_K}([p_x(i+1), p_x(i)])}\leq \frac 1 {\kappa_0}  \max_{i\leq l_x} \frac{ \pi_{\omega}(x) }{\pi_{\omega}(p_x(i))}  \leq C e^{C\lambda K},
 \end{equation}
with $l_x\leq 4(L+1) \abs{B_0(K)} +CK \leq C L+C$.

We are now ready to use the argument of Saloff-Coste~\cite{Saloff-Coste} for the modified walk. 
For $\abs{\abs{f}}_{L^2(\pi^{\omega_K})}=1$ and $f_{\mid (B(L,L^{\alpha})\cap K_{\infty})^c\cap  \text{GOOD}_K}=0$, we wrtie
\begin{align*}
1&=\sum_x f^2(x)\pi_{\omega_K}(x) = \sum_x \Bigl[\sum_i f(p_x(i+1))-f(p_x(i))\Bigr]^2\pi_{\omega_K}(x) \\
  & \leq \sum_x l_x\Bigl[\sum_i (f(p_x(i+1))-f(p_x(i)))^2\Bigr]\pi_{\omega_K}(x) .
  \end{align*}

Now, by~(\ref{BL_backtrack}), we obtain 
\begin{align*}
1& \leq C e^{C\lambda K} \sum_{x,y \text{ neighbors in }\omega_K} (f(z)-f(y))^2 c_{\omega_K}([x,y]) \\
  & \qquad \qquad  \times \max_{b\in E(\Z^d)} \sum_{x\in K_{\infty} \cap \text{GOOD}_K \cap  \tilde{B}(L,L^{\alpha}),b\in p_x} l_x, 
 \end{align*}
 where $b\in p_x$ means that $b=[p_x(i),p_x(i+1)]$ for some $i$. Furthermore, note that
 \begin{enumerate}
 \item $l_x\leq  CL+C $ for any $x\in \omega_K$, and that
 \item $b=[x,y]\in \omega_K$ can only be crossed by paths \lq\lq $p_z$\rq\rq~ if $b\in E(\Z^d)$ and $z\in B_{\Z^d}(x,CL+C)$,
 \end{enumerate}
which implies (since $CL+C\leq CL$) that
 \[
 \max_b \sum_{x\in  \tilde{B}(L,L^{\alpha}),b\in p_x} l_x  \leq CL^{d+1}.
 \]

Consequently,
 \[
 1\leq C(K) L^{d+1}e^{C\lambda K} \sum_{x,y \text{ neighbors in }\omega_K} (f(z)-f(y))^2 c_{\omega_K}([x,y]).
 \]

 Now, by~(\ref{prop_dirichlet_1}),
 \[
 \Lambda_{\omega_K}(\tilde{B}(L,L^{\alpha})) \geq c L^{-(d+1)}.
 \]
\qed \medskip

\subsection{Time to exit the box for the modified walk}

We introduce, for $\epsilon_0 > 0$ small enough,
\[
A(L)=\Bigl\{\max_{x\in \Z^d \cap \tilde{B}(L,L^{\alpha})} W(\text{BAD}_{K}(x)) \leq \epsilon_0L\Bigr\},
\]
which, as a trivial consequence of Lemma~\ref{BL_trap_size}, verifies  the following bound.

\begin{lemma}
\label{prob_AL}
For any $\epsilon_0>0$ we can choose $K$ large enough such that 
\[
P_p[A(L)^c] \leq c e^{-cL}.
\]
\end{lemma}

We need some special notions of boundary for $\tilde{B}(L,L^{\alpha})$ in $\omega_K$, illustrated in Figure~9. We set
\[
\tilde{\partial} \tilde{B}(L,L^{\alpha})=\tilde{B}(L,L^{\alpha}) \cap  \Bigl[ \bigcup_{\substack{ z\in \partial \tilde{B}(L,L^{\alpha})\\ \text{BAD}_K(z)=\emptyset}} \{z\} \cup \partial  \Bigl[  \bigcup_{\substack{ z\in \partial \tilde{B}(L,L^{\alpha})\\ \text{BAD}_K(z)\neq \emptyset}} \text{BAD}_K(z) \Bigr]\Bigr];
\]
to denote  the boundary in the direction of the drift, as well as
\[
\tilde{\partial}^+ \tilde{B}(L,L^{\alpha})=\tilde{B}(L,L^{\alpha}) \cap  \Bigl[ \bigcup_{\substack{ z\in \partial \tilde{B}^+(L,L^{\alpha})\\ \text{BAD}_K(z)=\emptyset}} \{z\} \cup \partial \Bigl[ \bigcup_{\substack{ z\in \partial \tilde{B}^+(L,L^{\alpha})\\ \text{BAD}_K(z)\neq \emptyset}} \text{BAD}_K(z) \Bigr] \Bigr];
\]
we also set
\[
\tilde{\partial}^- \tilde{B}(L,L^{\alpha})=\tilde{B}(L,L^{\alpha}) \cap  \Bigl[ \bigcup_{\substack{ z\in \partial \tilde{B}^-(L,L^{\alpha})\\ \text{BAD}_K(z)=\emptyset}} \{z\} \cup  \partial \Bigl[  \bigcup_{\substack{ z\in \partial \tilde{B}^-(L,L^{\alpha})\\ \text{BAD}_K(z)\neq \emptyset}} \text{BAD}_K(z) \Bigr] \Bigr],
\]
as well as
\[
\tilde{\partial}^{\text{s}} \tilde{B}(L,L^{\alpha})=\tilde{\partial} \tilde{B}(L,L^{\alpha}) \setminus \bigl(\tilde{\partial}^+ \tilde{B}(L,L^{\alpha}) \cup \tilde{\partial}^- \tilde{B}(L,L^{\alpha})\bigr).
\]

Note that 
all these sets are included in $ \text{GOOD}_K$. 

\begin{figure}[h]\label{modwalkpic}
\centering
\epsfig{file=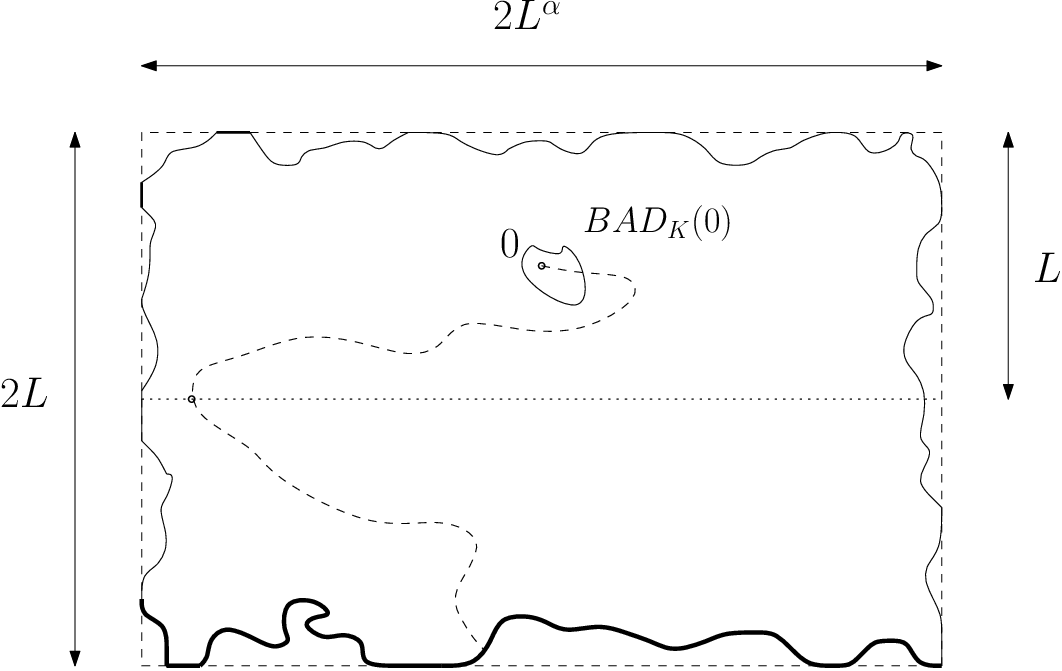, width=13cm}
\caption{The picture illustrates how, if the modified random walk reaches the bold border, i.e.~$\tilde{\partial}^+ \tilde{B}(L,L^{\alpha})$, before the rest of the border $\tilde{\partial}^- \tilde{B}(L,L^{\alpha}) \cup \tilde{\partial}^{\text{s}} \tilde{B}(L,L^{\alpha})$,  then the normal walk  exits $B(L,L^{\alpha})$ on $\partial^+ B(L,L^{\alpha})$.}
\end{figure}

\begin{proposition}
\label{BL_time}
There exists $C\in (0,\infty)$ such that, if $ K_{\infty}\cap \text{GOOD}_K \cap \tilde{B}(L,L^{\alpha})\neq \emptyset$, then
\[
\max_{x\in  K_{\infty}\cap \text{GOOD}_K  \cap B_{\Z^d}^{\infty}(0, \epsilon_0L)} P_x^{\omega_K}[T_{\tilde{\partial}\tilde{B}(L,L^{\alpha})}> C L^{d+2} ]\leq C \exp(-L),
\]
with the convention that the maximum over the empty-set is $0$ and $B_{\Z^d}^{\infty}$ denotes the ball for the infinity norm.
\end{proposition}

It will be necessary to consider the exit time of the modified walk starting from the border of the trapping zone at $0$. By Lemma~\ref{prob_AL}, this region is typically included in $K_{\infty}\cap \text{GOOD}_K  \cap B_{\Z^d}^{\infty}(0, \epsilon_0L)$, and this explains our convention to ignore the case when this set is empty.

\noindent{\bf Proof.}
Using $K_{\infty}\cap \text{GOOD}_K \cap \tilde{B}(L,L^{\alpha})\neq \emptyset$, we find from  Perron-Frobenius' theorem that the operator 
\begin{align*}
 & \tilde{P}^{\omega_K}_{\text{GOOD}}(L,\alpha)\\
:=& \1{\tilde{B}(L,L^{\alpha}) \cap  K_{\infty}\cap \text{GOOD}_K } P^{\omega_K} \1{\tilde{B}(L,L^{\alpha}) \cap  K_{\infty}\cap \text{GOOD}_K } 
\end{align*}
 has a norm on $L^2(\pi_{\omega_K})$ given by its maximum positive eigenvalue which, is $1-\Lambda_{\omega_K}(\tilde{B}(L,L^{\alpha}) )$. This yields that, for any $x\in  K_{\infty}\cap \text{GOOD}_K \cap \tilde{B}(L,L^{\alpha})$,
\begin{align*}
 & \pi_{\omega_K}(x)P_x^{\omega_K}[T_{\tilde{\partial}\tilde{B}(L,L^{\alpha})}>n]  \\
 = & (\1{x}, (\tilde{P}^{\omega_K}_{\text{GOOD}}(L,\alpha))^n  \1{\tilde{B}(L,L^{\alpha})\cap \text{GOOD}_K })_{L^2(\pi_{\omega_K})} \\  
 \leq & (\pi_{\omega_K}(x) \pi_{\omega_K}({\tilde{B}(L,L^{\alpha})}\cap \text{GOOD}_K ))^{1/2}(1-\Lambda_{\omega_K}(\tilde{B}(L,L^{\alpha}) ))^{n}.
\end{align*}

The first property of Proposition~\ref{BL_prop_walk} then implies that
\begin{align*}
&\pi_{\omega_K}(x) \pi_{\omega_K}({\tilde{B}(L,L^{\alpha})}\cap \text{GOOD}_K ) \\
\leq &\pi_{\omega}(x) \pi_{\omega}({\tilde{B}(L,L^{\alpha})} ) \leq  C L^{1+(d-1)\alpha}e^{4\lambda (1+\epsilon_0)L}.
\end{align*}

Hence, by Lemma~\ref{BL_Dirichlet},
\begin{align*}
& \max_{x\in K_{\infty}\cap \text{GOOD}_K  \cap B_{\Z^d}^{\infty}(0, \epsilon_0L))} P_x^{\omega_K}[T_{\tilde{\partial}\tilde{B}(L,L^{\alpha})}>n] \\ 
\leq & cL^{(1+(d-1)\alpha)/2}\exp(2\lambda (1+\epsilon_0) L-c n  L^{-(d+1)}),
\end{align*}
and the result follows.
\qed \medskip

\subsection{Linking exit probabilities in $\omega$ and $\omega_K$}


Let us now prove: 
\begin{lemma}
\label{qwerty}
For $\epsilon_0<1/4$, we have
\begin{align*}
& \PR\Bigl [A(L),~X_{T^{\text{ex}}_{ B(L,L^{\alpha})}} \notin \partial^+ B(L,L^{\alpha})\Bigr]\\
\leq &  C E_p \Bigl[A(L),~ \max_{x\in K_{\infty} \cap \text{GOOD}_K  \cap B_{\Z^d}^{\infty}(0, \epsilon_0L) } P_x^{\omega_K}[X_{T_{\tilde{\partial} \tilde{B}(L,L^{\alpha})}} \notin \tilde{\partial}^+ \tilde{B}(L,L^{\alpha})]\Bigr],
\end{align*}
with the convention that the last maximum is $0$ if $K_{\infty} \cap \text{GOOD}_K  \cap B_{\Z^d}^{\infty}(0, \epsilon_0L) =\emptyset$.
\end{lemma}

\noindent{\bf Proof.} 

{\it Step 1: Entering $\text{GOOD}_K$.}

\vspace{0.5cm}

 By the second part of Lemma~\ref{boxes}, we see that $T_{\text{GOOD}_K}=T_{\text{BAD}_K(0)}^{\text{ex}}$ under $P^{\omega}_0$. On $\mathcal{I}\cap A(L)$, we have that $\text{BAD}_K(0)\subseteq B_{\Z^d}^{\infty}(0, \epsilon_0L)$, so that, under $P^{\omega}_0$,
 \[
 T_{\text{GOOD}_K} \leq T_{B_{\Z^d}^{\infty}(0, \epsilon_0L)}^{\text{ex}}<T_{\partial B(L,L^{\alpha})}.
 \]
 
  Hence, on $\mathcal{I}\cap A(L)$, we have that, by Markov's property,
\begin{align}
\label{BL_stepi}
&P^{\omega}_0[X_{T_{\partial B(L,L^{\alpha})}} \notin \partial^+ B(L,L^{\alpha})]\\ \nonumber
 \leq & \max_{x\in K_{\infty}\cap \text{GOOD}_K \cap B_{\Z^d}^{\infty}(0, \epsilon_0L)} P^{\omega}_x [X_{T_{\partial B(L,L^{\alpha})}} \notin \partial^+ B(L,L^{\alpha})],
\end{align}
where we used $X_{T_{\omega_K}}\in  K_{\infty}\cap \text{GOOD}_K \cap B_{\Z^d}^{\infty}(0, \epsilon_0L)$.

\vspace{0.5cm}

{\it Step 2: Linking the exit distributions of the regular and the modified walk in $\text{GOOD}_K$.}

\vspace{0.5cm}

 We want to prove that, for any $x\in K_{\infty}\cap \text{GOOD}_K \cap B_{\Z^d}^{\infty}(0, \epsilon_0L)$, on $A(L)$,
\begin{equation}
\label{BL_ineg}
 P^{\omega}_x [X_{T_{\partial B(L,L^{\alpha})}} \notin \partial^+ B(L,L^{\alpha})]\leq P^{\omega_K}_x [X_{T_{\tilde{\partial}^+ \tilde{B}(L,L^{\alpha})}} \notin \tilde{\partial} \tilde{B}(L,L^{\alpha})].
\end{equation}

Indeed, on $A(L)$, we have $B_{\Z^d}^{\infty}(0, \epsilon_0L)\cap \tilde{\partial} \tilde{B}(L,L^{\alpha})=\emptyset$ since $\epsilon_0<1/2$ and $\tilde{\partial} \tilde{B}(L,L^{\alpha})\subset \cup_{x\in \tilde{B}(L,L^{\alpha})} B_{\Z^d}(x,\epsilon_0L)$. Hence, $T_{\tilde{\partial} \tilde{B}(L,L^{\alpha})}$ is finite under $P_x^{\omega_K}$. Trivially, 
\[
 P^{\omega_K}_x [X_{T_{\tilde{\partial} \tilde{B}(L,L^{\alpha})}} \in \tilde{\partial^+} \tilde{B}(L,L^{\alpha})]  =  P^{\omega_K}_x [T_{\tilde{\partial} \tilde{B}(L,L^{\alpha})} =T_{\tilde{\partial}^+ \tilde{B}(L,L^{\alpha})} ].
 \]
 
Moreover, the vertex $x$ and the sets $\tilde{\partial} \tilde{B}(L,L^{\alpha})$ and $\tilde{\partial^+} \tilde{B}(L,L^{\alpha})$ are in $\text{GOOD}_K$, so that, by the second property of Proposition~\ref{BL_prop_walk},
\[
 P^{\omega_K}_x [T_{\tilde{\partial} \tilde{B}(L,L^{\alpha})} = T_{\tilde{\partial^+} \tilde{B}(L,L^{\alpha})} ] =P^{\omega}_x [T_{\tilde{\partial} \tilde{B}(L,L^{\alpha})} = T_{\tilde{\partial^+} \tilde{B}(L,L^{\alpha})} ].
 \]
 
Furthermore, on $A(L)$, we have that (see Figure~8)
\[
\{T_{\tilde{\partial} \tilde{B}(L,L^{\alpha})}=T_{\tilde{\partial^+} \tilde{B}(L,L^{\alpha})}\} \subseteq \{T_{B (L,L^{\alpha})} = T_{ \partial^+ B(L,L^{\alpha})}\} \qquad P_x^{\omega}\text{-a.s.,}
\] 
 so that three last equation imply
\begin{align*}
P^{\omega_K}_x [T_{\tilde{\partial} \tilde{B}(L,L^{\alpha})} = T_{\tilde{\partial^+} \tilde{B}(L,L^{\alpha})} ] & =  P^{\omega}_x [T_{\tilde{\partial} \tilde{B}(L,L^{\alpha})} = T_{\tilde{\partial^+} \tilde{B}(L,L^{\alpha})}] \\
      & \leq P^{\omega}_x [T_{B (L,L^{\alpha})}= T_{ \partial^+ B(L,L^{\alpha})}].
 \end{align*}
 
This proves~(\ref{BL_ineg}).

\vspace{0.5cm}

{\it Step 3: Conclusion.}

\vspace{0.5cm}

Hence, by~(\ref{BL_stepi}) and~(\ref{BL_ineg}),
\begin{align*}
& \PR\Bigl [A(L),~X_{T_{\tilde{\partial} \tilde{B}(L,L^{\alpha})}} \notin \partial^+ B(L,L^{\alpha})\Bigr]\\
\leq &   \PR \Bigl[A(L),~ \max_{x\in  K_{\infty}\cap \text{GOOD}_K  \cap B(0, \epsilon_0L) } P_x^{\omega_K}[X_{T_{\tilde{\partial} \tilde{B}(L,L^{\alpha})}} \notin \tilde{\partial}^+ \tilde{B}(L,L^{\alpha})]\Bigr] \\
\leq &  \frac 1{ P_p[\mathcal{I}]}  E_p \Bigl[A(L),~ \max_{x\in  K_{\infty}\cap \text{GOOD}_K  \cap B(0, \epsilon_0L) } P_x^{\omega_K}[X_{T_{\tilde{\partial} \tilde{B}(L,L^{\alpha})}} \notin \tilde{\partial}^+ \tilde{B}(L,L^{\alpha})]\Bigr].
\end{align*}
\qed \medskip

\subsection{Proof of Theorem~\ref{BL} }      
 
We may now apply the reasoning (\ref{eqheur}) indicated at the outset, in the case of the modified walk. Before turning to the details, note that Proposition~\ref{BL_time} shows that we may in effect choose $n = L^{d+2}$ in the estimate (\ref{eqheur}). The origin being at distance at least $L^{\alpha}$ from the lateral sides of $B(L,L^\alpha)$, it is simply impossible that the original walk exit this box on those sides in time $n$ (and, as we will argue, that remains true for the modified walk, for which we work with chemical distance in the modified graph). That is, the modified walk moves so quickly that we need only concern ourselves with the possibility of exit on the top-side of the box, a matter that is easily handled by various means, of which we choose the Carne-Varopoulos bound since it is already at hand.

 Using Lemma~\ref{prob_AL} and Lemma~\ref{qwerty}, we find that, for each $n \in \N$,
\begin{align}
\label{BL_sum3}
& \PR[X_{T^{\text{ex}}_{B(L,L^\alpha)}}\notin \partial^+ B(L,L^\alpha)] \\\nonumber
\leq &  \PR\Bigl [A(L),~X_{T^{\text{ex}}_{ B(L,L^{\alpha})}} \notin \partial^+ B(L,L^{\alpha})\Bigr]   + CP_p [A(L)^c] \\\nonumber
\leq  & C E_p \Bigl[A(L),~ \max_{x\in K_{\infty} \cap \text{GOOD}_K  \cap B_{\Z^d}^{\infty}(0, \epsilon_0L) } P_x^{\omega_K}[X_{T^{\text{ex}}_{ \tilde{B}(L,L^{\alpha})}} \notin \tilde{\partial}^+ \tilde{B}(L,L^{\alpha})]\Bigr] +ce^{-cL},
\end{align}
where we recall that the maximum is non-zero only if $K_{\infty} \cap \text{GOOD}_K  \cap B_{\Z^d}^{\infty}(0, \epsilon_0L)  \neq \emptyset$, which implies that $K_{\infty} \cap \text{GOOD}_K \cap \tilde{B}(L,L^{\alpha})\neq \emptyset$.

Setting $n=\lfloor C L^{d+2} \rfloor$, Proposition~\ref{BL_time} implies that
\begin{align}
\label{BL_step2}
&E_p \Bigl[A(L),~ \max_{x\in  K_{\infty}\cap \text{GOOD}_K  \cap B(0, \epsilon_0L) } P_x^{\omega_K}[X_{T^{\text{ex}}_{ \tilde{B}(L,L^{\alpha})}} \notin \tilde{\partial}^+ \tilde{B}(L,L^{\alpha})]\Bigr], \\\nonumber
\leq & E_p \Bigl[A(L),~ \max_{x\in  K_{\infty}\cap \text{GOOD}_K  \cap B(0, \epsilon_0L) } P_x^{\omega_K}[X_{T_{\tilde{\partial}\tilde{B}(L,L^{\alpha})}} \notin \tilde{\partial}^+ \tilde{B}(L,L^{\alpha}), ~T_{\tilde{\partial}\tilde{B}(L,L^{\alpha})}\leq n ]\Bigr]  \\ \nonumber
   & \qquad \qquad \qquad \qquad \qquad \qquad \qquad \qquad + Ce^{-cL}.
\end{align}

Furthermore, for $x,y\in  K_{\infty}\cap \text{GOOD}_K $, the Carne-Varopoulos estimate~(\ref{CVB}) yields
\[
P_x^{\omega_K}[X_n=y] \leq 2 \Bigl(\frac{\pi_{\omega}(y)}{\pi_{\omega}(x)}\Bigr)^{1/2} \exp\Bigl(-\frac{d_{\omega_K}(x,y)^2}{2n}\Bigr),
\]
 where we used the first property of Proposition~\ref{BL_prop_walk}.

On $A(L)$, we have that
\[
d_{\omega_K}(x,y)\geq \lfloor d_{\omega}(x,y)/(2d\epsilon_0L) \rfloor,
\]
and
\[
 \tilde{\partial}^{\text{s}}\tilde{B}(L,L^{\alpha}) \subseteq \bigcup_{x \in \partial^{\text{s}} B(L,L^{\alpha})} B_{\Z^d}^{\infty}(x, \epsilon_0L).
\] 

This implies that for any $x\in  K_{\infty}\cap \text{GOOD}_K \cap B_{\Z^d}^{\infty}(0, \epsilon_0L)$, the distance (in $\omega_K$) from $x$ to any point of  $ \tilde{\partial}^{\text{s}}\tilde{B}(L,L^{\alpha})$ is  at least $L^{\alpha-1}-2\epsilon_0L$. In particular for $\alpha>d+3$,  the walk cannot reach $ \tilde{\partial}^{\text{s}}\tilde{B}(L,L^{\alpha})$  in time less than $n= \lfloor C L^{d+2} \rfloor$ for $L$ large enough.

Fix $\epsilon_0 <1/4$,  for any $x\in  K_{\infty}\cap \text{GOOD}_K \cap B_{\Z^d}^{\infty}(0, \epsilon_0L)$ the Carne-Varopoulos estimate yields that, on $A(L)$,
\begin{align*}
& P_x^{\omega_K}[X_{T^{\text{ex}}_{ \tilde{B}(L,L^{\alpha})}} \notin \tilde{\partial}^+ \tilde{B}(L,L^{\alpha}), ~T_{\tilde{\partial}\tilde{B}(L,L^{\alpha})}\leq n ]  \\
\leq & P_x^{\omega_K}[X_{T^{\text{ex}}_{ \tilde{B}(L,L^{\alpha})}} \in \tilde{\partial}^- \tilde{B}(L,L^{\alpha}), ~T_{\tilde{\partial}\tilde{B}(L,L^{\alpha})}\leq n ]  \\
\leq & C(p,\lambda) n \abs{\tilde{\partial}^-\tilde{B}(L,L^{\alpha})}  \max_{\substack{y\in\tilde{\partial}^-\tilde{B}(L,L^{\alpha}), \\ x\in B_{\Z^d}^{\infty}(0, \epsilon_0L)}  } e^{\lambda (y-x)\cdot \vec{\ell}}\\
\leq &  CnL^{C\alpha}e^{-\lambda L/2},
\end{align*}
where we used that, on $A(L)$,
\[
\tilde{\partial}^-\tilde{B}(L,L^{\alpha}) \subseteq \bigcup_{x \in \partial^- B(L,L^{\alpha})} B_{\Z^d}^{\infty}(x, \epsilon_0L).
\]

Hence, for $\alpha>d+3$, we have that, on $A(L)$, 
\begin{equation}
\label{BL_CV}
\max_{x\in K_{\infty}\cap \text{GOOD}_K \cap B_{\Z^d}^{\infty}(0, \epsilon_0L) } P_x^{\omega_K}[X_{T_{\tilde{\partial}\tilde{B}(L,L^{\alpha})}} \notin \tilde{\partial}^+ \tilde{B}(L,L^{\alpha}), ~T_{\tilde{\partial}\tilde{B}(L,L^{\alpha})}<C\lfloor L^{d+2} \rfloor ] \leq e^{-cL}.
\end{equation}

Taking $\alpha>d+3$ and $n= C\lfloor L^{d+2} \rfloor$, the equations~(\ref{BL_sum3}),~(\ref{BL_step2}) and~(\ref{BL_CV}) are verified, so that
\[
 \PR[X_{T^{\text{ex}}_{B(L,L^\alpha)}}\notin \partial^+ B(L,L^\alpha)]  \leq e^{-cL} .
  \]

 Hence, we obtain 
Theorem~\ref{BL}.
\qed \medskip

		\section{Percolation estimates}
\label{sect_perco}

 In this section, we provide the percolation estimates required in the paper, notably for understanding trap geometry. The principal result proved is
  Proposition~\ref{backtrack_exponent}. As a direct consequence, we will obtain Lemma~\ref{estim_height}. 
As we discussed in the introduction, trap geometry differs between the two- and higher-dimensional cases.
This difference is reflected in the style of the proofs in this section, which are more subtle in the case that $d \geq 3$.  
  The section is split into two parts: Section~\ref{D2} treats $d=2$ and Section~\ref{D3}, $d=3$. 
  

We introduce the notation $f(n) \approx g(n)$ to indicate that $\log f(n) /\log g(n) \to 1$ as $n \to \infty$. Recall also that we denote the component-wise integer-part of a vector $u\in \R^d$ by $\lfloor u \rfloor$.

For $x\in \Z^d$ we will write $K^{\omega}_{n}(x)$ for the open cluster of $x$ in $\mathcal{H}^{+}(x\cdot \vec{\ell}-n)$. We see then that
\begin{equation}
\label{waka_waka}
\mathcal{B}\mathcal{K}(x) >n \text{ if and only if } \abs{K^{\omega}_{n}(x)}<\infty \text{ and } \mathcal{I}_x.
\end{equation}

\subsection{Case of $\Z^2$.}\label{D2}

The dual lattice $\Z^2_*$ is formed by translating the nearest-neighbor lattice $\Z^2$
by $(1/2,1/2)$. 
To a percolation configuration 
$\omega \in \{0,1\}^{E(\Z^2)}$ 
is associated a dual configuration 
$\omega' \in \{0,1\}^{E(\Z^2_*)}$: for each $e' \in E(\Z^2_*)$, there is an element $e \in E(\Z^2)$
that meets $e'$ at the common midpoint of the two edges. We set $\omega'(e') = 1 - \omega(e)$.
In this way, the law $P_p$ on $E(\Z^2)$ is coupled with the law $P_{1-p}$ on $E(\Z^2_*)$.

As we will see, the event $\backtrack(0) > n$ is essentially characterized by the presence of two disjoint subcritical open dual paths that meet at a dual vertex underneath $0$ (which is to say, further in the direction of the bias than~$0$), each path reaching  the half-space $\mathcal{H}^-(-n)$; these two paths then create a trap around $0$. We are interested, then, in the probability of such subcritical dual connections. It follows from a standard subadditivity argument that for any $\vec u \in S^1$,
there exists $\xi^{(\vec{u})}\in (0,\infty)$ such that
\begin{equation}
\label{cor_length}
P_{1-p}[0 \leftrightarrow \lfloor -n\vec{u} \rfloor ] \approx \exp(-\xi^{(\vec{u})} n),
\end{equation}
a result that may also be recovered from the more detailed assertion given in Proposition 6.47 (alongside Theorem 6.15) of~\cite{Grimmett}.

We need further information than that offered by the previous assertion in order to obtain estimates on ${\bf P}_p[\mathcal{B}\mathcal{K}(0)>n]$: indeed, we must understand the typical large-scale geometry of traps.
The article \cite{civ} develops Ornstein-Zernike theory to strengthen~(\ref{cor_length}), with Theorem A of \cite{civ}  providing 
a sharp asymptotic estimate of the probability on the left-hand side of~(\ref{cor_length}). This theory also elucidates details on the geometry of traps which we need, and we now heuristically explain its key ideas and recall the results we will use in the rest of this paper.

\subsubsection{Ornstein-Zernike theory}

Fix a subcritical parameter $q \in (0, 1/2)$, and let $\vec{u} \in S^1$ (although the results we summarise here apply equally in higher dimensions).  In the Ornstein-Zernike approach, a subset of
 \lq\lq cut points\rq\rq  of a finite open cluster is identified, with the property that the removal of a cutpoint splits the cluster into two connected pieces, each one lying 
 in one of the two half-spaces whose boundary is the hyperplane containing the cutpoint and orthogonal to $\vec u$.
  See Figure~10.

\begin{figure}\label{piccutpoit}
\centering\epsfig{file=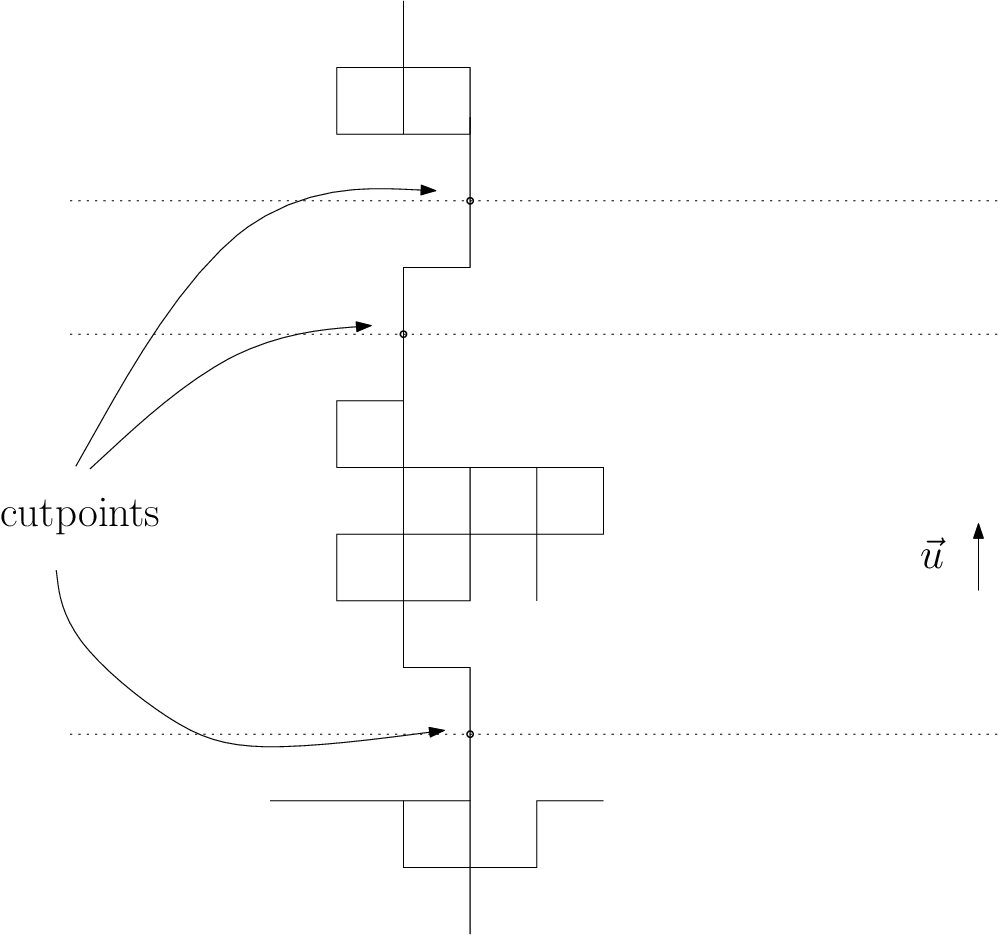, width=10cm}
\caption{The cutpoints splitting the subcritical cluster into irreducible components.}
\end{figure}

The removal of all cutpoints from the cluster results in a decomposition of the cluster into a succession of \lq\lq irreducible\rq\rq  components.
Under the law $P_q$, the successive irreducible components of the open cluster containing $0$ form an independent and identically distributed sequence. The probability that this cluster travels a high distance $n$ in the direction $\vec u$ may then be analysed by the renewal theorem; in effect, this is the method used in \cite{civ} to sharpen~(\ref{cor_length}).

The Ornstein-Zernike theory also gives a rather precise description of the geometry of the open cluster $K(0)$ containing $0$
under the conditional law $P_q \big( \cdot \big\vert 0 \leftrightarrow \lfloor n \vec u \rfloor   \big)$. Its irreducible components are small, the largest among them typically being of size of order $\log n$. More precisely, defining the spine $L = L_n$ of the cluster to be the union of the line segments that link the successive cutpoints, we have the following consequence 
(see (1.10) of \cite{civ}): for each $\epsilon > 0$, there exists $C > 0$ such that, for $n$ sufficiently high,
\begin{equation}\label{civoneten}
 P_q \Big( \textrm{$C \log n$-neighborhood of $L$ contains $K(0)$} \, \Big\vert \,  0 \leftrightarrow \lfloor n \vec u \rfloor  \Big) \geq 1 - \epsilon.
\end{equation}

Thus the geometry of $K(0)$ is intimately related to the geometry of the spine, which resembles asymptotically  a Brownian bridge. More precisely, Theorem $C$ of \cite{civ} derives a functional central limit theorem for the displacement of $L$. We will use the following consequence:  for each $\epsilon > 0$, there exists $\delta > 0$ such that, for $n$ sufficiently high,
\begin{equation}\label{civclt}
 P_q \Big(  \sup_{x \in L} d \Big( x, \big[ 0, \lfloor n \vec{u} \rfloor \big] \Big) \leq \epsilon n^{1/2} \,
  \Big\vert \, 0 \leftrightarrow \lfloor n u \rfloor 
 \Big) \geq \delta.  
\end{equation}

These two inequalities contain the information that we need to pursue the analysis of the two-dimensional case.

\subsubsection{Results for $d=2$}

Before stating our main result, we notice that the next lemma follows directly from~(\ref{cor_length}) and a union bound using the exponential decay for the size of subcritical connections (see~\cite{aiznew}).
\begin{lemma}
\label{cor_lengthnew}
Let $d = 2$ and $p > p_c$. We have that
\[
P_{1-p}[0 \leftrightarrow \mathcal{H}^{-}(n)] \approx \exp(-\hat{\xi}^{(\vec{\ell})} n),
\]
where $\hat{\xi}^{(\vec\ell)}\in (0,\infty)$ is given by
$\hat{\xi}^{(\vec\ell)} = \inf_{\vec{\ell'} \in S^1} \frac{\xi^{(\vec{\ell'})}}{\vec{\ell'} \cdot \vec\ell}$.
\end{lemma}

Our main result for  two-dimensional case is the following.
\begin{proposition}
\label{expo_2D}
Let $d = 2$ and $p > p_c$. We have that
\[
{\bf P}_{p}[\mathcal{B}\mathcal{K}(0) >n] \approx  \exp(-2 \hat{\xi}^{(\vec{\ell})}n),
\]
where $\hat{\xi}^{(\vec{\ell})} \in (0,\infty)$ is defined in Lemma~\ref{cor_lengthnew}.
\end{proposition}

Let $\vec\ell' \in S^1$ be such that the infimum   
$\inf_{\vec{u} \in S^1} \frac{\xi^{(\vec{u})}}{\vec{u} \cdot \vec\ell}$ 
is attained by the choice $\vec{u} = \vec\ell'$. 
For $x \in \Z^2$, let $C_x(n)$ denote the parallelogram with extremal vertices $x - n^{1/2} \vec\ell^\perp, x + n^{1/2} \vec\ell^\perp$ and the translates of these points by $- n \vec{\ell'}/(\vec{\ell}' \cdot \vec\ell)$. Call the interval   $\big[ x - n^{1/2} \vec\ell^\perp, x + n^{1/2} \vec\ell^\perp \big]$ the lower side of $C_x(n)$ and the interval between the other pair of extremal vertices the upper side. (The lower side should be considered to be at the level the base of a trap, and the upper side at the level of the trap's top, so that this orientation reflects our convention that the bias direction is downwards.) Two instances of $C_x(n)$ are depicted in the upcoming Figure 12  (on the left and the right); the lower side in each case is drawn in bold,  the other sides represented by dotted lines.

The next lemma will be used in the proof of Proposition~\ref{expo_2D}. 
\begin{lemma}\label{lemaxn} 
Let $x \in \Z^2$. Then 
\begin{eqnarray}
 & & P_{1-p} \Big( \textrm{an open path in $C_x(n)$ connects the lower and upper sides of $C_x(n)$} \Big) \nonumber \\
 & \approx & \exp(-\hat\xi^{(\vec{\ell})}n). \nonumber
\end{eqnarray}
\end{lemma}
\noindent{\bf Proof.} 
By translation invariance, we may assume that $x=0$. 
We argue that the conditions 
\begin{enumerate}
 \item  $0  \leftrightarrow  - \lfloor n \big( \vec\ell \cdot \vec\ell' \big)^{-1} \vec{\ell'} \rfloor$ ; 
 \item $\max_{y \in K^\omega_{n} (x)} \big\vert (y - x) \cdot \vec{\ell'}^\perp  \big\vert \leq \tfrac{1}{2} \vec{\ell'} \cdot \vec\ell \, n^{1/2}$.
\end{enumerate}
imply the event whose probability we seek to bound. Indeed, an open path realizing $0  \leftrightarrow  - \lfloor n \big( \vec\ell \cdot \vec\ell' \big)^{-1} \vec{\ell'} \rfloor$ (with, if necessary, a further open edge incident to  $- \lfloor n \big( \vec\ell \cdot \vec\ell' \big)^{-1} \vec{\ell'} \rfloor$) ensures that the upper and lower sides of $C_x(n)$ each intersect two edges lying in a common open cluster; note then that (2) ensures that this cluster remains in the strip whose boundaries contain the two $\vec\ell'$-oriented sides of $C_x(n)$. This forces an open path in $C_x(n)$ from this parallelogram's lower to its upper side.

The probability of property (1) is estimated in (\ref{cor_length}), with $n \vec{u}$ replaced by  $n \big( \vec\ell \cdot \vec\ell' \big)^{-1} \vec{\ell'}$.
For the lower bound in the lemma's statement, it is enough then to argue that, conditionally on $P \big( \cdot \big\vert \, 0  \leftrightarrow     
   - \lfloor n \big( \vec\ell \cdot \vec\ell' \big)^{-1} \vec{\ell'} \rfloor  \big)$, property (2) has a uniformly positive probability.
This follows from~(\ref{civoneten}) and~(\ref{civclt}).

The upper bound in the lemma's statement follows from Lemma~\ref{cor_lengthnew}, whence the result. \qed \medskip

The next result is also an ingredient for the proof of (\ref{eqexpolb}). It is a generalization of Theorem 11.55 in~\cite{Grimmett} to the case of non-axial orientation, and it has the same proof.
\begin{lemma}\label{lemopeninh}
Let $d=2$ and $p > p_c$.
Let $H \subseteq \Z^2$ be given by $T = \big\{ x \in \Z^2: x \cdot \vec{\ell} \geq 0 \, , \vert x  \cdot \vec{\ell}^{\perp} \vert \leq (x\cdot \vec{\ell}) ^{1/4} +3\big\}$.
There is a positive probability that $0$ lies in an infinite open simple  path, all of whose vertices lie in $T$.
\end{lemma}

\noindent{\bf Proofs of Proposition~\ref{expo_2D} and Lemma~\ref{estim_height} in the case $d=2$.}

\vspace{0.5cm}

{\it Step 1: Upper bound in Proposition~\ref{expo_2D}.}

\vspace{0.5cm}

We show that,
 for $\epsilon > 0$, and for $n \in \N$ sufficiently high,
\begin{equation}\label{eqexpoub}
{\bf P}_{p}[\mathcal{B}\mathcal{K}(0) >n] \leq  \exp\big\{ - \big( 2 \hat{\xi}^{(\vec{\ell})} + \epsilon \big) n \big\}.
\end{equation}

We find an unusual event that necessarily occurs in the case that $\backtrack(0) > n$. 
Recall from~(\ref{waka_waka}) that the condition $\backtrack(0) > n$ implies that $K^{\omega}_n(0)$ is finite. 

To any finite open cluster in $\Z^2$, there is associated an open dual circuit that encircles the set. 
This is a standard construction that is defined rigorously on page~386 of \cite{Kesten}. In our case, the dual circuit corresponding to $K_n^{\omega}(0)$ is not entirely open; however, any edge of this dual circuit with both endpoints in  $ \mathcal{H}^+\big( - (n-1) \big) $ necessarily is open.

\begin{figure}[h]
\centering
\includegraphics[width=12cm]{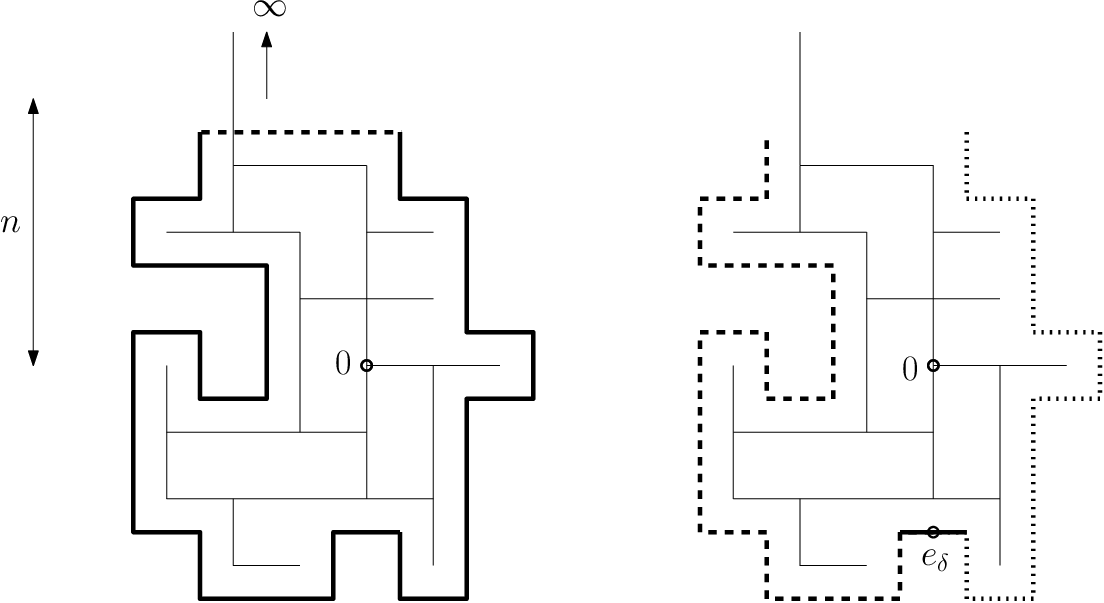}
\caption{The dual circuit in a region where $\backtrack(0)> n$.}
\label{primal}
\end{figure}

In Figure 11, the situation is illustrated for the axial choice $\vec\ell = e_1$ (with this bias direction having a downwards orientation): the set $K_n^{\omega}(0)$ consists of the primal vertices below the dotted line that are reachable in that region from $0$, while the dual circuit (which we call $\Gamma$) is comprised of the bold-face curve and this dotted line. By moving in the direction $e_1$ from $0$, the circuit $\Gamma$ must be encountered (under the assumption that $\backtrack(0) > n$ occurs) after a finite number of steps (which number we call $D$). The point of first encounter is the midpoint of a dual edge of $\Gamma$. We call that edge $e_{\delta}$.

Now removing from $\Gamma$ the edge $e_{\delta}$ and the edges with both endpoints in $\big\{ y \in \R^2: y \cdot \vec\ell \leq -(n-1) \big\}$, we are left with two disjoint open paths $\mathcal{P}_1$ and $\mathcal{P}_2$ starting at each side of $e_{\delta}$ and ending in $\big\{ y \in \R^2: y \cdot \vec\ell \leq -(n-1) \big\}$.

 Note then that 
$\big\{ \backtrack(0) \geq n \big\} \cap \big\{  D = k \big\}$ entails the existence of two disjoint open dual paths from $e_{\delta}$
each of which reaches to a distance in the direction $-\vec\ell$ of at least $n + k$. By summing over $k \in \N$, and making use of the BK-inequality (see chapter 2 in~\cite{Grimmett}) and Lemma~\ref{cor_lengthnew}, we obtain
\[
  P \big( \backtrack(0) \geq n \big) \leq
 \sum_{k=0}^\infty \exp \big\{ - 2 \big( \hat{\xi}^{(\vec{\ell})}- \epsilon  \big) (n+k) \big\}.
\]

This yields~(\ref{eqexpoub}).

\vspace{0.5cm}

{\it Step 2: Lower bound in Proposition~\ref{expo_2D}.}

\vspace{0.5cm}

We show that, for $\epsilon > 0$, and for $n \in \N$ sufficiently high,
\begin{equation}\label{eqexpolb}
{\bf P}_{p}[\mathcal{B}\mathcal{K}(0) > n] \geq  \exp\big\{ - \big( 2  \hat{\xi}^{(\vec{\ell})} + \epsilon \big) n \big\}.
\end{equation}

This is best understood by consulting Figure~12. The figure depicts the case where $\vec\ell$ is not axial. 
\begin{figure}[h] 
\centering
\epsfig{file=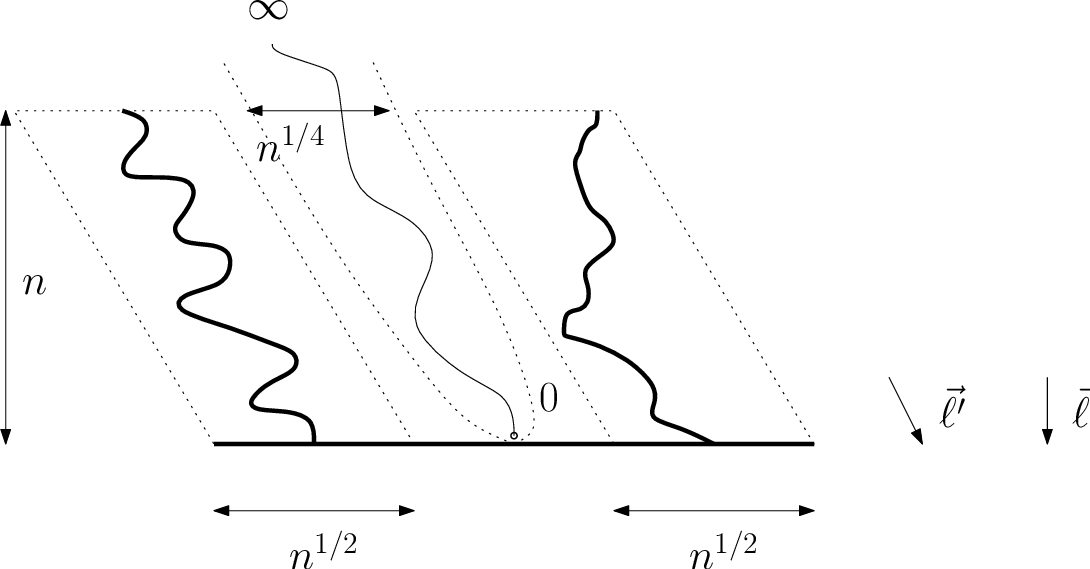, width=13cm}
\caption{A scenario such that $\mathcal{B}\mathcal{K}(0) > n$, the bold lines are dual walls and the line starting at $0$ represents an open infinite path.}\label{LB_2D1}
\end{figure}

We describe a sufficient condition for the event $\{\backtrack(0) > n\}$. Our construction of that event is depicted in Figure~\ref{LB_2D1}. We provide a precise formulation only in the case where $\vec\ell=e_1$ for the notation becomes more cumbersome in the non-axial case.  The proof for the general case is in essence the same.

Recall that an open dual path is subcritical and presents an obstacle to particle motion.
Let  $\FLOOR$ denote the set of dual edges of the form $\big( (-1/2,1/2)+i  \vec{\ell}^{\perp}, (-1/2,1/2)+  (i+1)  \vec{\ell}^{\perp} \big)$ with $\abs{i}\leq n^{1/2} + n^{1/4}\}$. We note that, under the assumption $\vec\ell=e_1$, $\vec{\ell}^{\perp}$ is simply the Euclidean unit vector $e_2$. 
We insist that each element in $\FLOOR$ be open (in the dual), an event of probability $\exp \big\{ - c n^{1/2} \big\}$. 

We also insist that an open dual path emanating from some endpoint of an element in $\FLOOR$ move upwards for a distance $n$ in the chamber indicated on the left; likewise, on the right. Each of these events has probability at least  $\exp \big\{ \big( - \hat{\xi}^{(\vec{\ell})} - \epsilon \big) n \big\}$, by Lemma~\ref{lemaxn} applied in the dual. These conditions are enough to prevent that $0 < \backtrack(0) \leq n$. Recall, however, that to exclude the possibility of $\backtrack(0)= 0$, it is necessary to establish that $0$ lies in an infinite open primal cluster. For this, we further require that the conclusion of Lemma~\ref{lemopeninh} apply for the primal percolation at $0$. By that lemma, this happens with positive probability. The four conditions being satisfied independently, all are satisfied with probability at least $c \exp \big\{ - c n^{1/2} \big\} \exp \big\{ -2 \big(  \hat{\xi}^{(\vec{\ell})}+ \epsilon \big) n \big\}$, whence~(\ref{eqexpolb}).

We have obtained Proposition~\ref{expo_2D}. Note that the statement of that proposition implies that the exponent appearing in Lemma~\ref{estim_height} is well-defined in the case $d=2$.

\vspace{0.5cm}

{\it Step 3: Proof of Lemma~\ref{estim_height}.}

\vspace{0.5cm}

We need to pursue a little further the construction used in proving (\ref{eqexpolb}): see Figure~13. We continue to assume that $\vec\ell = e_1$ since this simplifies notation.

\begin{figure}[h]\label{LB2D2}
\centering
\epsfig{file=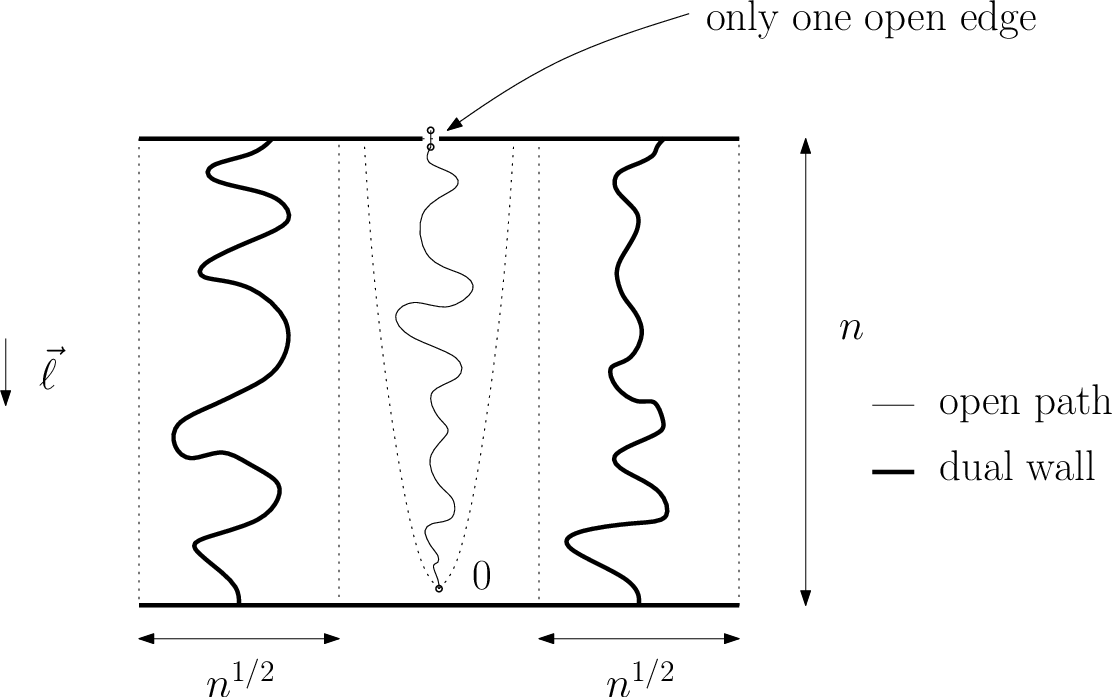, width=13cm}
\caption{The construction in Figure~12 is developed to force the presence of a one-headed trap. This picture is depicted in the axial case.}
\end{figure}

 The horizontal open dual path, the two open dual paths in the left- and right-chambers, and the infinite open primal path from the origin in $H$
force $\backtrack(0) > n$ to occur; we now introduce another requirement so that a one-headed trap of depth $n$ may be formed. 

Let $\LID$ denote the set of dual edges given by translating those in $\FLOOR$ by $-n\vec\ell$  ( which is $-n e_1$ in our case). The infinite open cluster of $0$  crosses some element in $\LID$; let $x$ be an endpoint of one of the primal edges intersecting some element of $\LID$ that is connected to $0$ by a simple open path that visits such an endpoint only at $x$. There might be several choices of $x$, in which case we choose this vertex according to some predetermined order on $\Z^2$. Suppose that, in addition to the events depicted in Figure~12, each edge in $\LID$ except that intersecting $[x,x+e_2]$ is open (in the dual), while this remaining edge is closed (in the dual).
This event has probability at least $\exp \big\{ - C n^{1/2}\big\}$.

The event is enough to ensure that the vertex $x-e_1$ is the head of a one-headed trap of depth at least $n-1$ where $x$ is one vertex among the set $\LID$ (We mention that, although the cluster of $0$ may now be finite, this does not pose a problem because of our definition for $\mathfrak{D}$). This means that
\[
\sum_{x\in \LID}P_p \big(\mathcal{T}^1(x-e_1) \neq \emptyset \text{ and } \mathfrak{D}(x)\geq n-1  \big)  \geq \exp( - C n^{1/2}) \exp(- \big( 2 \hat{\xi}^{(\vec{\ell})}+ \epsilon \big) n )
\]
there being at most $2n^{1/2} + n^{1/4}$ possible locations for $x$ to be in $\LID$, we find, by translation-invariance, that
\[
 3n^{1/2} P_p \big(\mathcal{T}^1(0) \neq \emptyset \text{ and } \mathfrak{D}(e_1)\geq n-1  \big) \geq \exp( - C n^{1/2}) \exp(- \big( 2 \hat{\xi}^{(\vec{\ell})}+ \epsilon \big) n ),
\]
where the last term, which arises from (\ref{eqexpolb}), provides a lower bound on the probability of the events depicted in Figure~13. 
This completes the proof of Lemma~\ref{estim_height} in the case $d=2$.
\qed \medskip

\subsection{Case $d\geq 3$}\label{D3}

The case of $d\geq 3$ is more complicated than the two dimensional case. Firstly, we need to introduce the concept of sausage-connection, a notion which  is intimately related to trap geometry for $d\geq 3$.

\subsubsection{Sausage-connections}

For $x,y\in\Z^d$ with $x\cdot \vec{\ell} > y\cdot \vec{\ell}$, we define $\tilde{K}^{\omega}_{x,y}(x)$ to be the open cluster of $x$ obtained using edges  both of whose endpoints lie in $ \mathcal{H}^-_x$ and at most one of whose endpoints lies in $\mathcal{H}^-_y$. 

We say that $x$ is sausage-connected to $y$ (with $x\cdot \vec{\ell} > y\cdot \vec{\ell} +2$) if
\begin{enumerate}
\item $y \in \tilde{K}^{\omega}_{x,y}(x)$,
\item $\abs{\tilde{K}^{\omega}_{x,y}(x)}< \infty $,
\item $\tilde{K}^{\omega}_{x,y}(x)\cap  \mathcal{H}^+_{x-e_1} = \{x,x-e_1\}$ and $\tilde{K}^{\omega}_{x,y}(x)\cap  \mathcal{H}^-_{y+e_1} = \{y,y+e_1\}$.
\end{enumerate}

We call $x$ the bottom of the sausage-connection and $y$ the top. See Figure~14.

\begin{figure}\label{picsausage}
\centering\epsfig{file=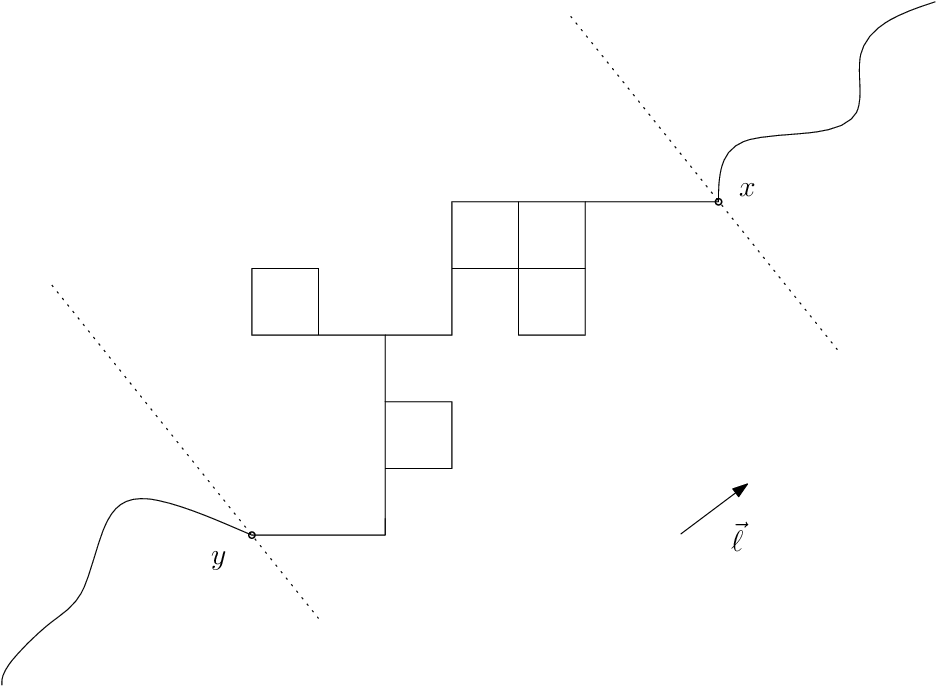, width=7cm}
\caption{An example of a sausage-connection.}
\end{figure}

This event is measurable with respect to $\omega(e)$ for edges which have both endpoints in $\mathcal{H}^-_x$ and at most one endpoint in $\mathcal{H}^-_y$.

We introduce the event that $0$ is sausage-connected to  $\mathcal{H}^-(-n)$ that there exists a unique $y\in \partial \mathcal{H}^-(-n)$ such that $0$ sausage-connected to $y$. This event is measurable with respect to $\omega(e)$ for edges which have at least one endpoint in $\mathcal{H}^-(-n)$.

 Central to our arguments is the concept of surgery, the modification of
limited regions of percolation configurations with the aim of producing a
given event.

Beyond a substantial use of surgery in Subsection \ref{secsurgery}, we will also obtain a number of other bounds by straightforward surgeries. To illustrate, we claim that
\begin{eqnarray}
 & & P_p[\text{$0$ and $ne_1$ are in the same finite cluster}] \label{eqsurgeryex} \\
 & \geq & (1-p)^{4d} P_p[\text{$0$ \, 
and $ne_1$ are sausage-connected}]. \nonumber
\end{eqnarray}

Indeed, consider a procedure which samples the law $P_p$, and then independently updates the status of the $4d$ edges which are either 
 incident to $0$ with
the other endpoint in $\mathcal{H}_0^-$, or incident to $ne_1$ with the other endpoint
not in $\mathcal{H}_{ne_1}^-$. Note that the output of this procedure is $P_p$-distributed. Were the copy of $P_p$ first realized by the procedure to satisfy the event that $0$  
and $ne_1$ are sausage-connected, and were the independent update to close each of the at most $4d$ edges, then under the resulting configuration,  
$0$ and $ne_1$ would in the same finite cluster. In this way, we see that~(\ref{eqsurgeryex}) holds. On a few occasions, we will obtain similar bounds as~(\ref{eqsurgeryex}), and will refer when we do so to finite surgery at a vertex, meaning the closure of a finite number of the vertex's incident edges.

\begin{lemma}
\label{expoa}
We have
\[
P_p[0\text{ is sausage-connected to } \mathcal{H}^-(-n)]\approx \exp(-\hatxi n),
\]
with $\hatxi\in(0,\infty)$ depends on $d$, $\vec{\ell}$ and $p$.
\end{lemma}
\noindent{\bf Proof.}
For this proof we use the notation
\[
\mu_n=P_p[0\text{ is sausage-connected to } \mathcal{H}^-(-n)].
\]

Choose $n,m \in \N$. On $\{0\text{ is sausage-connected to } \mathcal{H}^-(-n)\}$, we define a random variable $i_n$ to be the unique $y \in  \partial \mathcal{H}^-(-n)$ such that $0$ is sausage-connected to  $i_n$. In the case that there is not a unique such $y$, we set $i_n=\infty$.

The events $\{i_n=y\}$ and $\big\{ y\text{ is sausage-connected to } \mathcal{H}^-(-y\cdot \vec{\ell}- m) \big\}$ are $P_p$-independent by the definition of sausage-connection.
Hence,
\begin{equation}\label{eqdone}
 P_p[i_n=y,y\text{ is sausage-connected to } \mathcal{H}^-(-y\cdot \vec{\ell}- m)] = P_p[i_n=y] \mu_m.
\end{equation}

Now, we can see that on $\{i_n=y,y\text{ is sausage-connected to } \mathcal{H}^-(-y\cdot \vec{\ell}- m)\}$, we have that $0$ is sausage-connected to
 $\mathcal{H}^-(-y\cdot \vec{\ell}- m)$. By a finite surgery at the top of a sausage-connection from $0$ to $z\in  \mathcal{L}(-y\cdot \vec{\ell}- m)$, it is always possible to construct a sausage-connection from $0$ to some $z'\in 
\mathcal{L}(-n-m-2)$ (where $z'$ is $z$, $z+e_1$ or $z+2e_1$). Hence, for any $z\in  \mathcal{L}(-y\cdot \vec{\ell}- m)$ we have
\begin{eqnarray}
 & & P_p[i_n=y,~0\text{ is sausage-connected to } z] \nonumber \\
 & \leq & CP_p[i_n=y,~0\text{ is sausage-connected to } z'], \nonumber
\end{eqnarray}
where the constant can be chosen independently of $z$ and $z'$, which implies that
\begin{align}\label{eqdtwo}
   & \, P_p[i_n=y,y\text{ is sausage-connected to } \mathcal{H}^-(-y\cdot \vec{\ell}- m)]  \\ \nonumber
   = & \, \sum_{z\in  \partial \mathcal{H}^-(-y\cdot \vec{\ell}- m)} P_p[i_n=y,~0\text{ is sausage-connected to } z]  \\ \nonumber 
  \leq & \, 3C \sum_{z'\in  \partial \mathcal{H}^-(-n- m-2)} P_p[i_n=y,~0\text{ is sausage-connected to } z']\\ \nonumber 
  \leq & \, C P_p[i_n=y,~0\text{ is sausage-connected to } \mathcal{H}^-(-n- m-2)], 
 \end{align}
 where $C > 0$ may be chosen independently of $y$, $n$ and $m$.
 
 Using (\ref{eqdone}) and (\ref{eqdtwo}), and the fact that  $\{i_n=y\}_{z\in \Z^d}$ partitions the event $\{0\text{ is sausage-connected to } \mathcal{H}^-(-n)\}$, we obtain
\[
\mu_{n+m+2}\geq C^{-1} \mu_n \mu_m,
\]
which, by a super-additivity theorem (see (II.10) p.400 in~\cite{Grimmett}), implies that
\[
P_p[0\text{ is sausage-connected to } \mathcal{H}^-(-n)]\approx \exp(-\hatxi n),
\]
with $\hatxi\in [0,\infty)$.

It is simple to exclude the case $\hatxi=0$: indeed, by a finite surgery at both ends similar to that which provides~(\ref{eqsurgeryex}), a sausage-connection may be turned into a finite cluster of size at least $n$, so that 
\begin{eqnarray*}
& & P_p[0\text{ is sausage-connected to } \mathcal{H}^-(-n)]  \\
 & = & \sum_{y\in \partial \mathcal{H}^-(-n)} P_p[i_n=y] \\
 & \leq & C\sum_{y\in \partial \mathcal{H}^-(-n)} P_p[i_n=y, n\leq \abs{K^{\omega}(0)}<\infty] \\
 & \leq & C P_p[n\leq \abs{K^{\omega}(0)}<\infty],
\end{eqnarray*}
where we used that $\{i_n=y\}$ for $y\in\mathcal{H}^-(-n)$ form a partition of the event $\{0\text{ is sausage-connected to } \mathcal{H}^-(-n)\}$.

The probability that the open cluster of a given point in supercritical percolation has a finite size exceeding $n$ has an exponential decay in $n$   (see Theorem~8.21 of~\cite{Grimmett}). Hence, we are done. 
\qed \medskip

\subsubsection{Heuristics identifying the backtrack exponent for $d\geq 3$}

We will show the following result, which is Proposition~\ref{backtrack_exponent} for $d\geq 3$. 
\begin{proposition}
\label{expo_3D}
In $\Z^d$ with $d\geq 3$, we have
\[
{\bf P}_p[\mathcal{B}\mathcal{K}(0) > n] \approx \exp(-\hatxi n),
\]
where $\hatxi$ is defined in Lemma~\ref{expoa}.
\end{proposition}

We explain the main ideas of the proof before continuing. We wish to show that the probability of $\backtrack(0) > n$ and the event that $0$ is sausage-connected to $\mathcal{H}^-(-n)$ have the same rate of exponential decay.

 Of the two required bounds, one is easy; indeed, given a sausage-connection from $0$ 
to  $\mathcal{H}^-(-n)$, a finite surgery near $0$ will force $\backtrack(0) \geq n$.

  The harder bound is the other one, where we aim to turn a deep trap (which realizes the event $\{ \backtrack(0)\geq n \}$) into a sausage-connection at low probabilistic cost. To prove it, we notice that a typical region satisfying $\backtrack(0) >n$ contains a trap surface, which blocks particle motion, of length at least $n$ in  the direction $\vec\ell$. To create a long sausage-connection, we want to take that trap surface
  \begin{enumerate}
  \item make an entry at the bottom of the trap, and
  \item seal off all but one point of the top of the trap.
  \end{enumerate}
  
 The real difficulty for us is sealing off the top of the trap. Indeed, we have no {\it a priori} bound on the size of the top of the trap surface, and if it were too large then the cost of the sealing surgery might be too high.  However, if this surface were typically narrow near the top, we could indeed perform the surgery to close down this bottleneck in the trap at a low probabilistic cost. The result would be a finite open cluster that is long in the $\vec\ell$-direction. This is in essence the same as finding a long sausage-connection, as we may see by using finite surgery arguments.

We now explain why we expect the trap surface to be narrow in many places, under the law 
$P_p \big(\cdot \big\vert \backtrack(0) > n \big)$.  At the trap surface, we may find pairs of points, one on the outside and one on the inside of the trap. These points have the property that they are close in the $L^1$-distance; however, unless they happen to be close to the top of the trap, they are far in the chemical distance. Such pairs of points are known to be very unlikely, by the regularity of the supercritical phase of percolation.

Recall that two vertices $u,v \in \Z^d$ are said to be $*$-neighbors if $\max_{1 \leq i \leq d}\vert u_i - v_i \vert \leq 1$.
Let $F \subseteq \Z^d$, and let $A \subseteq F$. We write
\begin{enumerate}
\item $\extg^F(A)$ for the set of vertices in $F$ that $*$-neighbor $A$ and from which there exists an infinite simple  path in $F$ that does not meet $A$, and
\item $\text{INS}^F(A)$, the {\em inside} of $A$ in $F$, for the subset of vertices in $F \setminus A$ from which every infinite simple  path in $F$ meets $A$.  
\end{enumerate}

Any set of the form $\{z\in \Z^d,~(z-x)\cdot \vec u \geq 0\}$, with $x\in \R^d$ and $\vec u \in S^{d - 1}$, will be called a half-space.

We firstly need a section on geometrical results in $\Z^d$ with $d\geq 3$; the reader may skip the proofs on a first reading.

\subsubsection{Some results about connectedness in half-spaces}


\begin{lemma}\label{lemstarconn}
 Let $H$ be a half-space, and let $C \subseteq H$. Then 
 $\partial_{\rm ext}^H(C)$ is nearest-neighbor connected. 
 \end{lemma}
\noindent{\bf Remark.}
In the case that $H = \Z^d$, H. Kesten provided a careful argument for the statement that made use of algebraic topology (see Lemma 2.23 in~\cite{Kesten1}). Recently, A. Tim\'ar \cite{timar} has provided an elementary proof of this result that may be adapted more readily to other contexts.

\noindent{\bf Proof.}
For this proof, we slightly abuse notation, denoting by $H$ the graph obtained from the half-space $H$ with the nearest-neighbor notion of
adjacency and by $H^*$ this graph with the $*$-connected notion of adjacency.

The result follows from Theorem $3$ of \cite{timar}, whose hypotheses
we must verify. In the notation of that theorem, we take $G = H$ and $G^+= H^*$.

We have two conditions to verify. For the first one, we need to show that any cycle $C$ in
$H$ may be generated by basic $4$-cycles in $H$, (where a basic $4$- cycle means a non-trivial cycle of length $4$). As noted in \cite{timar}, this property is verified in the case of $H=\Z^d$ with the nearest neighbor topology.  We will reduce to this case. To do so,  note that, in $\Z^d$, each element of the sequence of basic $4$-cycles generating $C$ is included in $C+B_C$, where $B_C$ is a bounded region. We may choose $B_{C'}=B_C$ for any translate $C'$ of $C$.

Now given any cycle $C$ in $H$, we consider the translate $C_M$ of $C$ by $M$ units in the $e$-direction, where $e$ is a standard unit vector selected so that $e$ is directed into $H$ (that is $d_{\Z^d}(ne,H^c)\to \infty$ as $n \to \infty$). We know that $C_M$ can be generated by basic $4$-cycles included in $C_M+B_C$; by choosing $M$ large enough, $C_M+B_C$ is included in $H$. Now, using basic $4$-cycles in $H$, we can transform $C_M$ into $C$. Hence we can generate $C$ with basic $4$-cycles in $H$.

For the second condition of Theorem $3$, we must show that,
for any $*$-adjacent points $x,y$ in $H$, we can find a nearest-neighbor path in $H$ from $x$ to $y$ such that
any two vertices in this path are $*$-adjacent. 
By relabeling the unit vectors, we may assume that the coordinates of $y - x$ are each either zero or one.
Let us denote $1\leq w_1<w_2 \leq \ldots \leq w_k \leq d$ the coordinates of $y - x$ which are equal to one. Let $l(z):\Z^d
\to \R$ denote the distance to the half-space  $\{x\in \R^d: x\cdot \vec{\ell}<0\}$ in $\R^d$.
Then $\min \{ l(x),l(y) \} > 0$, since $x,y\in H$. A nearest-neighbor path from $x$ to
$y$ may be formed by making the successive displacements $e_{w_1},\ldots,e_{w_k}$; and others may be formed by choosing a different
order on the displacements. We choose the permutation so that all moves $e_{w_i}$ for which
$l(e_{w_i})>0$ are taken first. In this way, $l(v) > 0$ holds for each $v$ in the
resulting path. This implies that the resulting path, which goes from $x$ to $y$, lies in $H$.
\qed \medskip

\begin{lemma}
\label{lemhtrap}
Let $d \geq 2$, and let $H$ denote $\Z^d$, or a half-space in $\Z^d$. 
Let $J \subseteq H$ be finite.
Let $C_1,\ldots,C_r$ denote the $*$-connected components in $H$ of $J \cap H$. 
Then
$$
  \ins^H(J) = \cup_{i=1}^r \ins^H(C_i).
$$
\end{lemma}

\noindent{\bf Proof.}
We give the proof for the case $H = \Z^d$, since the proof is the same in the other case. Clearly $\ins(C_i) \subseteq \ins(J)$ for each
$i \in \{ 1,\ldots, r \}$. To show the opposite inclusion, consider 
\begin{equation}\label{xassump}
x \not\in  \cup_{i=1}^r \ins(C_i). 
\end{equation}
Let $(\mathcal{P}(i))_{i\geq 0}$ be an arbitrary simple infinite path starting at $x$. If $\mathcal{P}$ does not  meet $J$, then $\mathcal{P}$ demonstrates that $x \not\in \ins(J)$, and we are done.

 In the other case, let $j \in \N$ be minimal such that $\mathcal{P}(j) \in J$ and $k \in \{1,\ldots ,r\}$ be such that $\mathcal{P}(j) \in C_k$. We claim that $\mathcal{P}(j-1) \in \partial_{\rm ext}(C_k)$.
 
 Indeed,~(\ref{xassump}) implies that there exists an infinite simple path emanating from $x$ that does not meet $C_k$. Concatenating this path to the reversed path $[\mathcal{P}(j-1),\ldots,\mathcal{P}(0)]$ yields an infinite path from $\mathcal{P}(j-1)$ that has at most finitely many self-intersections and never meets $C_k$. The vertex $\mathcal{P}(j-1)$ being a nearest neighbor of $C_k$, this path demonstrates that $\mathcal{P}(j-1) \in \partial_{\rm ext}(C_k)$.

Now take $j' \in \N$ to be the index of the last visit $\mathcal{P}(j')$ by $\mathcal{P}$ to $C_k$. The subpath of $\mathcal{P}$ starting at $\mathcal{P}(j' + 1)$ demonstrates that $\mathcal{P}(j' + 1) \in  \partial_{\rm ext}(C_k)$.

The set $\partial_{\rm ext}(C_k)$ being nearest-neighbor connected by Lemma \ref{lemstarconn}(i), there exists a  path $[q_0,\ldots,q_m]$ from $\mathcal{P}(j-1)$ to $\mathcal{P}(j' + 1)$ in this set. Since $C_k$ is one of the $*$-connected components of $J$, we see that $\partial_{\rm ext}(C_k)$ is disjoint from~$J$. 

We consider the path $\mathcal{P}_1$ which is modified from the original path $\mathcal{P}$ by substituting for $[\mathcal{P}(j-1),\ldots,\mathcal{P}(j' + 1)]$ the new path  $[\mathcal{P}(j-1) = q_0,\ldots,q_m = p_{j' + 1}]$. Note that $P_1$ restricted to the index set $\big\{ 0,\ldots, p_{j' + 1} \big\}$ is a simple path that is disjoint  from $C_k$ and also from each of the other $C_i$. 

If $\mathcal{P}_1$ does not meet $J$, then $\mathcal{P}_1$ demonstrates that $x \not\in \ins^H(J)$. Otherwise we simply repeat the above procedure, choosing $\mathcal{P}_1$ in place of $\mathcal{P}$ in the above, and so defining a new path $\mathcal{P}_2$. Iterating this at most $r$ times we obtain a path $\mathcal{P}'$ that is disjoint from $J$ which demonstrates that $x \not\in J$. This completes the proof. \qed \medskip

\subsubsection{The coarse-grained trap surface is a one-dimensional object}

For this section we use $\dtrap$ as a shorter notation for $K_n^{\omega}(0)$ which recall is the open cluster of $x$ in $\mathcal{H}^{+}(x\cdot \vec{\ell}-n)$. 
In particular $\mathcal{B}\mathcal{K}(0) > n$ implies that
the set  $\dtrap$ is non-empty and finite. We will refer to $\dtrap$ from time to time as a trap.

Let $K \in \N$. A $K$-box will refer to any $B_K(z)=K z + [-K,K]^d$ for some $z \in \Z^d$.

Let $D \subseteq \Z^d$. Then $\overline{D}$ is defined to be the set of $K$-boxes intersecting $D$.

Using the one-to-one correspondence between $K$-boxes and $\Z^d$, we may extend previous notions to the case of $K$-boxes. We recall that $\text{INS}$ and $\extg$ were defined after Proposition~\ref{expo_3D}. The notation $\extg$, $\text{INS}$,  $\extg^H(A)$ and $\text{INS}^H(A)$ extends to $\extk$, $\text{INS}_K$, $\extk^{\overline{H}}$ and $\text{INS}^{\overline{H}}_K(A)$ where $\overline{H}$ is a half-space in $K$-boxes. In fact, we will henceforth write $\overline{H}$ for the half-space in $K$-boxes defined by $\{B_z(K): z\in \mathcal{H}^+(-\lfloor n/K \rfloor + 10)\}$. The presence of the $+ 10$ term in this definition ensures that any point in the trap $\dtrap$ lying in some $K$-box in $\overline{H}$ is not close to the top of the trap; any such point will thus be at least of order $10K$ in chemical distance from any nearby point not lying in $\dtrap$, because an open path connecting such a pair of points must leave the trap.  

We will now prove the existence of a set $\overline{\Fsurf}$ which should be considered to be a coarse-grained trap surface. 
\begin{lemma}
\label{lemetrap}
On $\{\mathcal{B}\mathcal{K}(0) > n\}$, there exists a $*$-connected component $\overline{\Fsurf}$
of $\partial_{{\rm ext},K} \overline{\dtrap} \cap \overline{H}$ with the property that 
\[
B_K(0) \in \text{INS}_K^{\overline{H}} \big( \overline{\Fsurf} \big).
\]
\end{lemma}

\noindent{\bf Proof.}
On $\{\mathcal{B}\mathcal{K}(0) >n\}$, we claim that 
\begin{equation}
\label{bkcl}
B_K(0) \in \text{INS}_K^{\overline{H}} \big(\partial_{{\rm ext},K} \overline{\dtrap} \cap \overline{H} \big).
\end{equation}
To derive (\ref{bkcl}), note firstly that $B_K(0)\in \overline{\dtrap}$ implies that $B_K(0) \notin \partial_{{\rm ext},K} \overline{\dtrap}$. Hence, it is enough to show that  a simple infinite  path in $K$-boxes in $\overline{H}$ that begins at $B_K(0)$ necessarily meets $\partial_{{\rm ext},K} \overline{\dtrap}$. On $\{\mathcal{B}\mathcal{K}(0) >n\}$, the set $\dtrap$ is finite, so any simple infinite  path will eventually leave $\overline{\dtrap}$, thus intersecting $\partial_{{\rm ext},K} \overline{\dtrap}$.

The existence of $\overline{\Fsurf}$ now follows from Lemma~\ref{lemhtrap}.
\qed \medskip

It is $\overline{\Fsurf}$ which we will argue is small by means of the renormalization argument invoking regularity of supercritical percolation to which we have already alluded and which is now further illustrated by Figure 15.

\begin{figure}
\centering\epsfig{file=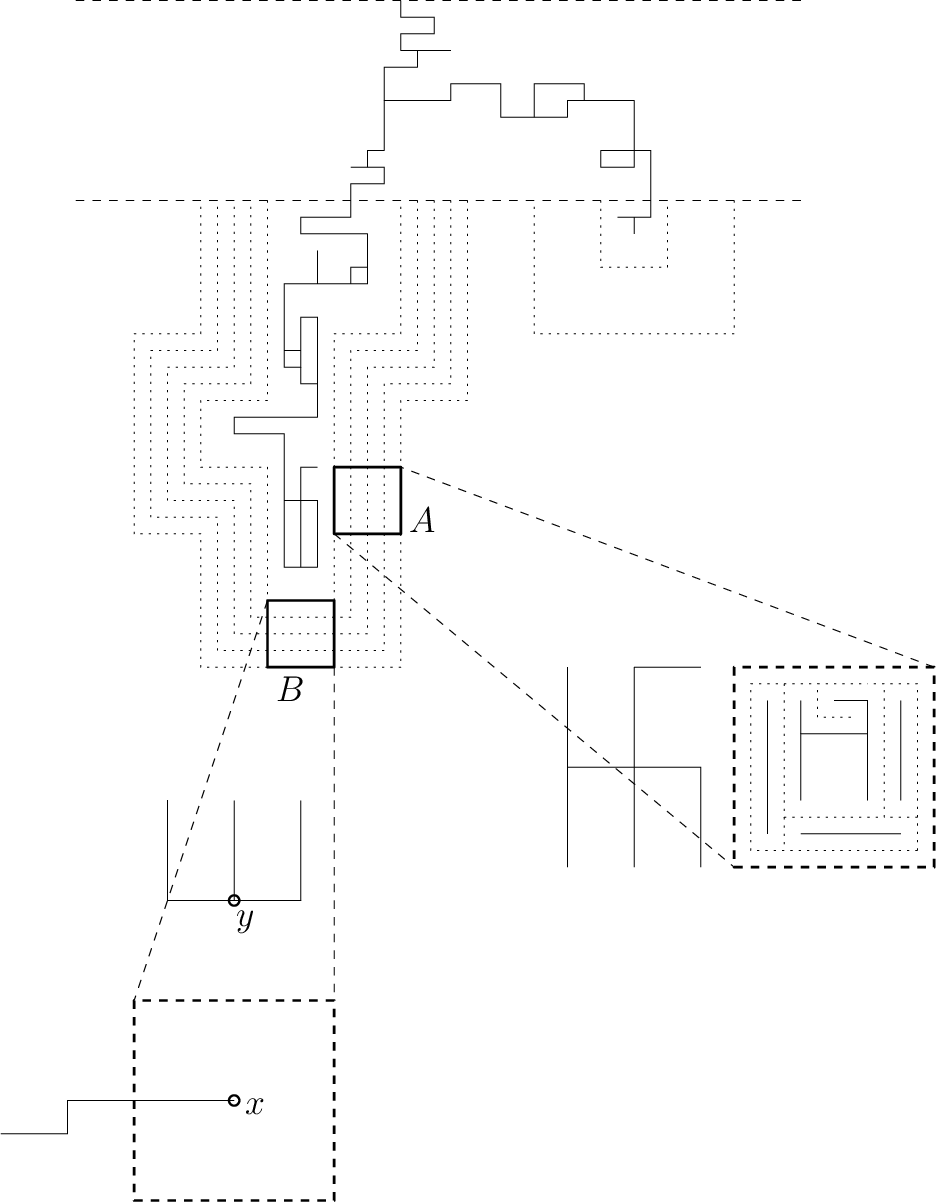, width=8cm}
\caption{The coarse-grained trap surface $\overline{\Fsurf}$ and how its $K$-box elements are bad. The dimension is chosen to be two because this case is easier to depict. In the main sketch, the trap $\dtrap$ is shown with solid lines, $A$ and $B$ are examples of $K$-boxes, and the region below the lower of the two long horizontal dashed lines is the union of the $K$-boxes that comprise $\overline{H}$. In fact, the slice between these two long horizontal lines is of width ten in $K$-boxes (and not three as depicted). The region indicated by parallel dotted lines is $\overline{\Fsurf}$. Each of the $K$-boxes therein is bad, and the enlarged details illustrate the two reasons that this may be so for a given such $K$-box. For the enlargement of the box $A$, note that the dotted lines indicate the open dual edges. In the region of $A$, the trap surface is comprised of many small bubbles, one on top of another, so that no long open primal path emanates from inside $A$. In the case of $B$, the neighboring trap surface is smoother, and a long open path emanates from $x \in B$; it must travel far to reach the nearby $y$ inside the trap, making $B$ a bad box as well.}
\end{figure}

A $K$-box $B_K(x)$ is said to be bad if either of the following conditions is satisfied:
\begin{enumerate}
\item there does not exist an open path from $Kx+[-K,K]^d$ to $Kx+\partial [-3K,3K]^d$,
\item there is some $y$ with $\abs{\abs{x-y}}_{\infty}\leq1$ such that there exist open paths, one from $Kx+ [-K,K]^d$  to $Kx+\partial [-3K,3K]^d$, and the other from $Ky+\partial [-K,K]^d$ to $Kx+\partial [-3K,3K]^d$, which lie in disjoint open clusters of $Kx + [-3K,3K]^d$.
\end{enumerate}
A $K$-box that is not bad will be called good.
 \begin{lemma}
 \label{lembadbox}
 Let $\epsilon > 0$. For $K \in \N$ large enough, the set of bad $K$-boxes is stochastically dominated by an independent percolation on $K$-boxes of parameter $\epsilon > 0$.
 \end{lemma}
\noindent{\bf Proof.}
 The set of bad $K$-boxes forms a $6$-dependent percolation. The probability that a given $K$-box is bad tends to zero with $K \to \infty$ by Lemma~7.89 of~\cite{Grimmett}. Hence, the statement follows from Theorem 0.0 (i) in~\cite{LSS}.
 \qed \medskip
 
We also require the following result.

\begin{lemma}
On $\{\mathcal{B}\mathcal{K}(0)>n\}$, for $n$ large enough, every element of $\partial_{{\rm ext},K} \overline{\dtrap} \cap \overline{H}$ is a bad $K$-box.
\end{lemma}
\noindent{\bf Proof.}
Let $B \in \partial_{{\rm ext},K} \overline{\dtrap} \cap \overline{H}$ and choose $z\in\Z^d$ to be such that $B=B_z(K)$. Aiming to prove the lemma by contradiction, we assume that $B$ is a good $K$-box. This assumption implies that
\[
\text{there exists an open path $\mathcal{P}$ from $Kz+[-K,K]^d$ to $Kz+\partial [-3K,3K]^d$}.
\]

By $B\in \partial_{{\rm ext},K} \overline{\dtrap}\cap \overline{H} $, we see that there exists $z_1\in \Z^d$ for which $\abs{\abs{z-z_1}}_{\infty}\leq 1$ with $B_K(z_1) \in \overline{\dtrap}$. Hence, we may take $y\in B_K(z_1)\cap \dtrap$. By the definition of $\dtrap$, $y$ is connected to $0$ by an open path which in turn is thus connected to $\mathcal{H}^-(-n)$ since $\mathcal{B}\mathcal{K}(0) >n$ occurs. This implies that
\[
y\in Kz_1+ [-K,K]^d \text{ is connected by an open path $\mathcal{P}'$ to } Kz+\partial [-3K,3K]^d.
\]

Now, the box $B$ being good, we find that $\mathcal{P}$ and $\mathcal{P}'$ intersect, so that $y \in  \dtrap$ is connected to $B$ inside $Kz+[-3K,3K]^d$. Moreover, the box $B=B_z(K)$ belonging to $\overline{H}$, we have that $Kz+[-3K,3K]^d\subset \mathcal{H}^+(-n)$; thus, an open path in $\mathcal{H}^+(-n)$ connects some element of $B$  to $y$. 

Now recalling that $y\in \mathcal{T}$, this means that $B\in \overline{\dtrap}$, which contradicts $B \in  \partial_{{\rm ext},K} \overline{\dtrap}$.

This establishes that $B$ must be bad.
 \qed \medskip

We now reach the principal statement of this subsection which asserts that the $K$-box  trap surface has a size growing linearly in the trap height.
Recall that $\overline{\Fsurf}$ is a random set of $K$-boxes whose distribution depends on $K$.
\begin{lemma}
\label{lemfcn}
There exists $K$ large enough such that the following property holds: for any $C_2<\infty$, there exists $C_1<\infty$ such that
\[
P_p \Big( 
 \big\vert \overline{\Fsurf} \big\vert \geq  C_1n \Big\vert \mathcal{B}\mathcal{K}(0) > n \Big)
 \leq \exp \bigl( -C_2 n  \bigr).
\]
\end{lemma}   

\begin{remark}\label{decompo_regen}
Invoking Lemma~\ref{lemfcn}, it may then be argued that a uniformly positive fraction of cross-sections have an intersection with the trap surface whose size is bounded above uniformly. Surgeries in uniformly bounded regions then  imply the existence of a positive-density sequence of cutpoints whose consecutive elements are sausage-connected.
The trap thus has a renewal structure under which the irreducible component -- the sausage-connection -- has size of bounded mean. This inference is the key input needed to analyse trap geometry and probability by means of Ornstein-Zernike theory.
\end{remark}

\noindent{\bf Proof.}
We recall that $e_1$ is one of the Euclidean unit vectors with maximal scalar product with $\vec{\ell}$.  The path $(ie_1)_{i\geq 0}$ being nearest-neighbor, we know by Lemma~\ref{lemetrap}
that $\overline{\Fsurf}$ necessarily contains a $K$-box of the form $B_K(i_0 e_1)$ for some $i_0 \in \N$. Hence, $\abs{\overline{\Fsurf}} \geq C_1 n$ entails the existence of a $*$-connected component of bad boxes starting at $B_K(i_0 e_1)$ that has size at least $\max(C_1n, i_0-10)$. 
By Lemma 5.1 in~\cite{Kesten}, there exists $C_3 > 1$ such that 
 the number of $*$-connected components of size $n$ is less than $C_3^n$. Taking $K$ large enough (for the choice of 
$\epsilon= (3C_3)^{-1}$ in Lemma~\ref{lembadbox}), we see that 
\[
P_p \Big( \abs{\overline{\Fsurf}} \geq C_1n \Big) 
  \leq  C_1n 3^{-C_1 n}+\sum_{i_0\geq C_1n-10} 3^{-i_0} \leq C 2^{-C_1n}.
\]

 By considering an explicit trap, such as a succession of open vertical edges surrounded by closed edges, we have the bound 
\[
P_p \Big( \backtrack(0) >n \Big) \geq \exp \big\{ -Cn \big\},
\]
for some $C <\infty $. Then using the two preceding equations, we have that
\[
P_p \Big(  \big\vert \overline{\Fsurf} \big\vert \geq  C_1n \, \Big\vert \, \mathcal{B}\mathcal{K}(0) > n \Big)
 \leq  C2^{-C_1n}\exp \big\{ Cn \big\},
\]
so that choosing $C_1$ large enough we are done.
 \qed \medskip

\subsubsection{Surgery closes the narrow top and produces a sausage-connection}\label{secsurgery}
We now give definitions and lemmas needed to carry out surgery to close down the edges in a narrow bottleneck at small macroscopic distance from the top of the trap. The output of surgery is in essence a sausage-connection travelling over a distance comparable to the trap height.

From now on $K$ and $C_1$ will be fixed so that Lemma~\ref{lemfcn} is verified with $C_2=1$. 
\begin{definition}
We partition the $K$-box half-space $\overline{H} = \big\{ \overline{H}_1,\overline{H}_2,\cdots \big\}$,
where $\overline{H}_i = \big\{ B \in \overline{H}: B = B_y(K), \,  -n+10+3i \leq y \cdot \vec{\ell}  < -n+10+3(i+1) \big\}$. 
\end{definition}

For a set $\overline{A}$ of $K$-boxes, we introduce the notation $[\overline{A}] =\{y\in B \text{ with } B\in \overline{A}\}$ for the points contained in the $K$-boxes in $\overline{A}$.

\begin{lemma}
\label{blaaa}
Let $\epsilon \in (0,1/6)$. There exists a non-random $M_n \in \big\{ 1,\ldots, \lfloor \epsilon n \rfloor \big\}$ such that
\[
P_p \Big(  \big\vert \overline{\Fsurf} \cap \overline{H}_{M_n}  \big\vert \leq 2 C_1  \epsilon^{-1} , [\overline{\Fsurf}]\subseteq B_{\Z^d}(0,Cn) \, \Big\vert \, \backtrack(0) > n \Big) \geq \frac 1{3\epsilon n}.
\]
\end{lemma}
\noindent{\bf Proof.}
By Lemma~\ref{lemfcn}, it holds that
\begin{equation}\label{eqfnotbig}
P_p \Big( [\overline{\Fsurf}] \cap B_{\Z^d}(0,C_1n)^c \not = \emptyset \, \Big\vert \, \backtrack(0)> n \Big) \leq \exp (- n ),
\end{equation}
and
\[
P \Big(  \big\vert \overline{\Fsurf} \big\vert \leq C_1 n \, \Big\vert \, \backtrack(0) >n \Big) \geq 1/2.
\]

Let $m \in \N\cup \{\infty\}$ be minimal such that
$\big\vert \overline{\Fsurf} \cap \overline{H}_m \big\vert \leq 2 C_1  \epsilon^{-1}$.
Note that $\vert \overline{\Fsurf} \vert \leq C_1 n$
implies that $m \leq \epsilon n$, because $\big\{ \overline{\Fsurf}\cap\overline{H}_i: i\in \N \big\}$ is a partition of $\overline{\Fsurf}$.
Let $M_n \in \big\{ 1,\ldots, \lfloor \epsilon n \rfloor \big\}$
denote the value such that $\{m = M_n\}$ is the most probable under 
$P \big( \cdot \big\vert \backtrack(0) > n \big)$ (taking say the smallest among such values if there is a choice to be made). 
We find then that
\[
P \Big(  \big\vert \overline{\Fsurf} \cap \overline{H}_{M_n} \big\vert \leq 2 C_1  \epsilon^{-1} \, \Big\vert \, \backtrack(0) >n  \Big) \geq \frac 1{2\epsilon n}.
\]

The statement of the lemma follows from (\ref{eqfnotbig}).
\qed \medskip

For $x\in \Z^d$, we define the slab
\[
\mathcal{SL}_x = \Big\{ z \in \Z^d:   -3< (z-x) \cdot \vec{\ell} \leq 0  \Big\},
  \]
 Let $u \in \R^d$, $u \not= 0$, satisfy $u \cdot \vec\ell = 0$. 
Write $R_u(x) = x + \big\{ t u: t \geq  0 \big\}$.
    
 The proofs of the next two lemmas are simple and we omit them. They essentially state that we can find infinite paths in tubes of radius larger than three.
  
  \begin{lemma}\label{lemray}
 For $x\in \Z^d$, let $y \in \mathcal{SL}_x$.  There exists an infinite simple  path $\mathcal{P}_{y,u}$ in $\mathcal{SL}_x$ that starts at $y$ and is contained in the $3$-ball, in the Euclidean metric, about $R_u(x)$.
\end{lemma}

A simple consequence of this is the following.
  \begin{lemma}
  \label{connectivity}
For any point $x \in\Z^d$ and any $y\in\mathcal{SL}_x $,   there is an infinite simple  path in $\mathcal{SL}_x $ starting  from $y$.
  \end{lemma}

Let us prove the following.

\begin{lemma}
\label{isop_ineg}
There exists a universal constant $C$, such that, for any $x\in\Z^d$, $F \subseteq  \mathcal{SL}_x$ implies that  $\vert \text{INS}^{\mathcal{SL}_x }(F) \vert \leq C \vert F \vert^{d+1}$. 
\end{lemma}
\noindent{\bf Proof.}
First, let us show that there exists a constant $C$ such that for any~$n$ and for any $y\in \mathcal{SL}_x$ there exist $n$  paths in $\mathcal{SL}_x$ starting from $y$ which are disjoint in $\mathcal{SL}_x\setminus B_{\Z^d}(y,Cn)$.

 We can choose $\vec u_1,\ldots, \vec u_n \in S^{d-1}$, each orthogonal to $\vec{\ell}$, where the angle between any two vectors is at least $2\pi/n$. This may be seen because there is a copy of $S^1$ embedded in any plane orthogonal to $\vec{\ell}$. Consider the paths $\mathcal{P}_{y,\vec u_i}$ provided by Lemma~\ref{lemray}. For $C$ large enough, outside of the ball $B_{\Z^d}(y,Cn)$, each pair of rays $R_{\vec u_i}(y)$ are at distance at least six. By the containment stated in Lemma~\ref{lemray}, the paths $\mathcal{P}_{y,\vec u_i}$ are mutually disjoint outside of $B_{\Z^d}(y,Cn)$.   

To prove the lemma, take $y\in \text{INS}^{\mathcal{SL}_x }(F)$. 
We apply the result that we just proved with $n= \abs{F}+1$; by the pigeon hole principle  one of these paths (say $\mathcal{P}$) starting at $y$ cannot contain a point of $F$ in its intersection with $\mathcal{SL}_x\setminus B_{\Z^d}(y,C(\abs{F}+1))$ (where all the $\abs{F}+1$ paths are disjoint). However, the path $\mathcal{P}$ must intersect $F$, because $y\in \text{INS}^{\mathcal{SL}_x }(F)$; we see that any point of intersection must lie in $B_{\Z^d}(y,C(\abs{F}+1))$. Hence, we know that any $y\in\text{INS}^{\mathcal{SL}_x }(F) $ is distance less than $C(\abs{F}+1)$ of a point in $F$; by a relabelling of $C$, we conclude that $\abs{\text{INS}^{\mathcal{SL}_x }(F) }<C \abs{F}^{d+1}$.
\qed \medskip

We now find a width-three $K$-box slab near the trap top whose intersection with the coarse-grained trap interior has uniformly bounded size:
\begin{lemma}
\label{we_have_M}
Let $\epsilon \in (0,1/6)$. Set $\const = C \big( 2 C_1 \epsilon^{-1} \big)^{d+1}$, where the constants $C_1$ and $C$ appear in Lemmas~\ref{blaaa} and~\ref{isop_ineg}. Then
\[
P \Big(   \big\vert \text{INS}^{\overline{H}_{M_n}}\big( \overline{\Fsurf} \big)  \big\vert \leq \const, [\overline{\Fsurf}]\subseteq B_{\Z^d}(0,C_1n) \, \Big\vert  \, \backtrack(0) > n \Big) \geq \frac 1{3\epsilon n},
\]
where $M_n$ is defined in Lemma~\ref{blaaa}.
\end{lemma}
\noindent{\bf Proof.} A consequence of Lemmas~\ref{blaaa} and~\ref{isop_ineg}. \qed

Under $P \big( \cdot \big\vert \backtrack(0) > n \big)$, we have identified a $K$-box slab $\overline{H}_{M_n}$ through which  the inside of the trap passes in a very narrow region with probability of order $1/n$. This narrow region is about to be closed down by surgery. 

 We introduce the event of the trap having a narrow top, setting
\[
NT(n)=\Big\{ \big\vert  \text{INS}_K^{\overline{H}_{M_n}}(\overline\Fsurf) \big\vert 
   \leq \const \Big\}
\cap \Big\{ [\overline{\Fsurf}] \subseteq B_{\Z^d}(0,C_1n) \Big\} 
\cap \Big\{ \backtrack(0) >n \Big\}.
\]
The plan is to close all of the edges lying in the $K$ boxes comprising   $\text{INS}_K^{\overline{H}_{M_n}}(\overline\Fsurf)$, with the cost of doing so being reasonable when $NT(n)$. The next lemma shows that in the outcome of such surgery, the origin is in a finite cluster.
\begin{lemma}
\label{mnbvc}
Under the event $\big\{ \backtrack(0) >n \big\}$, any infinite simple open  path in $\Z^d$ starting at $0$ passes through $\big[ \text{INS}_K^{\overline{H}_{M_n}}(\overline\Fsurf) \big]$.
\end{lemma}
\noindent{\bf Proof.}   On $\big\{ \backtrack(0) >n \big\}$, any infinite simple open  path in $\Z^d$ starting at $0$ must pass through one of $K$-boxes in $\overline{H}_{M_n}$. Fixing such a path $\mathcal{P}$, let $B$ be the such first box so encountered. By $\Fsurf \subseteq \extk \overline{\dtrap}$ and $B\in \overline{\mathcal{T}}$ we know that $B \not\in \Fsurf$. Moreover, on the event $\big\{ \backtrack(0) > n \big\}$, every infinite simple $K$-box path from $B$ in $\overline{H}_{M_n}$ must meet $\overline\Fsurf$. This demonstrates that $B\in  \text{INS}_K^{\overline{H}_{M_n}}(\overline\Fsurf)$ and establishes the lemma. \qed \medskip

Lemma \ref{mnbvc} permits the following definition.
\begin{definition}
On the event $\big\{ \backtrack(0) > n \big\}$, choose a deterministic procedure for selecting an  open path $\mathcal{O}$ leaving $0$ that ends on its first visit to  $\big[ \text{INS}_K^{\overline{H}_{M_n}}(\overline\Fsurf) \big]$. Let $\bigxlower \in \big[ \text{INS}_K^{\overline{H}_{M_n}}(\overline\Fsurf) \big]$ denote the endpoint of $\mathcal{O}$. 
\end{definition}

The next definition specifies $\mathfrak{B}$, the set of $K$ boxes the edges of which will be closed during surgery.

\begin{definition}
 Let $\big( \mathcal{B}_1, \ldots, \mathcal{B}_\const \big)$ be the list of $\const$ distinct $K$-boxes with the property that the conditional probability under $P \big( \cdot \big\vert NT(n) \big)$ of
 $\text{INS}_K^{\overline{H}_{M_n}}(\overline\Fsurf) \subseteq \big\{ \mathcal{B}_1,\ldots,\mathcal{B}_\const \big\}$ 
 is maximal over such lists. We abbreviate $\mathfrak{B} = \big\{ \mathcal{B}_1, \ldots, \mathcal{B}_\const \big\}$.
 
 Let $\xlower \in \Z^d$ have the property that  the conditional probability under
 $P \big( \cdot \big\vert NT(n) \cap \big\{ \text{INS}_K^{\overline{H}_{M_n}}(\overline\Fsurf) \subseteq \mathfrak{B} \big\} \big)$ of $\bigxlower = x$ is maximal over  $x \in \Z^d$ for the choice $x = \xlower$. 
\end{definition}
Note that, since on $NT(n)$ we know that $ \text{INS}_K^{\overline{H}_{M_n}}(\overline\Fsurf)$ has at most $Q$ points that are located in $B_{\Z^d}(0,C_1n)$, we have
\begin{equation}\label{eqntone}
P \Big( \text{INS}_K^{\overline{H}_{M_n}}(\overline\Fsurf) \subseteq \mathfrak{B} \, \Big\vert \, NT(n)   \Big) 
 \geq \abs{B_{\Z^d}(0,C_1n)}^{-Q} \geq  cn^{-C\const},
\end{equation}
and
\begin{equation}\label{eqnttwo}
P \Big(  \bigxlower = \xlower  \, \Big\vert \, \text{INS}_K^{\overline{H}_{M_n}}(\overline\Fsurf) \subseteq \mathfrak{B} , NT(n)   \Big) 
 \geq cK^{-d} \const^{-1}.
\end{equation}
For the latter, we used that, under the conditional measure, $\bigxlower$ necessarily lies in one of the boxes comprising $\mathfrak{B}$. 

\begin{definition} \label{resample}
Let $A \subseteq \Z^d$. Let $\sigma_A$ be the random map that, applied to 
$\omega \in \{ 0,1 \}^{E(\Z^d)}$, returns
\[
\sigma_A(\omega) = \begin{cases}  \omega_e' & \textrm{if $e \in E(\Z^d)$, with both endpoints of $e$ in $A$}, \\ 
  \omega_e & \textrm{otherwise,}  \end{cases}
 \]
where $\big\{ \omega_e': e \in E(\Z^d) \big\}$ is an independent collection of independent Bern$(p)$ random variables.
\end{definition}
We now define the surgical operation $\chi$.
\begin{definition}
Let $\chi$ denote the random procedure formed by successively and independently applying 
$\sigma_{\mathcal{B}_1}, \ldots, \sigma_{\mathcal{B}_\const}$.
 Note that $\chi$ leaves invariant the law $P_p$.
\end{definition}
The role of $\chi$ is to take a "satisfactory" input configuration -- in which a reasonably typical trap of depth $n$ is present -- and to randomly alter it in such a manner that there is a reasonable prospect that this surgery acts "successfully", producing an output configuration in which a sausage connection of length close to $n$ appears.  
\begin{definition}
We say that the input $\omega \in \{0,1\}^{E(\Z^d)}$ is satisfactory if it satisfies
$NT(n)  \cap  \big\{ \text{INS}_K^{\overline{H}_{M_n}}(\overline\Fsurf) \subseteq \mathfrak{B} \big\}  \cap \big\{ \bigxlower = \xlower \big\}$.

We say that the action of $\chi$ is successful if, for each $1 \leq i \leq \const$, each edge whose status is resampled by $\sigma_{\mathcal{B}_i}$ is closed.
\end{definition}

Consider the random operation $\chi$ applied to the law $P_p$. The probability of successful action is bounded away from zero uniformly in $n$. Furthermore, we claim that successful action on satisfactory input yields an outcome such that $B_{\Z^d}(0,C_1n)$ contains
a finite open cluster containing two points $a,b$ such that $a\cdot \vec{\ell}> b\cdot \vec{\ell} + (1-\epsilon)n$. 
Indeed, we may choose $a = 0$ and $b=\xlower$, since these two points are connected by $\mathcal{O}$. It remains to argue that the cluster containing $0$ is finite. To see this, note that, in the configuration before surgery, any simple infinite open nearest-neighbor path passes through $\big[ \text{INS}_K^{\overline{H}_{M_n}}(\overline\Fsurf) \big]$ by Lemma~\ref{mnbvc}, which is a subset of $\mathfrak{B}$ in the case that input is satisfactory. Successful action results in a configuration where every edge in $\mathfrak{B}$ is closed, so that the open cluster containing $0$ is finite after resampling. This procedure of resampling is what we refer to as surgery.

Set $V_n$ equal to the event that there exists a finite open cluster in $B_{\Z^d}(0,C_1n)$ containing two points such that $a\cdot \vec{\ell}> b\cdot \vec{\ell} + (1-\epsilon)n$. Hence,
\[
cP_p\Big( NT(n),     \text{INS}_K^{\overline{H}_{M_n}}(\overline\Fsurf) \subseteq \mathfrak{B}  , \backtrack(0) >n,  \bigxlower = \xlower \Big) \leq  P_p \big( V_n \big).
\]

By Lemma \ref{we_have_M}, (\ref{eqntone}) and~(\ref{eqnttwo}), we find that
$$
 P_p \big( V_n \big) \geq   \tfrac{1}{3\epsilon n}  (Cn)^{-\const} K^{-d} \const^{-1} P_p\big( \backtrack(0) > n \big),
$$
or, more simply, 
\begin{equation}
\label{prob_vn}
P_p \big( V_n \big) \geq   n^{-C} P_p\big( \backtrack(0) > n \big).
\end{equation}

We have managed to turn a trap into a long finite open cluster in the direction $\vec{\ell}$ at relatively low probabilistic cost. We now wish to turn this cluster into a sausage-connection at low probabilistic cost.

\noindent{\bf Proof of Proposition~\ref{expo_3D}.}
Given $V_n$, we can perform a finite surgery, depicted in Figure~16, at some given points of maximal and minimal scalar product with $\vec{\ell}$ in the finite cluster of $0$ (those two points are in $B_{\Z^d}(0,C_1n)$) to turn that cluster into a sausage-connection between two points: 
\begin{eqnarray}\label{prob_vn1}
P_p \big( V_n \big)   & \leq &  C n^C P_p \Big( \textrm{there exist sausage-connection between some pair}    \nonumber  \\
        &  & \qquad  \textrm{of points  $a,b \in B_{\Z^d}(0,C_1n)$ satisfying $(a-b)\cdot\vec{\ell}\geq(1-\epsilon)n$} \Big). \nonumber
\end{eqnarray}

\begin{figure}[h]
\centering
\epsfig{file=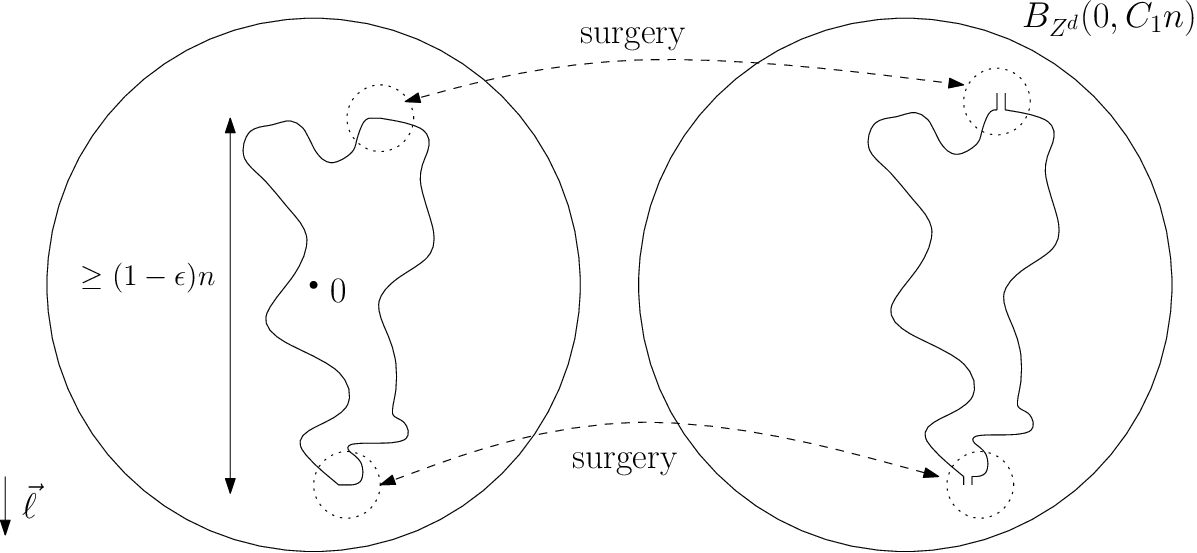, width=13cm}
\caption{A surgery procedure to turn a finite cluster into a sausage-connection}
\end{figure}

By $a,b \in B_{\Z^d}(0,C_1n)$,  translation invariance implies that
\begin{eqnarray}
 & & P_p \Big( 0  \textrm{ is sausage-connected to some $\mathcal{H}^-(-m)$ for $m\geq (1-\epsilon)n$} \Big) \nonumber  \\
 & \geq & n^{-C}  P_p\big( \backtrack(0) > n \big). \nonumber
\end{eqnarray}

By Lemma~\ref{expoa}, this means that for any $\epsilon>0$ we have
\[
{\bf P}_p[\backtrack(0)> n]  \leq C n^C \exp \big\{ -(1-\epsilon)\hatxi n \big\},
\]
which means that 
\[
\lim  n^{-1} \log {\bf P}_p[\backtrack(0)>n] \leq -(1-\epsilon) \hatxi, 
\]
and by taking $\epsilon$ to $0$ 
\begin{equation}
\label{limsupbk}
\lim  n^{-1} \log {\bf P}_p[\backtrack(0)> n]  \leq -\hatxi,
\end{equation}
this proves the upper bound on ${\bf P}_p[\backtrack(0)> n]$.

Take $x\in \partial\mathcal{H}^-(-n)$, we introduce  the event 
$$
\big\{0\text{ is only sausage-connected to } x\big\}
$$ 
when we have that $\{0\text{ is sausage-connected to } x\}$ and that $0$ is not sausage-connected to any $y\in \partial \mathcal{H}^-(-n)$ with $y\neq x$. We note that the events $\{0\text{ is only sausage-connected to } x\}$, for $x\in \partial \mathcal{H}^-(-n)$, form a partition of the event $\{0\text{ is only sausage-connected to } \mathcal{H}^-(-n)\}$.

We now show the lower bound. Starting from the event 
$$
\big\{ 0\text{ is only sausage-connected to } x \big\},
$$
using a finite surgery similar to that which shows (\ref{eqsurgeryex}) at the bottom of the sausage-connection, and imposing that the event $\mathcal{I}_x^-$ (which has positive probability) happens at the top $x$ of the sausage-connection, we ensure the occurrence of $\{\mathcal{B}\mathcal{K}(0) >n\}$ and $\mathcal{I}$. This yields, for any $x\in \partial\mathcal{H}^-(-n)$, that
\begin{eqnarray*}
& & P_p[0\text{ is only sausage-connected to } x] \\
& \leq & C P_p[0\text{ is only sausage-connected to } x, \mathcal{I},\mathcal{B}\mathcal{K}(0) >n],
\end{eqnarray*}
and hence, using that on a given $\omega$ $0$ is connected to a finite (non random and non $n$-dependent) number of points in $\mathcal{L}(-n)$, 
\begin{eqnarray*}
& & P_p[0\text{ is sausage-connected to } \mathcal{H}^-(-n)] \\
& = & \sum_{x\in \partial \mathcal{H}^-(-n)} P_p[0\text{ is only sausage-connected to } x] \\
& \leq & C\sum_{x\in\partial\mathcal{H}^-(-n)} P_p[0\text{ is only sausage-connected to } x, \mathcal{I},\mathcal{B}\mathcal{K}(0) >n] \\
& \leq & CP_p[\mathcal{I},\mathcal{B}\mathcal{K}(0) >n] \\
& \leq & {\bf P}_p[\mathcal{B}\mathcal{K}(0) >n],
\end{eqnarray*}
so that, by Lemma~\ref{expoa},
\[
\liminf \frac{\log {\bf P}_p[\mathcal{B}\mathcal{K}(0) > n]}{n} \geq -\hatxi.
\]

Alongside~(\ref{limsupbk}), we obtain
\[
\lim \frac{\log{\bf P}_p[\mathcal{B}\mathcal{K}(0)> n]}{ n} = -\hatxi,
\] 
which proves Lemma~\ref{expo_3D}. \qed \medskip

\noindent{\bf Proof of Lemma~\ref{estim_height} for $\Z^d$ with $d \geq 3$.}
Recalling~(\ref{prob_vn}) and~(\ref{prob_vn1}), 
\begin{eqnarray}
 & & P_p \big( \exists a,b \in B_{\Z^d}(0,C_1n): \, \textrm{$a,b$ are sausage-connected},  (a-b)\cdot\vec{\ell}\geq(1-\epsilon)n \big) \nonumber \\
 & \geq &  n^{-C} P_p\big( \backtrack(0) > n \big), \nonumber
\end{eqnarray}
for every $\epsilon>0$. By translation invariance, we see that
\[
P_p \big( \text{$0$ is sausage-connected to $x$} \big) 
\geq   n^{-C} P_p\big( \backtrack(0) > n \big),
\]
where $x\in B_{\Z^d}(0,C_1n)$  with $x\cdot\vec{\ell}\geq (1-\epsilon)n$ is chosen to maximize the left-hand side of this inequality.

 Using a finite surgery at the bottom $x$ of the sausage-connection, we can form a one-headed trap at $0$; (the notion of one-headed trap is introduced at the start of Section $3$). Hence, we see that
 \[
P_p \big(\mathcal{T}^1(0) \neq \emptyset \text{ and } \mathfrak{D}(e_1)\geq (1-\epsilon) n  \big) \geq n^{-C} P_p\big( \backtrack(0) > n \big).
\]

Using Lemma~\ref{expo_3D} and letting $\epsilon$ go to $0$ yields Lemma~\ref{estim_height} for $\Z^d$ with $d \geq 3$. \qed 
\end{document}